\documentclass[12pt]{article}
\usepackage{amssymb}
\usepackage{hyperref}
\usepackage{graphicx}
\usepackage{caption}
\usepackage{graphicx,latexsym,amssymb,amsmath}
\usepackage{tikz}
\definecolor{myblack}{RGB}{173, 216, 230}
\DeclareMathOperator*{\argmin}{arg\,min}
\usepackage[french,english]{babel}
\usepackage{here} 
\newcommand{\blind}{1}
\usepackage{bbm}
\usepackage{multirow}
\usepackage{amsmath}
\usepackage{amssymb}
\usepackage{amsfonts}
\newtheorem{assumptionalt}{Assumption}
\newenvironment{assumptionp}[1]{
  \renewcommand\theassumptionalt{#1}
  \assumptionalt
}{\endassumptionalt}

\def\R{{\mathbb{R}}}

\def\N{{\mathbb{N}}}

\def\P{{\mathbb{P}}}

\newcommand{{\convp}}{{\buildrel p\over\longrightarrow}}

\newcommand{{\Vs}}{{\cal V}}
\newcommand{{\Ps}}{{\cal P}}
\newcommand{{\Ss}}{{\cal S}}
\newcommand{{\Xs}}{{\cal X}}
\newcommand{{\Ls}}{{\cal L}}
\newcommand{{\Ns}}{{\cal N}}
\newcommand{{\Zs}}{{\cal Z}}
\newcommand{{\Fs}}{{\cal F}}

\usepackage{pdflscape}
\usepackage{afterpage}
\usepackage{capt-of}

\newtheorem{Lemma}{Lemma}[section]
\newtheorem{Assumption}{Assumption}[section]
\newtheorem{Theorem}{Theorem}[section] 
\newtheorem{Definition}{Definition}[section] 
\newcommand{\proofend}{$\hfill\Box{~}$}
\newenvironment{Proof}{\noindent {\em{\bf Proof.}}}{\proofend\\}

\newtheorem{corollary}{Corollary}[section]


\usepackage{natbib}
\usepackage{enumerate}
\def\spacingset#1{\renewcommand{\baselinestretch}%
{#1}\small\normalsize} \spacingset{1}

\usepackage{color}
\newcommand{\tr}{\textcolor{black}}
\newcommand{\red}{\textcolor{black}}
\usepackage{adjustbox}

\begin{document} 
\if1\blind
{
 \title{\bf Instrumental variable estimation of dynamic treatment effects on a duration outcome}
 \author{ Jad Beyhum \thanks{
   Financial support from the European Research Council (2016-2021, Horizon 2020 / ERC grant agreement No.\ 694409) is gratefully acknowledged.}\\
  ORSTAT, KU Leuven and CREST, ENSAI \vspace*{5pt}\\ 
   Samuele Centorrino \thanks{{\color{black}Samuele Centorrino would like to thank Stony Brook Research Computing and Cyberinfrastructure, and the Institute for Advanced Computational Science at Stony Brook University for access to the high-performance SeaWulf computing system, which was made possible by an NSF grant (\#1531492).}}\\
  Economics Department, Stony Brook University \vspace*{5pt}\\
  Jean-Pierre Florens \thanks{Jean-Pierre Florens acknowledges funding from the French National Research Agency (ANR) under the Investments for the Future program (Investissements d'Avenir, grant ANR-17-EURE-0010).}\\
  Toulouse School of Economics, Universit\'e Toulouse 1 Capitole\\
  and\\
   Ingrid Van Keilegom $^*$\\
  ORSTAT, KU Leuven}
 \maketitle
} \fi

\if0\blind
{
 \bigskip
 \bigskip
 \bigskip
 \begin{center}
  {\LARGE\bf Instrumental variable estimation of dynamic treatment effects on a survival outcome}
\end{center}
 \medskip
} \fi
\bigskip

\begin{abstract}
This paper considers identification and estimation of the time $Z$ until a subject is treated on a duration $T$. The treatment is not randomly assigned, $T$ is randomly right censored by a random variable $C$, and the time to treatment $Z$ is right censored by $\min(T,C)$. The endogeneity issue is treated using an instrumental variable explaining $Z$ and independent of the error term of the model. We study identification in a fully nonparametric framework. We show that our specification generates a system of integral equations, of which the regression function of interest is a solution. We provide identification conditions that rely on this identification equation. We assume that the regression function follows a parametric model for estimation purposes. We propose an estimation procedure and give conditions under which the estimator is asymptotically normal. The estimators exhibit good finite sample properties in simulations. Our methodology is applied to find evidence supporting the efficacy of a therapy for burnout.
\end{abstract}

 \noindent%
{\it Keywords:} Dynamic treatment; Endogeneity; Instrumental variable; Censoring; Nonseparability.  \\
\vfill

\newpage
\spacingset{1.4}

\setcounter{footnote}{0}
\setcounter{equation}{0}

\section{Introduction} 

Consider subjects in a certain state, such as unemployment or medical leave. The policymaker is sometimes interested in how the timing of some treatment affects the time the subjects remain in the state. An example comes from labor economics, where we may want to evaluate the effect of reducing unemployment benefits after some time on the unemployment duration. Also, in health economics, one may seek to find the optimal time to give therapy for workers on medical leave because of burnout. This setting corresponds to our empirical application.

{\color{black} The present paper considers the problem of estimating the causal effect of a treatment. The treatment is dynamic because it can be given at any point in time, and its effect can vary depending on its timing. Moreover,} the time $Z$ until the treatment is given is not randomly assigned. We solve the endogeneity issue thanks to an instrumental variable $W$ independent of the error term of the model and sufficiently related to $Z$. The outcome duration of interest $T$ is randomly right censored by a censoring time $C${\color{black}, assumed independent of the other variables of the model.} The timing of the treatment $Z$ is itself censored by $\min(T,C)$. {\color{black} The censoring of $Z$ by $T$ is endogenous, since the latter variables are dependent.} This context corresponds to studies where follow-up stops when the subjects leave the state of interest, or the treatment cannot be given to participants who leave the state of inflow. In the latter case, $Z$ corresponds to some latent duration to treatment, which is realized only when $Z\le T$. In the labor economics example mentioned above, this is justified because unemployment benefits cannot be reduced if the subject is no longer unemployed. In the illustration from health economics, there is no reason to treat cured workers.


The contributions of the paper are as follows. \red{We study a dynamic duration model where we are interested in the hazard rate of $T(z)$, where $T(z)$ is the potential outcome of the duration when the treatment happens at time $z$. We make a number of assumptions on this duration model. One assumption, called no anticipation, allows us to handle the issue that we never observe treatment times $Z$ which are larger than $T$ since $Z$ is right censored by $\min(T,C)$. A second assumption is a rank invariance condition common in the nonparametric instrumental variable (NPIV) literature. This assumption allows us to identify and estimate the hazard rate of $T(z)$ using an instrumental variable $W$ independent of the error term of the model. Specifically, we rewrite the duration model as a nonseparable NPIV model and adapt tools from the NPIV literature to study the identification of the regression function of the rewritten model in a fully nonparametric framework. The hazard rates of the potential outcomes are functionals of the regression function of the nonseparable model. Our identification results are not straightforward applications of the results from the nonseparable NPIV model literature since our nonseparable model is different from their standard quantile model because of the dynamic nature of the problem. For estimation purposes, we assume that the regression function follows a parametric model. We propose an estimation procedure and give conditions under which the semiparametric estimator is asymptotically normal. The finite sample properties of the estimator are assessed through simulations. We apply our methodology to evaluate the effect of the timing of therapy on the duration of medical leave for burnout.}


There exists an extensive literature on instrumental variable methods with randomly right-censored duration outcomes, where both the treatment and the instrument are time-independent. We only cite here some references to avoid lengthening the paper. Many works study semiparametric models, see, e.g., \cite{tchetgen2015instrumental} for an additive hazard model, \cite{chernozhukov2015quantile} for a quantile regression model, or \cite{martinussen2019instrumental} for the Cox model. Other works 
(\cite{frandsen2015treatment}, \cite{sant2016program}, \cite{richardson2017nonparametric}, \cite{blanco2019bounds}, \cite{sant2021nonparametric}) provide nonparametric estimation results for the average treatment effects on the compliers (see \cite{angrist1996identification}) with binary treatment and instrument. Moreover, \cite{2021nonparametric} nonparametrically estimate the average quantile treatment effect over the whole population when both the treatment and the instrument are categorical. Also, \cite{centorrinoflorens2021} study nonparametric estimation when the treatment is continuous, and the model is additive. They allow the instrument to be categorical, even in this case. Our paper differs from this literature because it allows the treatment to be dynamic. Our dynamic setting is more complex than those previously analyzed since the treatment $Z$ is not always observed and endogenously censored by $\min(T,C)$. Among the papers above, \cite{2021nonparametric} is the most closely related to our paper since it relies on a nonseparable NPIV model as in \citet{CH}. \red{There are two main differences with this paper. First, here $Z$ is dynamic, which creates endogenous censoring of $Z$ by $T$. Second, in the present context, $Z$ is a continuous variable. This makes the inverse problem ill-posed and requires a different identification analysis and a new estimation procedure.}

{\color{black}Some papers study the evaluation of the effect of the time to treatment on a survival outcome when the treatment is endogenous without relying on an instrumental variable. The first set of papers assumes conditional unconfoundedness of the treatment. See \cite{sianesi2004evaluation}, \cite{AVdB2,lechner2009sequential}, \cite{lechner2010identification}, \cite{vikstrom2017dynamic} and \cite{van2020policy}, \cite{kastoryano2022dynamic} for the literature in econometrics. \cite{hernan2010causal} give an excellent overview of the methods available in the biostatistics literature to infer causal effects in this context. These studies do not rely on an instrumental variable to solve the endogeneity problem. Another approach is that of \cite{AVdB}, where it is assumed that both the timing of the treatment and the duration outcome follow a mixed proportional hazards model.

Finally, like us, some papers use instrumental variables to estimate the effect of an endogenous time-varying treatment on a right-censored outcome. \cite{bijwaard2005correcting} considers a setting where the instrument is binary, control group subjects are never treated and a mixed proportional hazard model is assumed. \cite{heckman2007} study a dynamic treatment effect model for outcomes that can be durations. However, they do not allow for censoring of $T$ by $C$ and of $Z$ by $T$, and they rely on a control function approach requiring additional structural assumptions, and no estimation theory is proposed. Another paper is \cite{VdBBM}, which studies a one-sided noncompliance setting, where the time to treatment is either equal to the instrument or $\infty$ (untreated). {\color{black} They identify local average treatment effects on a subset of the population analogous to the compliers in static problems (see \cite{angrist1996identification}). In contrast, the present work does not assume one-sided noncompliance, and our method allows us to estimate effects over the whole population. \cite{VdBBM} also applies their estimation strategy to evaluate the causal effect of a reform of the unemployment insurance system in France. Our approach could also be applied to their dataset. Under our assumptions, we would identify the effect of the reform on the whole population rather than only part of it. Finally, note a recent line of work in biostatistics (\cite{tchetgen2018marginal, cui2017instrumental, michael2020instrumental}) which studies the case where the instrument varies with time, ruling out the case with a static instrument studied in our paper. They have a local (time by time) estimation approach leveraging the time variation of the instrument. This estimation strategy can not be straightforwardly adapted to allow for a time-independent instrument.}}

{\color{black}\textbf{Outline.} }This paper is organized as follows. The model specification is given in Section 2. We study identification in Section 3. Section 4 is devoted to estimation and inference. Section 5 describes our simulations and the empirical application. Concluding remarks are given in Section 6. The proofs of the results and additional numerical experiments are in the online appendix. 

\section{The model}
\subsection{The duration model}
\red{Let $T(z)$ be the potential outcome of the duration when the treatment time is set to $z\in\bar{\R}_+=\R_+\cup\{\infty\}$. When $z=\infty$, $T(z)=T(\infty)$ corresponds to the duration that would have been realized if the subject of interest were never treated. The treatment time is a random variable $Z$, with support $\mathcal{Z}\subset\R_+$. We impose the consistency condition $T=T(Z)$. In our empirical application, studied in Section \ref{empirical_application}, we want to evaluate the effect of a therapy for burnout on the duration of medical leave. We possess a dataset where each observation corresponds to a worker on medical leave for burnout. The variable $Z$ is the time until a therapy for burnout starts, $T(z)$ is the duration of medical leave that would have been realized if the therapy for burnout had been started at time $z$ and $T$ is the actual duration of medical leave. }

\red{We assume that $T(z)$ is a continuous random variable and let $\lambda(z,\cdot):\R_+\mapsto\R_+$ be its hazard rate, that is 
$$\lambda(z,t)= \lim_{dt \to 0}\frac{\P(T(z)\in[t,t+dt]|T(z)\ge t)}{dt}.$$
We call $\lambda$ the ``structural hazard". Remark that $\lambda(z,\cdot)$ differs from the hazard rate of $T$ conditional on  $Z=z$ in general (this is the endogeneity issue). The main goal of the paper is to identify and estimate this structural hazard. }

\begin{figure}[H]
 \centering
 \includegraphics[width=80mm]{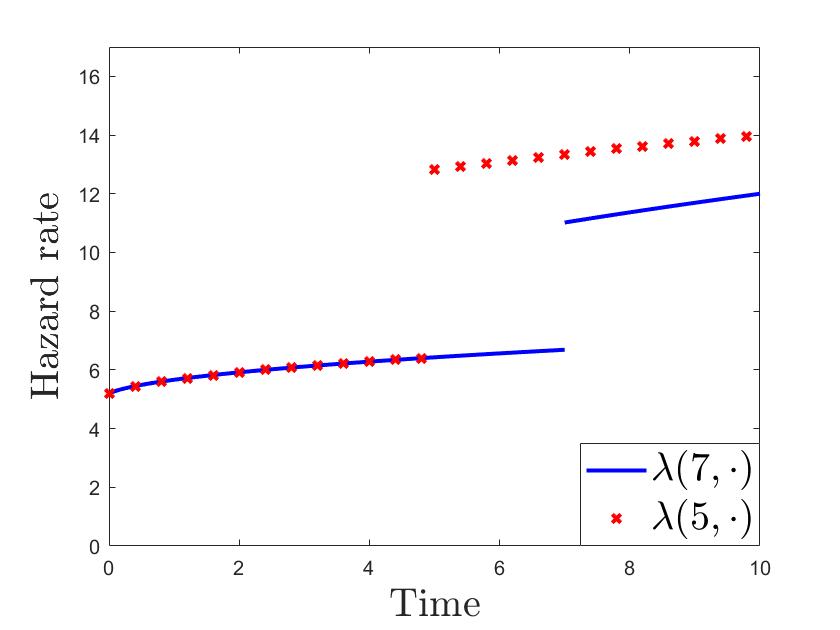}
 
 \caption{\red{Structural hazard rates under treatment levels $z=5$ and $z=7$ obtained by setting $\lambda(5,t)=( 5+0.2\sqrt{t})I(t<5)+(10+0.4\sqrt{t})I(5\le t)$ and $\lambda(7,t)=( 5+0.2\sqrt{t})I(t<7)+(6+0.6\sqrt{t})I(7\le t)$.}}
  \label{fig:illustration} 
\end{figure}

\red{As mentioned in the introduction, one of the problems considered in this paper is that $Z$ is censored by $T$. When $Z>T$, we never know to which potential outcome $T$ corresponds. To circumvent this issue, we impose the so-called no anticipation assumption (see \cite{AVdB, VdBBM, van2020policy}), which is standard in the literature on dynamic treatment effects. The no anticipation assumption can be formally stated as follows. 
\begin{Assumption}\label{haz}
For all $z,z'\in\bar{\R}_+$ and $0\le t\le \min(z,z')$, we have $\lambda(z,t)=\lambda(z',t)$.
\end{Assumption}
This assumption means that the hazard rate $\lambda(z,t)$ of $T(z)$ does not depend on $z$ when $z>t$. Under Assumption \ref{haz}, all observations for which $Z>T$ correspond to the same structural hazard (equal to $\lambda(\infty,\cdot)$), which allows us to solve the aforementioned issue that we do not observe $Z$ when $Z>T$. Denote by $\Lambda(z,t)=\int_0^t \lambda(z,s)ds$ the structural cumulative hazard of $T(z)$. By integration and derivation, Assumption \ref{haz} is also equivalent to $\Lambda(z,t)=\Lambda(z',t)$ when $t<\min(z,z')$, which is the way the no anticipation assumption is stated in \cite{AVdB} (except that the strict inequality $t<\min(z,z')$ is replaced by the weak inequality $t\le\min(z,z')$ in \cite{AVdB}). We illustrate the no anticipation assumption in Figure \ref{fig:illustration}, where we draw the structural hazard rates under treatment levels $z=5$ and $z=7$ (in the caption, we use $I(\cdot)$ to denote the indicator function). In this illustration, the treatment increases the hazard rate, but our model also allows the treatment to reduce the hazard rate. Before the treatment, the hazard rate does not depend on the treatment time, while after the treatment it does. Finally, note that the no anticipation assumption does not rule out that subjects actually anticipate receiving future treatment. It rather means that the counterfactuals that we are interested in correspond to a setting where the intervention does not change the anticipations but only the actual value of the treatment (see the discussion in \cite{AVdB} for more details).}

\red{To be able to identify treatment effects over the full population under endogeneity, we impose constraints on the unobserved heterogeneity of the model. For $z\in\bar{\R}_+$ and $t\in\R_{+}$, recall that $\Lambda(z,t)=\int_0^t \lambda(z,s)ds$ is the structural cumulative hazard under treatment $z$ and let $U(z)=\Lambda(z,T(z))$ be the hazard of $T(z)$ evaluated at $T(z)$. The $\{U(z)\}_{z\in\bar{\R}_+}$ can be thought of as the unobserved heterogeneity of the model. In our empirical application, it could correspond to the underlying mental health of the subject. We impose the following conditions on them.
\begin{Assumption} \label{RI}The following holds: 
\begin{itemize}
\item[(i)] There exists a random variable $U$ such that $U(z)=U,\ \text{for all $z\in\bar{\R}_+$;}$
\item[(ii)] For all $z\in\bar{\R}_+$, $\lambda(z,\cdot)$ is continuous on $[0,z)$ and $[z,\infty)$.
\end{itemize}
\end{Assumption}
Condition (ii) allows the hazard rate $\lambda(z,\cdot)$ to be discontinuous at the time of treatment $z$. This permits behaviors similar to that of Figure \ref{fig:illustration}, where treatment makes the hazard rate ``jump". We show in the online appendix (Lemma \ref{ls11}) that Assumption \ref{RI} (ii) implies that the cumulative hazard $\Lambda(z,\cdot)$ maps the support of $T(z)$ to $\R_+$ and is strictly increasing on the support of $T(z)$. Therefore we can define $\Lambda(z,\cdot)^{-1}$, the inverse of the mapping $\Lambda(z,\cdot)$ restricted to the support of $T(z)$ (that is $\Lambda(z,\cdot)^{-1}(u)$ is the unique element $t$ of the support of $T(z)$ such that $\Lambda(z,t)=u$). This implies that $T(z)=\Lambda(z,\cdot)^{-1}(U(z))$ is strictly increasing in $U(z)$. Then, Assumption \ref{RI} (i) yields that, for two subjects $i$ and $j$, $T_i(z)>T_j(z)$ implies $U_i(z) = U_i>U_j = U_j(z)$, which leads to $T_i(z')>T_j(z')$, for all $z,z'\in\bar{\R}_+$. In other words, Assumption \ref{RI} (i) implies that the rank in the outcome of any two subjects is the same across all potential outcomes. Assumption \ref{RI} (i) is therefore a rank invariance assumption as in \cite{CH}. This assumption restricts the heterogeneity of the treatment effects on the duration: the treatment can change the quantiles of the distribution of the potential outcomes but it cannot change the rank that a subject has in this distribution.  Moreover, remark that the rank invariance assumption does not restrict the possible values of the structural hazard $\lambda(z,t)$ beyond the continuity condition in Assumption \ref{RI} (ii) and therefore, in this sense, rank invariance is not a constraint on the marginal distribution of $T(z)$. It only imposes limits on the joint distribution of potential outcomes, that is the distribution of $(T(z))_{z\in \mathcal{Z}}$. Note also that we could relax the rank invariance assumption into a rank similarity assumption as in \cite{CH} while keeping all results valid.}

\red{Another important fact is that Assumption \ref{RI} (ii) implies that for all $z\in\bar{\R}_+$, $U(z) \sim\text{Exp}(1)$, where $\text{Exp}(1)$ is the unit exponential distribution (see Lemma \ref{ls11}). In virtue of Assumption \ref{RI} (i), we can therefore write $\Lambda(Z,T)=U,$ with $U\sim \text{Exp}(1)$. }

\red{The variables $U$ and $Z$ may be dependent, which creates an endogeneity issue, that is $\lambda(z,\cdot)$ may differ from the hazard rate of $T$ given $Z=z$. In the context of the empirical application, the endogeneity may, for instance, be due to the fact that subjects with worse burnout are more likely to be treated early. There exists an instrument $W$ allowing to solve this issue. In the real data, the instrument is related to the medical center to which workers on medical leave are assigned (see Section \ref{empirical_application} for more details). The support of $W$ is denoted $\mathcal{W}$. For simplicity, we limit ourselves to the case where $W$ is scalar, that is $\mathcal{W}\subset \R$. \red{We impose the following assumption:
\begin{Assumption}\label{IVind}
$W$ is independent of $U$.
\end{Assumption}}
The duration $T$ is randomly right censored by a random variable $C$ with support in $\bar{\R}_+$, so that we do not observe $T$ but $Y=\min(T,C)$. In the application, $C$ corresponds to the duration during which subjects on medical leave are followed in the data (the follow-up stops after two years of medical leave or at the end of 2020). The observables are $(Y, \delta, \tilde Z, \tilde D, W)$, where $\delta =I(T\le C)$, $\tilde Z=\min(Z,Y)$, is a censored version of $Z$ and $\tilde D= I(Z\le Y)$ is a censored version of $D=I(Z\le T)$, the treatment indicator. Note that $D=1$ for treated observations only (in the sense that they receive treatment before the end of their spell). In the burnout data, we have $\delta =1$ if the duration of medical leave is observed and $0$ otherwise, $\tilde Z$ is the minimum between the treatment time, the censoring time and the duration of medical leave, $D$ is equal to $1$ if the subject is treated before the end of its medical leave and $0$ otherwise and $\widetilde{D}=1$ if we observe the treatment time in the dataset and $0$ otherwise.}
\subsection{Reformulation as a nonseparable NPIV model}
\red{As mentioned in the introduction, this paper makes use of tools from the NPIV model to solve the dynamic problem. Recall that for all $z\in\bar{\R}_+$, $\Lambda(z,\cdot)^{-1}(u)$ is the unique element $t$ of the support of $T(z)$ such that $\Lambda(z,t)=u$ and that $I(\cdot)$ is the indicator function. The reformulation as a nonseparable NPIV model is based on the following lemma. }
\begin{Lemma}\label{reformulate} Under Assumptions \ref{haz} and \ref{RI}, we can write
\begin{equation}\label{model}T= \varphi(Z,U)=I(Z>\varphi_0(U)) \varphi_0(U) +I(Z\le \varphi_0(U))\varphi_1(Z,U)\ \text{a.s.},\end{equation} 
where $\varphi_0:\R_+\mapsto \R_+$ is equal to $\Lambda(\infty,\cdot)^{-1}$ and $\varphi_1(z,\cdot): \R_+\mapsto \R_+$ is equal to $\Lambda(z,\cdot)^{-1}$. Moreover, for all $z\in\R_+$, we have $\varphi_1(z,\varphi_0^{-1}(z))=z$.\end{Lemma}
Equation \eqref{model} means that there exists mappings $\varphi_0$ and $\varphi_1$ such that $T$ is equal to $\varphi_0(U)$ when $Z<\varphi_0(U)$ and equal to $\varphi_1(Z,U)$ otherwise. Intuitively, this signifies that there are two regressions functions: one of the ``not yet treated" corresponding to $\varphi_0$, and one for the ``already treated" which is $\varphi_1$. Writing the NPIV model as functions of these $\varphi_0$ and $\varphi_1$ allows us to define a parameter space which is a vector space (see Section \ref{sec.idwcens}).  Note that, by definition, $\varphi_0$ and $\varphi_1(z,\cdot)$ are strictly increasing. Moreover, the fact  that $\varphi_1(z,\varphi_0^{-1}(z))=z$ for all $z\in\R_+$ (see the end of Lemma \ref{reformulate}) implies that $\varphi$ is strictly increasing too.

\red{Equation \eqref{model} defines a nonseparable NPIV model similar to that of \cite{CH}. Hence, we can use tools from the literature on nonseparable NPIV models to obtain identification results on $\varphi$. Throughout the paper, we will therefore identify and estimate $\varphi$ rather than $\lambda$ because this approach simplifies the mathematical analysis. Using Lemma \ref{reformulate} and Assumption \ref{haz}, one can show that $\varphi(z,\cdot)=\Lambda(z,\cdot)^{-1}$ (see the proof of Lemma \ref{reformulate} in the online appendix for more details). The structural hazard $\lambda(z,t)$ is therefore the derivative of the inverse of $\varphi(z,\cdot)$ at $t$. As a result, $\varphi$ is a one-to-one transformation of $\lambda$ and identification of $\varphi$ implies identification of $\lambda$ and therefore of many quantities of interest in duration models (such as the survival function or the cumulative hazard). Note that, the present model involves two additional complexities with respect to the original model of \cite{CH}. First, $T$ is right censored by $C$. Second, $Z$ is right censored by $\min(T,C)$ and to (partially) solve this problem we have imposed a no anticipation assumption which implies that $\varphi$ follows the specific functional form of the right-hand side of \eqref{model}, which is different from the one used in the standard quantile model of \cite{CH}. As we will see later, this implies that we have to adapt results of the nonseparable NPIV model to our nonstandard model.}


\section{Identification} 

We study identification in a fully nonparametric setting. First, we show that our model generates a system of integral equations. Then we derive identification results based on this system of equations.

\subsection{Identification equation} In order to formulate the identification equation, we introduce the following reduced form quantities. For $t\in\R_+ $, $w\in\mathcal{W}$, let 
\begin{align*}
F_0(t,w)=\P\left( T\le t,D=0,W\le w\right);\ F_W(w)=\P(W\le w).
\end{align*}
Moreover, for a mapping $\phi:\mathcal{Z}\mapsto \R_+$ and $w\in\mathcal{W}$, we define
$$F_1(\phi,w)= \P(T\le \phi(Z), D=1,W\le w).$$
The following theorem states the system of equations that we use to obtain identification results.
\begin{Theorem}\label{th.main_eq} Let Assumptions \ref{haz}, \ref{RI} and \ref{IVind} hold. 
For the true $\varphi_0$ and $\varphi_1$ (defined below equation \eqref{model}) and all $w\in \mathcal{W}$, we have
\begin{equation}\label{main_eq}F_0(\varphi_0(u),w) +F_1(\varphi_1(\cdot,u),w)=(1-e^{-u})F_W(w).\end{equation}
\end{Theorem}
\red{In Section \ref{altchar} of the online Appendix, we provide an alternative characterization of the model in terms of reduced-form (conditional) hazard rates and survival functions of the observed duration $T$, which are natural quantities in the context of duration models. }

\subsection{Identification without censoring}\label{sec.idwcens}

In this section, we discuss the simple case where there is no censoring, that is $C=\infty$ a.s. Then, $Y=T$ and $\tilde D= D$, which implies that $F_0(t,w)$ is identified for all $t\in\R_+$ and $w\in\mathcal{W}$. Moreover, when $D=1$, we have $\tilde Z=Z$, and, hence, $F_1(\phi,w)$ is identified for all $\phi:\mathcal{Z}\mapsto \R_+$ and $w\in\mathcal{W}$. \red{Therefore, $F_0,F_1,F_W$ in \eqref{main_eq} are all identified in the absence of censoring.} Hence, uniqueness of the solutions to \eqref{main_eq} implies identification. This allows us to derive identification results in the next two subsubsections.

\red{We focus on identification of $\varphi(\cdot,u)$ for a given $u\in\R_+$.
At this point, it is useful to define the set to which $(\varphi_0(u),\varphi_1(\cdot,u))$ belongs. Let the parameter space $\mathcal{P}$ be the set of $(\psi_0,\psi_1)$ such that $\psi_0\in \R$, and $\psi_1$ is a bounded mapping from $\mathcal{Z}$ to $\R$. This set $\mathcal{P}$ is a vector space and we endow it with the norm $\|\cdot\|_{\mathcal{P}}$, where $\left\|(\psi_0,\psi_1)\right\|_{\mathcal{P}}^2= E[(\psi(Z))^2|U=u],$ with $\psi(z)=\psi_0I(z> \psi_0)+\psi_1(z)I(z\le\psi_0).$
The function $\psi$ is the mapping "induced" by $(\psi_0,\psi_1)$.
It is similar to the object that we want to identify ($\varphi(\cdot,u)$). The norm $\left\|(\psi_0,\psi_1)\right\|_{\mathcal{P}}$ is finite
because $\psi_1$ is bounded. We assume that $\psi_1$ is bounded, because, in practice, to identify $\varphi_1(\cdot,u)$ in the presence of censoring with finite support, we will have to assume that $\varphi_1(\cdot,u)$ is bounded by the upper bound of the support of the censoring variable (see the discussion in Section \ref{cens.sec.id} and Assumption \ref{paramid}(ii) later in the paper). The fact that we restrict the parameter space to bounded functions $\psi_1$ also allows us to weaken the completeness conditions for identification.}

\red{Remark also that since (i) $\mathcal{P}$ is not a standard $L^2$-space with respect to a continuous distribution and (ii) what we wish to identify is not $(\varphi_0(u),\varphi_1(\cdot,u))$ but $\varphi(\cdot,u)$, we cannot rely on the high-level identification theory for nonseparable NPIV models as in \cite{CH,chen2014local}. Therefore to obtain our identification results, we adapt the proofs of these papers to our specific model.}
\subsubsection{Local identification}
 \red{We start by local identification.
\begin{Definition}The regression function $\varphi(\cdot,u)$ is locally identified in a set $\mathcal{N}\subset \mathcal{P}$ if for all $(\psi_0,\psi_1)\in\mathcal{N}$, \begin{equation}\label{restatemaineq}F_0(\psi_0,w) +F_1(\psi_1,w)=(1-e^{-u})F_W(w),\text{for all}\ w\in\mathcal{W},\end{equation}
implies that $\psi(Z)=\varphi(Z,u)$ almost surely.
\end{Definition}}
\noindent \red{Notice that we only seek to identify $\varphi(z,u)$ for all $z\in\mathcal{Z}$ and not $(\varphi_0(u),\varphi_1(z,u)), z\in \mathcal{Z}$. This is because it is not possible (and not interesting) to identify $\varphi_1(z,u)$ for $z>\varphi_0(u)$ because $\varphi_1(z,u)$ never generates the data for such $z$.}
The main assumption for local identification is the following bounded completeness condition.
\begin{Assumption}\label{CC}For all $(\psi_0,\psi_1)\in \mathcal P$,
\begin{align*} &E[\psi_0I(Z>\varphi_0(u))+I(Z \le\varphi_0(u)) \psi_1(Z)|U=u,W]=0\ a.s. \\
&\Rightarrow \P\left(\psi_0I(Z>\varphi_0(u))+I(Z \le\varphi_0(u)) \psi_1(Z)=0|U=u\right)=1.\end{align*}
\end{Assumption}
Intuitively, this condition means that $Z$ and $W$ are sufficiently dependent given $U=u$. 
Assumption \ref{CC} is implied by the bounded completeness of $Z$ given $W,U=u$, that is for all {\color{black} bounded} functions \begin{equation}\label{complete1}\text{$m:\mathcal{Z}\mapsto \R$, $E[m(Z)|U=u,W]=0\ a.s. \Rightarrow \P(m(Z)=0|U=u)=1$.}\end{equation} Such a condition is imposed in \cite{cazals2016nonparametric}. In separable NPIV models, the related bounded completeness condition
\begin{equation}\label{complete2}\text{$m:\mathcal{Z}\times \R_+\mapsto \R$, $E[m(Z)|W]=0\ a.s. \Rightarrow m(Z)= 0\ a.s.$}\end{equation} for $m$ belonging to some class of functions, is often imposed \citep[see][among others]{newey2003,Darolles}. Condition \eqref{complete1} is a ``conditional on $U=u$" version of condition \eqref{complete2}. Several authors have provided sufficient conditions for bounded completeness as stated in \eqref{complete2} \citep[see][among others]{newey2003,xavier2011, hu2018nonparametric, andrews2011}. These sufficient conditions are restrictions on the family of distributions $\{Z|W=w\}_{W=w}$ where ``$Z|W=w$" stands for the distribution of $Z$ given $W=w$. To obtain sufficient conditions for \eqref{complete1}, it therefore suffices to impose the sufficient conditions from \cite{newey2003,xavier2011, hu2018nonparametric, andrews2011} on the family of distributions $\{Z|W=w,U=u\}_{W=w}$ rather than $\{Z|W=w\}_{W=w}$. As a last remark, notice that restricting conditions \eqref{complete1} to bounded functions makes it more likely to hold (see the aforementioned papers for further details).

In addition to Assumption \ref{CC}, we also impose some regularity conditions. Since these conditions are technically involved, they are stated in the online appendix (see Assumption \ref{asss1} in Section \ref{prooflocident}). Remark that these regularity conditions require Gâteaux differentiability of some operator but no Fréchet differentiability is needed.
\red{We have the following local identification result. 
\begin{Theorem}\label{th.loc_ident} Let Assumptions \ref{haz}, \ref{RI}, \ref{IVind}, \ref{CC}, \ref{asss1} hold and assume that $(\varphi_0(u),\varphi_1(\cdot,u))\in\mathcal{P}$. Then, for all $(\psi_0,\psi_1)\in\mathcal{P}$,
 there exists $\epsilon>0$ such that $\varphi(\cdot,u)$ is locally identified on 
$$\mathcal{N}_{\epsilon}=\{(\varphi_0(u)+\delta \psi_0,\varphi_1(\cdot,u)+\delta \psi_1):\ \delta\in [-\epsilon,\epsilon]\}.$$ 
\end{Theorem}
We have shown here local identification on a segment $\mathcal{N}_{\epsilon}$. As noted in \cite{chen2014local}, in nonparametric nonlinear structural models (as the one of the present paper), it is often not possible to derive local identification results when $\mathcal{N}$ is an open ball (in the topology defined by $\|\cdot\|_{\mathcal{P}}$). We therefore focused on a smaller set, which is not an open ball.}
\subsubsection{Global identification}
\red{Next, we discuss global identification, which we define as follows.
\begin{Definition}The function $\varphi(\cdot,u)$ is globally identified if, for all $(\psi_0,\psi_1)\in\mathcal{P}$, the fact that \eqref{restatemaineq} holds
implies that $\psi(Z)=\varphi(Z,u)$ almost surely.
\end{Definition}
We adapt the theory of \cite{CH} to the case of the present model. We introduce $\epsilon=T-\varphi_0(u)(1-D)- \varphi_1(Z,u)D.$
Let $f_{\epsilon|D,W}(\cdot|0,w)$ be the density of $\epsilon$ given $D=0$ and $W=w$, $f_{\epsilon|D,Z,W}(\cdot|1,z,w)$ be the density of $\epsilon$ given $D=1$, $Z=z$ and $W=w$ and $f_{Z|D,W}(\cdot|d)$ be the density of $Z$ given $D=d$ (their existence is guaranteed under our assumptions).
Let us make the next Assumption:
\begin{Assumption}
\label{Rglob}
The distribution of $(U,Z,W)$ is absolutely continuous with continuous density, $0<\P(D=1)<1$, and there exists a constant $K>0$ such that $f_Z(z)I(z\le \varphi_0(u))/f_{Z|D}(z|1)\le K$.
\end{Assumption}
 For $(\Delta_0,\Delta_1)\in\mathcal{P}$, we define
\begin{align*} \omega_\Delta(Z,D,W)&=\left[\int_0^1f_{\epsilon|D,W}(\delta\Delta_0|0,W)d\delta\right] (1-D)\\
&\quad + \left[\int_0^1f_{\epsilon|D,Z,W}(\delta\Delta_1(Z)|1,Z,W)d\delta\right]D.\end{align*}
We make the following hypothesis:
\begin{Assumption}\label{as.global} For all $(\Delta_0,\Delta_1)\in\mathcal{P}$, we have 
\begin{align*}&E[(\Delta_0(1-D)+\Delta_1(Z)D) \omega_\Delta(Z,D,W)|W]=0\ a.s. \\
&\Rightarrow \Delta_0(1-D)+\Delta_1(Z)D=0\ a.s.
\end{align*}
\end{Assumption}
This is a type of bounded strong completeness condition. It is the counterpart of Assumption L1$^*$ of \cite{CH} in our model. This assumption holds when $Z$ and $W$ are sufficiently dependent given $U$. To the best of our knowledge, the only sufficient conditions known for this type of assumption correspond to condition L2$^*$ in \cite{CH}. In Section \ref{suffglobident} of the online appendix, we give sufficient conditions for Assumption \ref{as.global} which are in the spirit of condition L2$^*$ in \cite{CH}. These conditions include assuming that a family (in $w$) of distributions is boundedly complete and therefore relate Assumptions \ref{as.global} to standard bounded completeness conditions (as in \eqref{complete2}). 
We have the following global identification result.
\begin{Theorem} \label{globident}
Let Assumptions \ref{haz}, \ref{RI}, \ref{IVind}, \ref{Rglob}, \ref{as.global} hold and assume that $(\varphi_0(u),\varphi_1(\cdot,u))\in\mathcal{P}$, 
then $\varphi(\cdot,u)$ is globally identified.
\end{Theorem}}

\subsection{Identification with censoring} \label{cens.sec.id}

Let us now consider the case where $T$ is right censored. In this case, $F_0$ and $F_1$ may not be identified everywhere because of censoring. We make the following assumption on the censoring:
\begin{Assumption}\label{as.censoring}
The censoring time $C$ is independent of $(U,Z,W)$.
\end{Assumption}
This assumption could be relaxed. For instance, the fact that $C$ and $U$ are independent given $Z,W$ would suffice for identification. However, we impose the stronger Assumption \ref{as.censoring} to simplify the exposition and the estimation. 

Let $c_0$ be the upper bound of the support of $C$. Identification of $F_0(t,w)$ and $F_1(\psi,w)$ for $\psi:\ \mathcal{Z}\mapsto [0,t]$ is only possible for $t\in[0,c_0]$. Hence, we can only check that $\varphi(\cdot,u)$ is the solution of the identification equation for $u\in[0,u_0]$, where 
$u_0 =\inf\{u\in\R_+ : \varphi(z,u)< c_0 \: \forall z\in\mathcal{Z}\}.$
Define $G(t)= \P(C\ge t)$, the survival function of $C$. Since $F_0(t,w)= E[I(T\le t,D=0,W\le w)]$, \red{using the law of iterated expectations and Assumption \ref{as.censoring}, we can show that}
\begin{equation}\label{F0cens} F_0(t,w)=E\left[\frac{\delta}{G(Y)} I(Y\le t,\tilde D=0,W\le w)\right],\end{equation}
for all $t\in[0,c_0)$. Similarly, for all $\phi:\mathcal{Z}\mapsto [0,c_0)$, it holds that 
\begin{equation}\label{F1cens} F_1(\phi, w)=E\left[\frac{\delta}{G(Y)} I(Y\le \phi(\tilde Z),\tilde D=1,W\le w)\right].\end{equation}
The proof of \eqref{F0cens} and \eqref{F1cens} is given in Section \ref{proofcens} of the online appendix. By standard arguments from the survival analysis literature, $G(t)$ is identified for all $t\in[0,\sup\{t\in\R_+:\ \P(T\ge t)>0\}]$ \red{(on this interval $G(t)$ is equal to the population analog of the Kaplan-Meier estimator of the survival function of $C$ which identifies it)}. Hence, the quantities on the right-hand side of \eqref{F0cens} and \eqref{F1cens} are identified and so are $F_0(t,w)$ for all $t\in[0,c_0)$ and $ F_1(\phi, w)$ for all $\phi:\mathcal{Z}\mapsto [0,c_0)$. Therefore, identification results on $\varphi(\cdot,u)$ for $u\in[0,u_0]$ can be obtained as in the case without censoring. \red{Note however that the fact that we can identify $\varphi$ only up to $u_0$ has several implications. First, it means that average treatment effects ($E[T(z)-T(z')]$) are not identified. Second, only some quantile treatment effects ($\varphi(z,u)-\varphi(z',u)$) are identified. Third, the structural hazard $\lambda(z,t)=(\varphi(z,\cdot)^{-1})'(t)$ is identified only for $t\le  \varphi(z,u_0)$.}

\section{Estimation}\label{sec.est}

\subsection{Parametric regression function}\label{subsec.param} 

Our strategy is to first estimate $F_0,F_1$ and $F_W$, and then to solve an estimate of equation \eqref{main_eq} for $\varphi_0$ and $\varphi_1$, which is obtained by plug-in of the estimates of $F_0,F_1$ and $F_W$. Since \eqref{main_eq} is a {complicated} integral equation, this method is unlikely to deliver precise nonparametric estimates of $\varphi$ on datasets of reasonable size. As a result, we decide to assume that $(\varphi_0,\varphi_1)$ follows a parametric model $\{\varphi_{\theta0},\varphi_{\theta1}\}_{\theta\in\Theta}$, that is $\varphi_0=\varphi_{\theta_*0}$ and $\varphi_1=\varphi_{\theta_*1}$ for some $\theta_*\in \Theta$, where $\Theta\subset\R^K$ is the parameter set. Here, for all $\theta\in\Theta$, $\varphi_{\theta0}$ is a mapping from $\R_+$ to $\R_+$ and $\varphi_{\theta1}$ is a mapping from $\mathcal{Z}\times \R_+$ to $\R_+$ For all $\theta\in\Theta$, we can also define $\varphi_{\theta}:\mathcal{Z}\times\R_+\mapsto\R_+$ such that $\varphi_{\theta}(z,u)=\varphi_{\theta0}(u)I(z>\varphi_{\theta0}(u))+\varphi_{\theta1}(z,u)I(z\le \varphi_{\theta0}(u))$ for all $z\in\mathcal{Z},u\in\R_+$. This approach has three additional advantages. First, it avoids the need for regularization, since the parameter set has finite dimension. Second, parametric shapes enable to summarize simply the properties of $\varphi$. Finally, they allow to estimate $\varphi(\cdot,u)$ even for $u> u_{0}$, which is useful when the study has insufficient follow-up. Note that, the model remains semiparametric since $F_0$ and $F_1$ are not parametrically constrained. Let us give some examples of parametric models for $\varphi$. \\

\noindent \textbf{Example 1: Weibull model.} A first example comes from the Weibull distribution. Recall that $T(z)$ is the potential outcome of $T$ when the treatment time is set to $z$. We assume that before $z$, the hazard rate of $T(z)$ corresponds to that of a Weibull distribution with parameters $\theta_{00},\theta_{01}$. After $z$, the hazard rate $T(z)$ is that of a Weibull distribution with parameters $\theta_{10},\theta_{11}$. The structural hazard of $T(z)$ at time $t$ is therefore given by 
\begin{equation}\label{hazweibull}
{\color{black} \lambda(z,t) = \theta_{00} \theta_{01}t^{\theta_{01} - 1} I( t< z) + \theta_{01} \theta_{11} t^{\theta_{11} - 1} I(t\ge z)}.
\end{equation}
By inverting the cumulative hazard, it can be shown that 
\begin{align*}\varphi_{\theta 0}(u)=\left( \frac{u}{\theta_{00}} \right)^{\frac{1}{\theta_{01}}};\ \varphi_{\theta 1}(z,u)= \left( \frac{u- \theta_{00}z^{\theta_{01}}}{\theta_{10}} +z^{\theta_{11}} \right)^{\frac{1}{\theta_{11}}}.\\
\end{align*}

\noindent \textbf{Example 2: Log-normal model.} The second example comes from the log-normal distribution. Before $z$, the hazard rate of $T(z)$ is assumed to be equal to that of a log-normal distribution with mean $\theta_{00}$ and variance $\theta_{01}$. After $z$, the hazard rate of $T(z)$ corresponds to that of a log-normal distribution with mean $\theta_{10}$ and variance $\theta_{11}$. The structural hazard rate of $T(z)$ at time $t$ is then given by
\begin{equation}\label{hazlognormal}
\lambda(z,t) = \frac{\phi \left(\frac{\log(t) - \theta_{00}}{\theta_{01}} \right)}{\theta_{01} t \left[ 1 - \Phi \left(\frac{\log(t) - \theta_{00}}{\theta_{01}} \right) \right]} I \left( t < z \right) + \frac{\phi \left(\frac{\log(t) - \theta_{10}}{\theta_{11}} \right)}{\theta_{11} t \left[ 1 - \Phi \left(\frac{\log(t) - \theta_{10}}{\theta_{11}} \right) \right]}  I \left( t \ge z \right),\end{equation}
where $\phi$ and $\Phi$ are respectively the density and the cumulative distribution function of a standard normal distribution. Inverting the cumulative hazard, we obtain
\begin{align*}
\varphi_{\theta 0}(u) =& \exp \left[ \theta_{00} + \theta_{01} \Phi^{-1}\Big( 1 - \exp(-u) \Big) \right]\\
\varphi_{\theta 1}(z,u) =& \exp \left[ \theta_{10} + \theta_{11} \Phi^{-1}\Big(1 - \exp(-u+ \log(R_z(\theta))) \Big) \right],
\end{align*}
where $R_z(\theta) = \left[1 - \Phi \left(\frac{\log(z) - \theta_{10}}{\theta_{11}} \right)\right]\left[1 - \Phi \left(\frac{\log(z) - \theta_{00}}{\theta_{01}} \right)\right]^{-1}.$

\subsection{Estimation of the integral equation} 

For $\theta\in \Theta$, let 
$$M_\theta(u,w) = F_0(\varphi_{\theta 0}(u),w) +F_1(\varphi_{\theta 1}(\cdot,u),w) - (1-e^{-u})F_W(w)$$
be the value of the identifying equation \eqref{main_eq} in $(\varphi_{\theta0},\varphi_{\theta1})$. In order to estimate $\theta_*$ using \eqref{main_eq}, it is necessary to estimate the unknown operator $M$. Assume that we possess an i.i.d. sample $\{Y_i, \delta_i,\tilde Z_i, \tilde D_i, W_i\}_{i=1}^n$. We estimate $F_0$ and $F_1$ using \eqref{F0cens} and \eqref{F1cens}.
 Let 
\begin{align*}
N(t)= \sum_{i=1}^n I(Y_i\le t,\delta_i=0);\ Y(t) = \sum_{i=1}^n I(Y_i\ge t).
\end{align*}
The Kaplan-Meier estimator of $G(t)$ is given by
$\widehat{G}(t) = \prod_{s< t}\left(1-\frac{dN(s)}{Y(s)}\right),$
where $dN(s)=N(s)-\lim\limits_{s'\to s,s'<s}N(s')$. In turn, $F_0$ and $F_1$ are estimated by
\begin{align*}\widehat{F}_{0}(t,w)&=\frac{1}{n}\sum_{i=1}^n \frac{\delta_i}{\widehat{G}(Y_i)} I(Y_i\le t,\tilde D_i=0, W_i\le w);\\\widehat{F}_{1}(\psi,w)&=\frac{1}{n}\sum_{i=1}^n \frac{\delta_i}{\widehat{G}(Y_i)}I(Y_i\le \psi(\tilde Z_i), \tilde D_i=1, W_i\le w) .\end{align*}
Then, $F_W$ is estimated by $\widehat{F}_W(w)= n^{-1}\sum_{i=1}^n I(W_i \le w)$. Finally, the estimator of $M$ is 
$$\widehat{M}_\theta(u,w)=\widehat{F}_0(\varphi_{\theta0}(u),w) +  \widehat{F}_1(\varphi_{\theta1}(\cdot,u), w)- (1-e^{-u})\widehat{F}_W(w). $$
Remark that, although $Z$ and $W$ are continuous random variables, we avoid smoothing because we use an unconditional identification equation.
\subsection{Estimator of $\theta_*$} \label{sec.selectu}

It is computationally impossible to solve the estimated identifying equation for every $u$ and $w$. Hence, we solve it on a grid $0\le u_1<\dots<u_m<u_0$ of values of $u$, where $m\in \N_*$ is fixed, and at $w=W_1,\dots, W_n$. 
 The estimator of $\theta_*$ is 
\begin{equation}\label{estimator}\widehat{\theta}\in\argmin_{\theta \in\Theta}\widehat{L}(\theta),\end{equation}
where $\widehat{L}(\theta)= (nm)^{-1}\sum_{i=1}^n\sum_{j=1}^mp(u_j) \widehat{M}_\theta(u_j,W_i)^2$ and $p(\cdot)$ is some weighting function, typically $p(u)=e^{-u}$. {\color{black}This type of estimators is akin to minimum distance from independence estimators as in \cite{brown2002weighted}. We can not rely directly on the theory of \cite{brown2002weighted}. Indeed, in our case, because of censoring, $G$ has to be estimated in a first step. Moreover, weights on $W$ are chosen according to the empirical measure, while they are set according to a fixed measure selected by the researcher in \cite{brown2002weighted}. This avoids the need to choose weights on $W$ and might lead to greater efficiency. }



\subsection{Asymptotic normality}
Now, we state the conditions that we impose to show asymptotic normality. The first assumption ensures identification, that is 
\begin{equation}\label{ident_main_eq} \theta_*=\argmin_{\theta\in\Theta} L(\theta),\end{equation}
where $L(\theta)=m^{-1}\sum_{j=1}^mp(u_j)E[M_\theta(u_j,W)^2].$
\begin{Assumption}\label{paramid} The following holds
\begin{itemize}
\item[(i)] The regression function $\varphi(\cdot,u_j)$ is globally identified for all $j=1,\dots,m$.
\item[(ii)] We have $\sup\limits_{\theta \in\Theta,z\in\mathcal{Z}}\varphi_{\theta 0}(u_m)\vee \varphi_{\theta 1}(z,u_m)<c_0$.
\item[(iii)] For all $\theta,\tilde\theta\in\Theta$ such that $\varphi_\theta(z,u_j)=\varphi_{\tilde \theta}(z,u_j),$ for all $j\in\{1,\dots, m\},\ z\in\mathcal{Z}$, we have $\theta =\tilde \theta$.
\end{itemize}
\end{Assumption}
Condition (i) was studied in the previous section, whereas condition (ii) restricts the choices of $\Theta$ and $u_m$, and \red{(iii) is a constraint on the parametric family that guarantees that $\theta_*$ is identified when $\varphi(\cdot,u_j)$ is known for all $j\in\{1,\dots, m\}$. Together (i) and (iii) ensure that $\theta$ is identified from \eqref{main_eq}, which yields \eqref{ident_main_eq}. Condition (ii) implies that $\varphi_{\theta_*1}(\cdot,u)=\varphi_1(\cdot,u)$ is bounded when $c_0<\infty$ (that is censoring has finite support).} Let $||\cdot||$ denote the Euclidean norm in $\R^K$. We also impose some regularity conditions:
\begin{Assumption}\label{paramregul} The following holds
\begin{enumerate}[\textup{(}i\textup{)}] 
\item\label{Riv} The true parameter $\theta_*$ is an interior point of $\Theta$.
\item \label{Rv} The parameter space $\Theta$ is compact.
\item\label{Ri} For all $u,w$, the mapping
$\theta\mapsto M_\theta(u,w)$ is three times differentiable and its third order derivative is bounded uniformly in $u,w,\theta$. 
\item\label{Rii} The matrix $\nabla^2L(\theta_*)$ is positive definite.
\item \label{Riii} The class $\{(t,z)\in \R_+\times\mathcal{Z} \mapsto I(t\le \varphi_{\theta d}(z,u)),\ \theta\in\Theta\}$ is Donsker for all $u\in \R_+,d\in\{0,1\}$. 
\item \label{Rvi} For all $u\in\R_+$, there exists a constant $C_u>0$ such that $|\varphi_{\theta d}(z,u)-\varphi_{\theta_*d}(z,u)|\le C_u||\theta-\theta_*||$ for all $\theta \in \Theta,d\in\{0,1\},z\in \mathcal{Z}$. 
\item \label{Rvii} The density of $T$ given $Z,D$ is uniformly bounded.
\end{enumerate}
\end{Assumption}
These standard and mild conditions depend simultaneously on the regularity of the mapping $\theta\mapsto \varphi_\theta$ and on the distribution of $(U,Z,W)$. 

\begin{Theorem} \label{AN}Under Assumptions \ref{haz}, \ref{RI}, \ref{IVind}, \ref{as.censoring}, \ref{paramid} and \ref{paramregul}, there exists a $K\times K$ asymptotic variance matrix $\Sigma$ such that $\sqrt{n}(\widehat{\theta}-\theta_*)\xrightarrow{d}\mathcal{N}(0,\Sigma).$
\end{Theorem}

\subsection{Bootstrap}

Since the asymptotic variance matrix $\Sigma$ has a complicated expression (see the proof of Theorem \ref{AN}), we rely on the nonparametric bootstrap for inference. Let $\{Y_{bi}, \delta_{bi},\tilde Z_{bi}, \tilde D_{bi}, W_{bi}\}_{i=1}^n$ be the bootstrap sample drawn with replacement from the original sample \linebreak $\{X_i=(Y_i, \delta_i,\tilde Z_i, \tilde D_i, W_i)\}_{i=1}^n$. Let also $\widehat{\theta}_b$ be the value of the estimator computed on the bootstrap sample $b$. The following result allows to build confidence intervals with the na\"ive bootstrap using the empirical distribution of the estimates in the bootstrap samples.

\begin{Theorem} \label{ANB}Under Assumptions \ref{haz}, \ref{RI}, \ref{IVind}, \ref{as.censoring}, \ref{paramid} and \ref{paramregul}, it holds that
$$\sqrt{n}(\widehat{\theta}_b-\widehat{\theta})\xrightarrow{d}\mathcal{N}(0,\Sigma) \quad[P],$$
where the convergence is for the law of $\widehat\theta_b$ conditional on the original sample, in probability with respect to the original sample. In other words, 
$$ \sup_t \Big|P^*(\sqrt{n}(\widehat\theta_b - \widehat\theta) \le t) - P(\sqrt{n}(\widehat\theta-\theta_*) \le t) \Big|= o_P(1), $$
where $P^*$ stands for the probability law conditionally on the original data.
\end{Theorem}

\section{Numerical experiments}

\subsection{Simulations} \label{sec:montecarlo}

For a sample of size $n$, and for $i = 1,\dots,n$, we generate the instrument $W_i$ and $U_i$ from two independent exponential distributions with parameter equal to one. We also generate an additional error term $r_i \sim Exp(1)$. The treatment time $Z_i$ is taken equal to $Z_i = \sqrt{2r_i U^{\alpha}_i W^{\beta}_i},$
where the parameter $\alpha$ controls the level of endogeneity, and the parameter $\beta$ controls the strength of the instrument. When $\alpha = 0$, the treatment time, $Z_i$, is independent of $U_i$, and when $\beta=0$, the instrument cannot explain any of the variation in the treatment time. Finally, 
$$
T_i = \varphi(Z_i,U_i) = \varphi_0(U_i) I(\varphi_0(U_i) < Z_i) + \varphi_1(Z_i,U_i) I(\varphi_0(U_i) \ge Z_i).$$

For $\varphi$, we consider the two parametric models of Section \ref{subsec.param}. 
For the Weibull model (example 1), we take the true parameter vector $\theta = (\theta_{00},\theta_{10},\theta_{01},\theta_{11})^\top = (1,2,1.5,2)^\top$. 
For the log-normal model (example 2), the true parameter vector is $\theta = (\theta_{00},\theta_{10},\theta_{01},\theta_{11})^\top = (0,1,1,1)^\top$. 
For both designs, we consider $\alpha \in \{ 0.25,0.75 \}$ to vary the level of endogeneity, and $\beta \in \{ 0.5,1 \}$ to vary the strength of the instrument. We further take censoring into account as follows:

\begin{itemize}
\item[(a)] $T_i$ is not censored ($C_i = \infty$). 
\item[(b)] In Setting 1 (Weibull) we take $C_i-0.3 \sim \mbox{Exp}(2)$, and in Setting 2 (log-normal), $\log(C_i) \sim N(1,1)$. The parameters of the distribution of $C_i$ are chosen in such a way that about $20\%$ of the observations are censored.
\end{itemize}

The sample size is fixed to $n \in \{ 500,1000,3000 \}$. We thus have a total of $3 \times 2^4 = 48$ simulation schemes, and we run $R = 1000$ replications for each scheme. The optimization algorithm is started at 100 random values. Each of these starting values yields a local minimum. The local minimum which leads to the lowest value of the objective function corresponds to the estimate. We use a grid $u_1,\dots,u_{100}$ of $100$ values of $U$, where $u_1$ (respectively $u_{100}$) is the 0.025 (respectively, 0.975) quantile of the unit exponential distribution and the points are equally spaced.

We report the bias and standard error of the estimator and the coverage of the bootstrap percentile confidence intervals at level $90\%,\ 95\%$ and $99\%$. As the computation time increases with the sample size, the coverage of the confidence intervals are evaluated using the Warp-Speed method of \citet{giacomini2013}, where only one bootstrap resampling is used for each simulated sample. We also provide the average number of treated units observed, $\bar{D}$, and the average number of uncensored observations, $\bar{\delta}$.

We summarize the simulations results for the Weibull design with censoring in Table \ref{tab:simW2}. Between $33$ and $43\%$ of observations are treated. With a strong instrument ($\beta=1$), the bias is low, even when $n=500$. The coverage of the confidence intervals improves when $n$ grows and are relatively close to nominal for $n=3000$. When the instrument is weaker ($\beta=0.5$), the estimator exhibits good performance in terms of bias when $n=1000$ or $n=3000$. \red{The results for the Weibull design without censoring and the log-normal model (with and without censoring) are reported in Section \ref{S.3} of the online appendix.}

\afterpage{
\clearpage
  \thispagestyle{empty}
  \begin{landscape}
\begin{table}[!p]
\centering
\adjustbox{max width=1.3\textwidth}
{
\begin{tabular}{l | c c c c | c c c c | c c c c | c c c c}
\hline \hline
~ & $\hat{\theta}_{00}$ & $\hat{\theta}_{10}$ & $\hat{\theta}_{01}$ & $\hat{\theta}_{11}$ & $\hat{\theta}_{00}$ & $\hat{\theta}_{10}$ & $\hat{\theta}_{01}$ & $\hat{\theta}_{11}$ & $\hat{\theta}_{00}$ & $\hat{\theta}_{10}$ & $\hat{\theta}_{01}$ & $\hat{\theta}_{11}$ & $\hat{\theta}_{00}$ & $\hat{\theta}_{10}$ & $\hat{\theta}_{01}$ & $\hat{\theta}_{11}$ \\ \hline \hline
\multicolumn{17}{c}{\textsc{Censoring}}\\ \hline
\multicolumn{17}{l}{\textbf{$n = 500$}}\\ \hline
~ & \multicolumn{4}{c|}{$\alpha = 0.25, \beta = 1, \bar{D} = 0.40, \bar{\delta} = 0.80$} & \multicolumn{4}{c|}{$\alpha = 0.75, \beta = 1, \bar{D} = 0.43, \bar{\delta} = 0.80$} & \multicolumn{4}{c|}{$\alpha = 0.25, \beta = 0.5, \bar{D} = 0.33, \bar{\delta} = 0.80$} & \multicolumn{4}{c}{$\alpha = 0.75, \beta = 0.5, \bar{D} = 0.36, \bar{\delta} = 0.79$}\\
 \hline
Bias & 0.022 & 0.149 & 0.022 & -0.023 & 0.009 & 0.066 & 0.026 & -0.041 & 0.039 & 0.922 & 0.029 & -0.094 & 0.018 & 0.719 & 0.031 & -0.193 \\ 
 SE & 0.247 & 2.366 & 0.247 & 0.551 & 0.261 & 1.079 & 0.262 & 0.548 & 0.271 & 7.003 & 0.257 & 0.913 & 0.282 & 4.552 & 0.267 & 0.923 \\ 
 $90\%$ & 0.808 & 0.818 & 0.816 & 0.812 & 0.823 & 0.844 & 0.805 & 0.761 & 0.720 & 0.698 & 0.780 & 0.634 & 0.687 & 0.673 & 0.749 & 0.538 \\ 
 $95\%$ & 0.868 & 0.944 & 0.867 & 0.902 & 0.879 & 0.948 & 0.872 & 0.881 & 0.811 & 0.854 & 0.853 & 0.788 & 0.787 & 0.820 & 0.842 & 0.724 \\ 
 $99\%$ & 0.954 & 0.991 & 0.944 & 0.987 & 0.967 & 0.997 & 0.990 & 0.989 & 0.973 & 0.969 & 0.976 & 0.923 & 0.938 & 0.936 & 0.967 & 0.955 \\ 
  \hline

\multicolumn{17}{l}{\textbf{$n = 1000$}}\\ \hline
~ & \multicolumn{4}{c|}{$\alpha = 0.25, \beta = 1, \bar{D} = 0.40, \bar{\delta} = 0.80$} & \multicolumn{4}{c|}{$\alpha = 0.75, \beta = 1, \bar{D} = 0.43, \bar{\delta} = 0.80$} & \multicolumn{4}{c|}{$\alpha = 0.25, \beta = 0.5, \bar{D} = 0.33, \bar{\delta} = 0.80$} & \multicolumn{4}{c}{$\alpha = 0.75, \beta = 0.5, \bar{D} = 0.36, \bar{\delta} = 0.79$}\\
 \hline
Bias & 0.009 & 0.010 & 0.000 & 0.004 & -0.009 & 0.026 & 0.011 & -0.015 & 0.028 & 0.182 & 0.004 & -0.003 & 0.017 & 0.160 & 0.002 & -0.072 \\ 
 SE & 0.227 & 0.393 & 0.213 & 0.382 & 0.245 & 0.401 & 0.230 & 0.391 & 0.229 & 1.924 & 0.207 & 0.654 & 0.247 & 1.841 & 0.231 & 0.610 \\ 
 $90\%$ & 0.873 & 0.853 & 0.861 & 0.852 & 0.853 & 0.825 & 0.854 & 0.814 & 0.762 & 0.653 & 0.814 & 0.661 & 0.747 & 0.653 & 0.817 & 0.678 \\ 
 $95\%$ & 0.901 & 0.926 & 0.899 & 0.936 & 0.882 & 0.929 & 0.874 & 0.915 & 0.858 & 0.820 & 0.891 & 0.836 & 0.842 & 0.818 & 0.885 & 0.830 \\ 
 $99\%$ & 0.966 & 0.999 & 0.982 & 0.999 & 0.960 & 0.997 & 0.955 & 0.998 & 0.959 & 0.975 & 0.969 & 0.955 & 0.976 & 0.968 & 0.933 & 0.984 \\ 
  \hline

\multicolumn{17}{l}{\textbf{$n = 3000$}}\\ \hline
~ & \multicolumn{4}{c|}{$\alpha = 0.25, \beta = 1, \bar{D} = 0.40, \bar{\delta} = 0.80$} & \multicolumn{4}{c|}{$\alpha = 0.75, \beta = 1, \bar{D} = 0.43, \bar{\delta} = 0.80$} & \multicolumn{4}{c|}{$\alpha = 0.25, \beta = 0.5, \bar{D} = 0.33, \bar{\delta} = 0.80$} & \multicolumn{4}{c}{$\alpha = 0.75, \beta = 0.5, \bar{D} = 0.36, \bar{\delta} = 0.79$}\\
 \hline
Bias & -0.003 & 0.002 & 0.011 & 0.009 & -0.008 & 0.007 & 0.018 & 0.001 & 0.002 & 0.012 & -0.002 & 0.005 & -0.004 & 0.031 & 0.010 & -0.010 \\ 
 SE & 0.210 & 0.277 & 0.209 & 0.265 & 0.218 & 0.266 & 0.208 & 0.269 & 0.220 & 0.483 & 0.201 & 0.368 & 0.215 & 0.494 & 0.208 & 0.365 \\ 
 $90\%$ & 0.891 & 0.872 & 0.888 & 0.871 & 0.894 & 0.863 & 0.885 & 0.865 & 0.866 & 0.801 & 0.873 & 0.809 & 0.863 & 0.763 & 0.875 & 0.835 \\ 
 $95\%$ & 0.903 & 0.929 & 0.905 & 0.930 & 0.902 & 0.922 & 0.901 & 0.925 & 0.908 & 0.883 & 0.906 & 0.895 & 0.898 & 0.877 & 0.905 & 0.892 \\ 
 $99\%$ & 0.960 & 0.973 & 0.969 & 0.989 & 0.970 & 0.967 & 0.958 & 0.984 & 0.966 & 0.999 & 0.972 & 0.999 & 0.963 & 0.990 & 0.966 & 0.992 \\ 
  \hline

\hline
\end{tabular}
}
\caption{Simulation results under the Weibull design with censoring.}
\label{tab:simW2}
\end{table}
  \end{landscape}
  \clearpage
}

\subsection{Empirical application}\label{empirical_application}

We use data from a large Belgian insurance company. This insurer offers a product to other companies, which consists in paying for the salaries of their client' workers who are on medical leave because of burnout. This insurance product also contains the possibility of following a free therapy for these workers suffering from burnout. The goal of the company is to reduce the duration of medical leave (our duration variable in this application) and, hence, it would like to know the effect of the start of the therapy (our treatment) on the duration of medical leave. 

We take advantage of the treatment assignment mechanism to evaluate the causal effect of the treatment. The employees on medical leave are assigned to one of several partner institutions for medical care, which are responsible for carrying out the therapy. These partners evaluate the patient and decide to propose or not the therapy on the basis of the expected medical benefits of the therapy. The employee may then accept or decline to be treated, which suggests that $Z$ is endogenous. The exact date of the start of the therapy varies depending on the availabilities of the partner and the patient, which makes the treatment time-varying. 

The assignment by the insurance company of the partner was only based on geographical distance, which should make it exogenous. Moreover, partners are more or less likely to offer a therapy to assigned patients, hence the propensity of the partner to offer treatment has an impact on the treatment time $Z$. We use the proportion of patients who are given the possibility to be treated by the partner in a given year as the instrumental variable. This choice of instrument is common in medical studies (see \cite{brookhart2006evaluating,chen2011use} for reviews). 

The sample consists of 838 individuals who entered medical leave for burnout between 2017 and 2020. Their age at the start of the medical leave was between 30 and 39 years old. They are observed until one of the following events happens: their medical leave ends, their medical leave exceeds 2 years, or the study is ended (at the end of 2020). In the latter two cases the individual is censored. Since the duration of follow-up depends on external factors, we expect censoring to be uninformative. Around 41\% of the observations are censored and 48\% are treated before censoring. The average (respectively, median) duration of medical leave for uncensored observations is 189 days (respectively, 159 days). For the observations for which the treatment time is uncensored, the average is 112 days and the median is equal to 90 days.

For both the Weibull model and the log-normal model, we estimate the parameters using 100 random starting values around the values corresponding to the na\"ive fit of the Weibull distribution or the log-normal distribution to the data of durations and censoring indicators. We use a grid $u_1,\dots,u_{100}$ of $100$ values of $U$, where $u_1$ is the 0.025 quantile of the unit exponential distribution and the points are equally spaced. The upper bound $u_{100}$ is chosen such that
\begin{equation}\label{postcheck}\varphi_{\widehat{\theta}0}(u_{100})<c_0\text{ and } \varphi_{\widehat{\theta}1}(\tilde Z_i,u_{100})< c_0,\ i\in\{1,\dots, n:\ \tilde D_i=1\},\end{equation}
where $c_0$, the maximum follow-up time, is equal to 2 years. Condition \eqref{postcheck} ensures that Assumption (I) (ii) is satisfied for some $\Theta$ in a neighborhood of $\widehat{\theta}$. Following this approach, we chose $u_{100}$ equal to the 0.9 (respectively, 0.8) quantile of a unit exponential for the Weibull (respectively, log-normal) model. The curves of the hazard rates for the individuals who are never treated ($z=\infty$) and those who received treatment at time 0 ($z=0$) corresponding to the estimated parameters are plotted in Figure \ref{fig:cov_design1} for the Weibull model and in Figure \ref{fig:cov_design2} for the log-normal model. 

These hazard rates are computed by plug-in of the estimates in \eqref{hazweibull} and \eqref{hazlognormal}. By definition of our Weibull and log-normal models, for an arbitrary value of $z$ the structural hazard rate at $t$ under treatment at time $z$ is equal to that of the never treated for $t \le z$ and is equal to that of the treated at time $0$ for $t>z$. Hence, in Figures \ref{fig:cov_design1} and \ref{fig:cov_design2}, the estimated structural hazard under treatment $z$ corresponds to the red curve before $z$ and then ``jumps'' to the blue dashed curve.

We see also that the therapy appears to increase the hazard rate for all possible treatment timings. As a result, the treatment should be administered at time $0$ in order to minimize the duration of medical leave.

 Bootstrap confidence intervals for the hazard rates are given in Section \ref{S.4} of the online appendix. Note that the estimated hazards exhibit different shapes under the two models. This is to be expected since the models assume different parametric forms. However, the treatment significantly increases the hazard rate in both models (see the bootstrap confidence intervals in the online appendix). The robustness of this conclusion to the choice of the model constitutes statistical evidence supporting the efficacy of the therapy.

\begin{figure}
\begin{minipage}{0.48\textwidth}
 \centering
 \includegraphics[width=80mm]{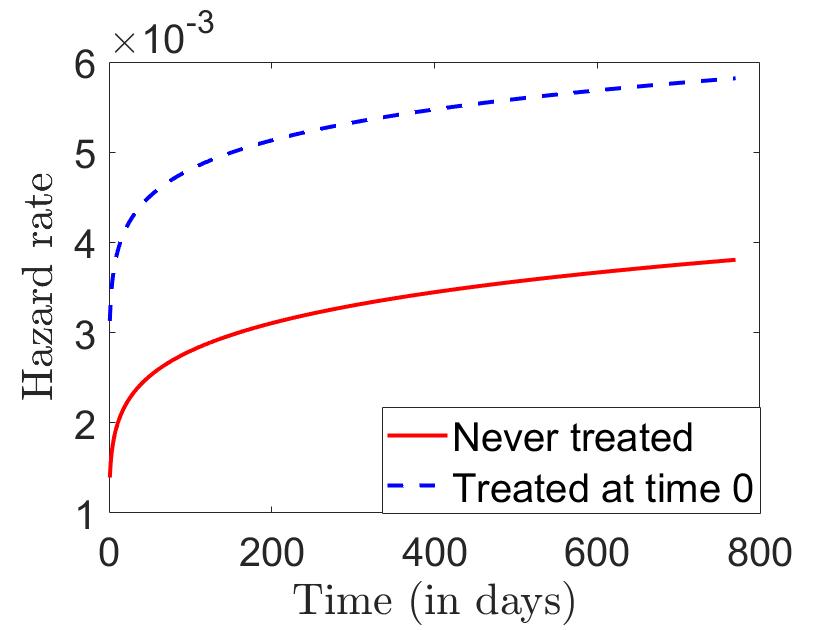}
 
 \caption{Estimated hazard rates with the Weibull model.}
  \label{fig:cov_design1} 
  \end{minipage}
  \quad\quad
  \begin{minipage}{0.48\textwidth}
    \centering
  \includegraphics[width=80mm]{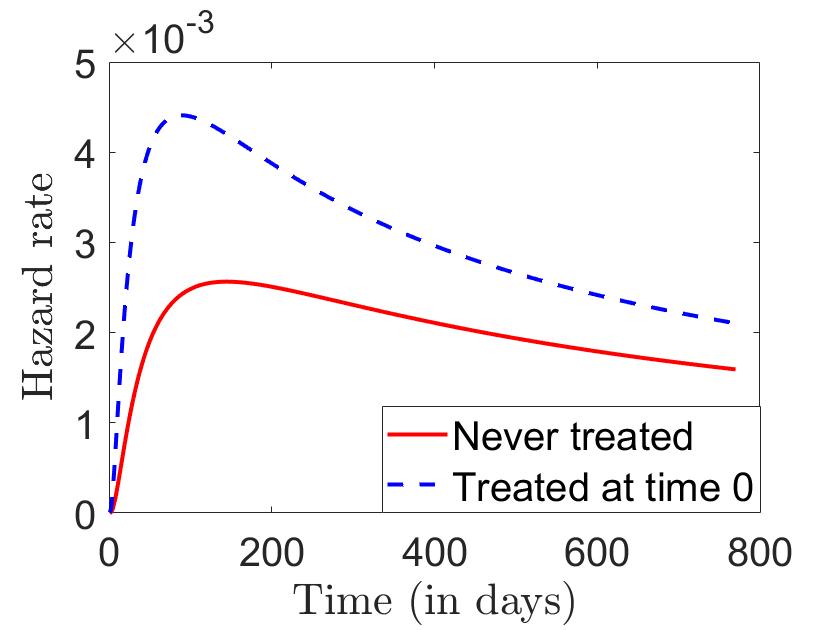}
 \caption{Estimated hazard rates with the log-normal model.}
   \label{fig:cov_design2}
   \end{minipage}
\end{figure}

\section{Concluding remarks}

This paper develops an instrumental variable approach to estimate the causal effect of the time until a treatment is started on a possibly right-censored duration outcome. Therefore, the treatment $Z$ corresponds to the jump of a counting process with single jump. As an extension, it would be of interest to consider procedures where $Z$ is a more general process. Another possible research direction could be the development of a control function approach in the context of the present paper. \red{This might allow us to use a discrete instrument, as in \cite{xavier2015,torgovitsky2015}.} 

\bibliographystyle{apalike}
\bibliography{JMP-Beyhum}

\begin{thebibliography}{}

\bibitem[Abbring and Van~den Berg, 2003]{AVdB}
Abbring, J.~H. and Van~den Berg, G.~J. (2003).
\newblock The nonparametric identification of treatment effects in duration
  models.
\newblock {\em Econometrica}, 71(5):1491--1517.

\bibitem[Abbring and Van~den Berg, 2005]{AVdB2}
Abbring, J.~H. and Van~den Berg, G.~J. (2005).
\newblock Social experiments and instrumental variables with duration outcomes.
\newblock {\em Technical Report}.

\bibitem[Andrews, 2017]{andrews2011}
Andrews, D. W.~K. (2017).
\newblock {Examples of {$L^2$}-Complete and Boundedly-Complete Distributions}.
\newblock {\em Journal of Econometrics}, 199(2):213--220.

\bibitem[Angrist et~al., 1996]{angrist1996identification}
Angrist, J.~D., Imbens, G.~W., and Rubin, D.~B. (1996).
\newblock Identification of causal effects using instrumental variables.
\newblock {\em Journal of the American Statistical Association},
  91(434):444--455.

\bibitem[Beyhum et~al., 2022]{2021nonparametric}
Beyhum, J., Florens, J.-P., and Keilegom, I.~V. (2022).
\newblock Nonparametric instrumental regression with right censored duration
  outcomes.
\newblock {\em Journal of Business \& Economic Statistics}, 40(3):1034--1045.

\bibitem[Bickel and Freedman, 1981]{bickel1981asymptotic}
Bickel, P.~J. and Freedman, D. (1981).
\newblock Asymptotic theory for the bootstrap.
\newblock {\em Annals of Statistics}, 9(6):1196--1217.

\bibitem[Bijwaard and Ridder, 2005]{bijwaard2005correcting}
Bijwaard, G.~E. and Ridder, G. (2005).
\newblock Correcting for selective compliance in a re-employment bonus
  experiment.
\newblock {\em Journal of Econometrics}, 125(1-2):77--111.

\bibitem[Blanco et~al., 2020]{blanco2019bounds}
Blanco, G., Chen, X., Flores, C.~A., and Flores-Lagunes, A. (2020).
\newblock Bounds on average and quantile treatment effects on duration outcomes
  under censoring, selection, and noncompliance.
\newblock {\em Journal of Business \& Economic Statistics}, 38:901--920.

\bibitem[Bose and Chatterjee, 2018]{bose2018u}
Bose, A. and Chatterjee, S. (2018).
\newblock {\em U-statistics, Mm-estimators and Resampling}.
\newblock Springer.

\bibitem[Brookhart et~al., 2006]{brookhart2006evaluating}
Brookhart, M.~A., Wang, P., Solomon, D.~H., and Schneeweiss, S. (2006).
\newblock Evaluating short-term drug effects using a physician-specific
  prescribing preference as an instrumental variable.
\newblock {\em Epidemiology (Cambridge, Mass.)}, 17(3):268.

\bibitem[Brown and Wegkamp, 2002]{brown2002weighted}
Brown, D.~J. and Wegkamp, M.~H. (2002).
\newblock Weighted minimum mean-square distance from independence estimation.
\newblock {\em Econometrica}, 70(5):2035--2051.

\bibitem[Cazals et~al., 2016]{cazals2016nonparametric}
Cazals, C., F{\`e}ve, F., Florens, J.-P., and Simar, L. (2016).
\newblock Nonparametric instrumental variables estimation for efficiency
  frontier.
\newblock {\em Journal of econometrics}, 190(2):349--359.

\bibitem[Centorrino and Florens, 2021]{centorrinoflorens2021}
Centorrino, S. and Florens, J.-P. (2021).
\newblock {Nonparametric estimation of accelerated failure-time models with
  unobservable confounders and random censoring}.
\newblock {\em Electronic Journal of Statistics}, 15(2):5333 -- 5379.

\bibitem[Chen et~al., 2014]{chen2014local}
Chen, X., Chernozhukov, V., Lee, S., and Newey, W.~K. (2014).
\newblock Local identification of nonparametric and semiparametric models.
\newblock {\em Econometrica}, 82(2):785--809.

\bibitem[Chen and Briesacher, 2011]{chen2011use}
Chen, Y. and Briesacher, B.~A. (2011).
\newblock Use of instrumental variable in prescription drug research with
  observational data: a systematic review.
\newblock {\em Journal of Clinical Epidemiology}, 64(6):687--700.

\bibitem[Chernozhukov et~al., 2015]{chernozhukov2015quantile}
Chernozhukov, V., Fern{\'a}ndez-Val, I., and Kowalski, A.~E. (2015).
\newblock Quantile regression with censoring and endogeneity.
\newblock {\em Journal of Econometrics}, 186(1):201--221.

\bibitem[Chernozhukov and Hansen, 2005]{CH}
Chernozhukov, V. and Hansen, C. (2005).
\newblock An {IV} model of quantile treatment effects.
\newblock {\em Econometrica}, 73(1):245--261.

\bibitem[Cui et~al., 2022]{cui2017instrumental}
Cui, Y., Michael, H., Tanser, F., and Tchetgen~Tchetgen, E. (2022).
\newblock Instrumental variable estimation of the marginal structural cox model
  for time-varying treatments.
\newblock {\em Biometrika (forthcoming)}.

\bibitem[Darolles et~al., 2011]{Darolles}
Darolles, S., Fan, Y., Florens, J.-P., and Renault, E. (2011).
\newblock Nonparametric instrumental regression.
\newblock {\em Econometrica}, 79(5):1541--1565.

\bibitem[D'Haultfoeuille, 2011]{xavier2011}
D'Haultfoeuille, X. (2011).
\newblock {On the Completeness Condition in Nonparametric Instrumental
  Problems}.
\newblock {\em Econometric Theory}, 27:460--471.

\bibitem[D'Haultf{\oe}uille and F{\'e}vrier, 2015]{xavier2015}
D'Haultf{\oe}uille, X. and F{\'e}vrier, P. (2015).
\newblock {Identification of Nonseparable Triangular Models With Discrete
  Instruments}.
\newblock {\em Econometrica}, 83(3):1199--1210.

\bibitem[Frandsen, 2015]{frandsen2015treatment}
Frandsen, B.~R. (2015).
\newblock Treatment effects with censoring and endogeneity.
\newblock {\em Journal of the American Statistical Association},
  110(512):1745--1752.

\bibitem[Giacomini et~al., 2013]{giacomini2013}
Giacomini, R., Politis, D.~N., and White, H. (2013).
\newblock A {W}arp-{S}peed method for conducting {M}onte {C}arlo experiments
  involving bootstrap estimators.
\newblock {\em Econometric Theory}, 29(3):567--589.

\bibitem[Gill, 1983]{gill1983large}
Gill, R. (1983).
\newblock Large sample behaviour of the product-limit estimator on the whole
  line.
\newblock {\em Annals of Statistics}, 11:49--58.

\bibitem[Heckman and Navarro, 2007]{heckman2007}
Heckman, J.~J. and Navarro, S. (2007).
\newblock {Dynamic discrete choice and dynamic treatment effects}.
\newblock {\em Journal of Econometrics}, 136(2):341 -- 396.

\bibitem[Hern{\'a}n and Robins, 2010]{hernan2010causal}
Hern{\'a}n, M.~A. and Robins, J.~M. (2010).
\newblock {\em Causal Inference: What If}.
\newblock Chapman \& Hall/CRC Boca Raton, FL.

\bibitem[Hu and Shiu, 2018]{hu2018nonparametric}
Hu, Y. and Shiu, J.-L. (2018).
\newblock Nonparametric identification using instrumental variables: sufficient
  conditions for completeness.
\newblock {\em Econometric Theory}, 34(3):659--693.

\bibitem[Kastoryano and van~der Klaauw, 2022]{kastoryano2022dynamic}
Kastoryano, S. and van~der Klaauw, B. (2022).
\newblock Dynamic evaluation of job search assistance.
\newblock {\em Journal of Applied Econometrics}, 37(2):227--241.

\bibitem[Kosorok, 2008]{kosorok2008introduction}
Kosorok, M.~R. (2008).
\newblock {\em Introduction to Empirical Processes and Semiparametric
  Inference.}
\newblock Springer.

\bibitem[Lechner, 2009]{lechner2009sequential}
Lechner, M. (2009).
\newblock Sequential causal models for the evaluation of labor market programs.
\newblock {\em Journal of Business \& Economic Statistics}, 27(1):71--83.

\bibitem[Lechner and Miquel, 2010]{lechner2010identification}
Lechner, M. and Miquel, R. (2010).
\newblock Identification of the effects of dynamic treatments by sequential
  conditional independence assumptions.
\newblock {\em Empirical Economics}, 39(1):111--137.

\bibitem[Lo and Singh, 1986]{lo1986product}
Lo, S.-H. and Singh, K. (1986).
\newblock The product-limit estimator and the bootstrap: some asymptotic
  representations.
\newblock {\em Probability Theory and Related Fields}, 71(3):455--465.

\bibitem[Martinussen et~al., 2019]{martinussen2019instrumental}
Martinussen, T., N{\o}rbo~S{\o}rensen, D., and Vansteelandt, S. (2019).
\newblock Instrumental variables estimation under a structural cox model.
\newblock {\em Biostatistics}, 20(1):65--79.

\bibitem[Michael et~al., 2020]{michael2020instrumental}
Michael, H., Cui, Y., Lorch, S., and Tchetgen, E.~T. (2020).
\newblock Instrumental variable estimation of marginal structural mean models
  for time-varying treatment.
\newblock {\em arXiv preprint arXiv:2004.11769}.

\bibitem[Newey and McFadden, 1994]{newey1994large}
Newey, W.~K. and McFadden, D. (1994).
\newblock Large sample estimation and hypothesis testing.
\newblock {\em Handbook of Econometrics}, 4:2111--2245.

\bibitem[Newey and Powell, 2003]{newey2003}
Newey, W.~K. and Powell, J.~L. (2003).
\newblock {Instrumental Variable Estimation of Nonparametric Models}.
\newblock {\em Econometrica}, 71(5):1565--1578.

\bibitem[Richardson et~al., 2017]{richardson2017nonparametric}
Richardson, A., Hudgens, M.~G., Fine, J.~P., and Brookhart, M.~A. (2017).
\newblock Nonparametric binary instrumental variable analysis of competing
  risks data.
\newblock {\em Biostatistics}, 18(1):48--61.

\bibitem[Sant'Anna, 2016]{sant2016program}
Sant'Anna, P.~H. (2016).
\newblock Program evaluation with right-censored data.
\newblock {\em arXiv preprint arXiv:1604.02642}.

\bibitem[Sant'Anna, 2021]{sant2021nonparametric}
Sant'Anna, P.~H. (2021).
\newblock Nonparametric tests for treatment effect heterogeneity with duration
  outcomes.
\newblock {\em Journal of Business \& Economic Statistics}, 39(3):816--832.

\bibitem[Sianesi, 2004]{sianesi2004evaluation}
Sianesi, B. (2004).
\newblock An evaluation of the {S}wedish system of active labor market programs
  in the 1990s.
\newblock {\em Review of Economics and Statistics}, 86(1):133--155.

\bibitem[Tchetgen et~al., 2018]{tchetgen2018marginal}
Tchetgen, E. J.~T., Michael, H., and Cui, Y. (2018).
\newblock Marginal structural models for time-varying endogenous treatments: A
  time-varying instrumental variable approach.
\newblock {\em arXiv preprint arXiv:1809.05422}.

\bibitem[Tchetgen et~al., 2015]{tchetgen2015instrumental}
Tchetgen, E. J.~T., Walter, S., Vansteelandt, S., Martinussen, T., and Glymour,
  M. (2015).
\newblock Instrumental variable estimation in a survival context.
\newblock {\em Epidemiology (Cambridge, Mass.)}, 26(3):402.

\bibitem[Torgovitsky, 2015]{torgovitsky2015}
Torgovitsky, A. (2015).
\newblock {Identification of Nonseparable Models Using Instruments With Small
  Support}.
\newblock {\em Econometrica}, 83(3):1185--1197.

\bibitem[Van~den Berg et~al., 2020a]{VdBBM}
Van~den Berg, G.~J., Bonev, P., and Mammen, E. (2020a).
\newblock Nonparametric instrumental variable methods for dynamic treatment
  evaluation.
\newblock {\em Review of Economics and Statistics}, 102(2):355--367.

\bibitem[Van~den Berg et~al., 2020b]{van2020policy}
Van~den Berg, G.~J., Bozio, A., and Costa~Dias, M. (2020b).
\newblock Policy discontinuity and duration outcomes.
\newblock {\em Quantitative Economics}, 11(3):871--916.

\bibitem[Van~der Vaart and Wellner, 1996]{VdVW}
Van~der Vaart, A.~W. and Wellner, J.~A. (1996).
\newblock {\em Weak Convergence and Empirical Processes with Applications to
  Statistics}.
\newblock Springer.

\bibitem[Vikstr{\"o}m, 2017]{vikstrom2017dynamic}
Vikstr{\"o}m, J. (2017).
\newblock Dynamic treatment assignment and evaluation of active labor market
  policies.
\newblock {\em Labour Economics}, 49:42--54.

\end{thebibliography}

\appendix

\numberwithin{equation}{section}

\numberwithin{Assumption}{section}
\renewcommand*{\thesection}{\Alph{section}}
\renewcommand*{\thesection}{S.\arabic{section}}
\setcounter{equation}{0}


\newpage

   \noindent {\LARGE\bf Online appendix to ``Instrumental variable estimation of dynamic treatment effects on a duration outcome"}\\

Section \ref{S.0} contains additional identification results and the proof of results from Sections 2 and 3 of the paper.  The proof of Theorem 4.1 is in Section \ref{S.1}. Section \ref{S.2} corresponds to the proof of Theorem 4.2. Additional simulation results are given in Section \ref{S.3}. Bootstrap confidence intervals for the empirical application can be found in Section \ref{S.4}.
\section{Proofs of identification results} \label{S.0}
\subsection{Consequences of Assumption 2.2}
We show the following lemma. 
\begin{Lemma} \label{ls11} Under Assumption \ref{RI} (ii), we have 
\begin{itemize}
\item[(i)] For all $z\in\bar{\R}_+$ and $u\in\R_+$, there exists a unique $t_{z,u}$ in the support of $T(z)$ such that $\Lambda(z,t_{z,u})=u$.
\item[(ii)] For all $z\in\bar{\R}_+$, $U(z)\sim\text{Exp}(1)$.
\end{itemize}
\end{Lemma}
\begin{Proof} First, we show (i). Let us denote by $\mathcal{T}_z$ the support of $T(z)$. Remember that, by definition, $\mathcal{T}_z$ is equal to $\{t\in\R_+:\ \lambda(z,t)>0\}$. 
 Let us show that for all $u\in\R_+$, there exists $t_{z,u}\in\mathcal{T}_z$ such that $\Lambda(z,t_{z,u})=u$. First, note that by Bayes' theorem, we have $\Lambda(z,t)=-S_z'(t)/S_z(t)$, where $S_z$ is the survival function of $T(z)$.
Since $\lim\limits_{t\to\infty} S_z(t)=0$, this yields $\lim\limits_{t\to\infty} \Lambda(z,t)=\infty$. Because $\Lambda(z,0)=0$ and $\Lambda(z,\cdot)$ is continuous, by the intermediate value theorem, there, therefore, exists $t_*\in\R_+$ such that $\Lambda(z,t_*)=u$. Let us now define $t_{z,u}=\sup\{t\in\mathcal{T}_z:\ \Lambda(z,t)\le \Lambda(z,t_*)\}$. We claim that $\Lambda(z,t_{z,u})=u$. To show it, we reason by contradiction. If we had $\Lambda(z,t_{z,u})<u=\Lambda(z,t_*)$, then we would have $\Lambda(z,t_*)-\Lambda(z,t_{z,u})=\int_{t_{z,u}}^{t_*}\lambda(z,s)ds >0$. This would imply that there exists $s_*\in(t_{z,u},t_*)$ such that $\lambda(z,s_*)>0$ (otherwise $\int_{t_{z,u}}^{t_*}\lambda(z,s)ds =0$). In this case, because $\Lambda(z,\cdot)$ is increasing, we would have $s_*\in\{t\in\mathcal{T}_z:\ \Lambda(z,t)\le \Lambda(z,t_*)\}$ which contradicts the fact that $t_{z,u}=\sup\{t\in\mathcal{T}_z:\ \Lambda(z,t)\le \Lambda(z,t_*)\}$. We conclude that $\Lambda(z,t_{z,u})=u$.

It remains to show that $t_{z,u}$ is unique. Let $t_{z,u}'\in\mathcal{T}_z$ such that $\Lambda(z,t_{z,u}')=u$. Without loss of generality assume $t_{z,u}\le t_{z,u}'$. We have 
\begin{equation}\label{uniquelam}
\int_{t_{z,u}}^{t_{z,u}'}\lambda(z,s)ds =0
\end{equation} Now, there are two different cases. 

\textit{Case 1: $t_{z,u}<z$.} Since $\lambda(z,\cdot)$ is continuous on $[0,z)$ and $\lambda(z,t_{z,u})>0$, $\lambda(z,\cdot)$ is larger than a constant $\mu>0$ on an interval $[t_{z,u},t_{z,u}+\eta],\eta>0$. This implies that $\int_{t_{z,u}}^{t_{z,u}'}\lambda(z,s)ds\ge \int_{t_{z,u}}^{\min(t_{z,u}', \eta)}\lambda(z,s)ds\ge \mu (\min(\eta,| t_{z,u}'-t_{z,u}|))$, which implies $ t_{z,u}'=t_{z,u}$ by \eqref{uniquelam}.

\textit{Case 2: $t_{z,u}\ge z$.} Since $\lambda(z,\cdot)$ is continuous on $[z,\infty)$ and $\lambda(z,t_{z,u})>0$, $\lambda(z,\cdot)$ is larger than a constant $\mu>0$ on a neighborhood $[t_{z,u},t_{z,u}+\eta],\eta>0$ of $t_{z,u}$. This implies that $\int_{t_{z,u}}^{t_{z,u}'}\lambda(z,s)ds\ge \mu (\min(\eta,| t_{z,u}'-t_{z,u}|))$, which implies $ t_{z,u}'=t_{z,u}$ by \eqref{uniquelam}.\\

Second, we show (ii). For $u\in\R_+$, let $t_{z,u}$ be the unique element of $\mathcal{T}_z$ such that $\Lambda(z,t_{z,u})=u$. We have
\begin{align*}
\P(\Lambda(z,T(z))\ge u)&=\P(\Lambda(z,T(z))\ge \Lambda(z,t_{z,u}))\\
&=\P(T(z)\ge t_{z,u})\\
&=S_z(t_{z,u})\\
&=e^{-\Lambda(z,t_{z,u})}=e^{-u},
\end{align*}
where $S_z$ is the survival function of $T(z)$ and in the fourth equality we used $\Lambda(z,t)=-\log(S_z(t))$.  To conclude, just recall that $u\in\R_+\mapsto e^{-u}$ is the survival function of a unit exponential distribution.
\end{Proof}
\subsection{Proof of Lemma 2.1}\label{proofweird}
We prove first equation \eqref{model} in the main text. To do so, we show the more general result, that for all $z\in\bar{\R}_+$, we have
\begin{equation}\label{moregeneral}T(z)=\varphi(z,U(z))=I(z>\varphi_0(U(z)))\varphi_0(U(z))+I(z\le\varphi_0(U(z)))\varphi_1(z,U(z)).\end{equation}
Equation \eqref{model} in the main text is a direct consequence of \eqref{moregeneral} and Assumption \ref{RI} (i). 
To prove \eqref{moregeneral}, it suffices to show that (i) on the event $\{z> \varphi_0(U(z))\}$ it holds that $T=\varphi_0(U(z))$ and (ii) on the event $\{z\le \varphi_0(U(z))\}$ it holds that $T(z)=\varphi_1(z,U(z))$. Condition (ii) follows directly from the definition of $\varphi_1$. Let us show (i). On the event $\{z> \varphi_0(U(z))\}$, we have $U(z)<\int_0^z\lambda(\infty,s) ds$ by definition of $\varphi_0$. Since $\Lambda(z,T(z))=U(z)$, this yields $\Lambda(z,T(z)) <\int_0^z\lambda(\infty,s) ds$. Using Assumption \ref{haz}, we get
$$\Lambda(z,T(z))=\int_0^{\min(T(z),z)}\lambda(\infty,s)ds +  \int_{\min(T(z),z)}^{T(z)}\lambda(z,s)ds <\int_0^z\lambda(\infty,s) ds,$$
 which is only possible if $T(z)<  z$. This yields $U(z)=\Lambda(z,T(z))=\int_0^{T(z)}\lambda(\infty,s)ds $ on the event $\{z> \varphi_0(U(z))\}$ and therefore $T(z)=\varphi_0(U(z))$ on this event by definition of $\varphi_0$. This concludes the proof of equation \eqref{model}. Remark that since $T(z)=\varphi(z,U(z))$ and $U(z)=\Lambda(z,T(z))$, we indeed have $\varphi(z,\cdot)=\Lambda(z,\cdot)^{-1}$.\\

It remains to show that $\varphi_1(z,\varphi_0^{-1}(z))=z$ for all $z\in\R_+$.
By Assumption \ref{haz}, we have $\Lambda(z,z)= \int_{0}^z\lambda(\infty,s)ds={\color{black} \varphi_0^{-1}(z)}$. By definition of $\varphi$, this directly implies $\varphi_1(z,\varphi_0^{-1}(z))= \Lambda(z,\cdot)^{-1}(\Lambda(z,z))=z$.

\subsection{Proof of Theorem \ref{th.main_eq}}
First remark that the events $\{D=0\}$ and $\{Z>\varphi_0(U)\}$ are equal. Indeed, we have
\begin{align*}
\{Z>\varphi_0(U)\}
&=\{\varphi_0^{-1}(Z)>U\}\\
&=\{\varphi_1(Z,\varphi_0^{-1}(Z))>\varphi_1(Z,U)\}\\
&=\{Z>\varphi_1(Z,U)\}\\
&= \{Z>T\}=\{D=0\},
\end{align*}
where the third equality is due to the fact that $\varphi_1(z,\varphi_0^{-1}(z))=z$ for all $z\in\bar{\R}_+$ (Lemma 2.1) and the fourth equality comes from $\Lambda(Z,T)=U$ and $\varphi_1(Z,\cdot)=\Lambda(Z,\cdot)^{-1}$ a.s. 

As a result, on the event $\{D=0\}$, we have $\varphi(Z,U)=\varphi_0(U)$ a.s. and on the event $\{D=1\}$, it holds that $\varphi(Z,U)=\varphi_1(Z,U)$ a.s. Then, we obtain\begin{align*}
&F_0(\varphi_0(u),w) + F_1(\varphi_1(\cdot,u),w) \\
&=\P(T\le \varphi_0(u), D=0,W\le w) + E[I(T\le \varphi_1(Z,u),D=1, W\le w) ]\\
&= \P(T\le \varphi(Z,u), D=0,W \le w) + E[I(T\le \varphi(Z,u),D=1, W\le w) ].
\end{align*}
Using Assumptions \ref{haz} and \ref{RI}, we get
\begin{align*}
F_0(\varphi_0(u),w) + F_1(\varphi_1(\cdot,u),w)
&=\P(U\le u, D=0,W\le w) + E[I(U\le u,D=1, W\le w) ]\\
&= \P(U\le u,D=0,W\le w)+\P(U\le u,D=1,W\le w)\\
&=\P(U\le u,W\le w).
\end{align*}
Since $U$ and $W$ are independent and $U\sim\text{Exp}(1)$, $\P(U\le u,W\le w)=\P(U\le u)F_W(w)=(1-e^{-u})F_W(w)$, which concludes the proof.

\subsection{Alternative characterization in terms of conditional hazard rates and survival functions of the observed duration $T$}\label{altchar}
In order to formulate this other characterization, we introduce further notations.
For $t\in\R_+$, $z\in \mathcal{Z}$, $w\in\mathcal{W}$, let
\begin{align*}h_0(t|w)&=\lim_{dt\to 0}\frac{\P\left( T\in[t,t+dt]| D=0, W=w\right)}{\P(T\ge t|D=0, W= w)dt};\\
h_1(t|z, w)&=\lim_{dt\to 0}\frac{\P\left( T\in[t,t+dt]| Z=z, D=1,W= w\right)}{\P(T\ge t|Z=z, D=1, W= w)dt}
\end{align*} be the hazard rate of \tr{the observed duration} $T$ at $t$ given $ \{D=0, W=w\}$ and $\{Z=z,  D=1, W= w\}$, respectively. Moreover, define \begin{align*}p_0( t|w)&= \P(T\ge t, D=0|W =w);\\  
p_1( t,z|w)&= \lim_{dz\to 0}\frac{1}{dz}\P(T\ge t,Z\in[z,z+dz], D=1 |W=w).\end{align*} Then, we have the following corollary
\begin{corollary}
If $h_0(t|w),h_1(t|z,w),p_1(t,z|w)$ exist for all $t\in\R_+,z\in\mathcal{Z},w\in\mathcal{W}$ and the distribution of $W$ is continuous, equation \eqref{main_eq} in the main text is equivalent to 
$$\int_0^{\varphi_0(u)} h_0(t|w) p_0(t|w)dt +\int\int_0^{\varphi_1(z,u)} h_1(t|z, w)p_1(t,z|w)dtdz =1-e^{-u}.$$
\end{corollary}
\begin{Proof}
Differentiating both sides of equation (2) in the main text with respect to $w$ and dividing by $dF_W(w)/dw$, we get that equation \eqref{main_eq} is equivalent to
\begin{equation}\label{reec}\P(T\le \varphi_0(u),D=0|W=w) +\P(T\le \varphi_1(Z,u),D=1|W=w)=1-e^{-u},\end{equation}
for all $w\in \mathcal{W}$. 
Next, by Bayes' theorem, we have $$p_0(t|w)= \P(T\ge t|D=0, W= w)  \P(D=0|W =w).$$ As a result,
\begin{align*}
h_0(t|w) p_0(t|w)
&= \left[\lim_{dt\to 0}\frac{1}{dt}\P\left( T\in[t,t+dt]| D=0, W=w\right)\right] \P(D=0|W =w)\\
&= \lim_{dt\to 0}\frac{1}{dt}\P\left( T\in[t,t+dt],D=0| W=w\right).
\end{align*}
Therefore, $\int_0^{\varphi_0(u)} h_0(t|w) p_0(t|w)dt =\P(T\le \varphi_0(u),D=0|W=w) $. Similarly, we have
\begin{align*}p_1( t,z|w)&= \lim_{dz\to 0}\frac{1}{dz}\P(T\ge t,Z\in[z,z+dz], D=1 |W=w)\\
&=\P\left( T\in[t,t+dt]| Z=z, D=1,W= w\right)\left[\lim_{dz\to 0}\frac{1}{dz}\P( Z\in[z,z+dz],D=1|W=w)\right].\end{align*}
Hence, we get 
\begin{align*}
&h_1(t|z, w)p_1(t,z|w)\\
&= \left[\lim_{dt\to 0}\frac{1}{dt}\P\left( T\in[t,t+dt]| Z=z,D=1, W=w\right)\right]\left[\lim_{dz\to 0}\frac{1}{dz}\P( Z\in[z,z+dz],D=1|W=w)\right]\\
&= \lim_{dt\to 0} \lim_{dz\to 0}\frac{1}{dtdz}\P\left( T\in[t,t+dt],Z\in[z,z+dz],D=1| W=w\right),
\end{align*}
which yields $\int\int_0^{\varphi_1(z,u)} h_1(t|z, w)p_1(t,z|w)dtdz =\P(T\le \varphi(Z,u),D=1|W=w)$. So, from \eqref{reec}, equation \eqref{main_eq} in the main text is equivalent to the equation in Corollary \ref{reec}.
\end{Proof}
\subsection{On Theorem \ref{th.loc_ident}}\label{prooflocident}
\subsubsection{Theoretical framework and regularity conditions}
In this subsection, we first outline a theoretical framework embedding the proof of Theorem \ref{th.loc_ident} and then state the regularity conditions that we make to obtain Theorem \ref{th.loc_ident}.

Let $f_Z$ be the density of $Z$, $f_{Z|U}(\cdot|u)$ be the density of $Z$ given $U=u$, $f_{Z|W}(\cdot|w)$ be the density of $Z$ given $W=w$ and $f_{T|Z,W}(\cdot|z,w)$ be the density of $T$ given $Z=z,W=w$. These densities exist under our assumptions (see Assumption \ref{asss1} (i) and Lemma \ref{justifl} below). Let $L^2(W)=\{f:\R\to\R:\ E[f(W)^2]<\infty\}$, which we endow with the norm $\|f\|_{L^2(W)}^2=E[f(W)^2]$. For $\psi_u=(\psi_0,\psi_1)\in \mathcal P$, such that equation \ref{restatemaineq} in the main text holds, we have
\begin{align*}&\int_0^w \int \int  I(t\le \psi_0)  I(z>t)   f_{T|Z,W}(t|z,w)f_{Z|W}(z|\omega) dzdtdF_{W}(\omega)\\
&\quad + \int_0^w \int \int   I(t\le \psi_1(z)) I(z\le t)  f_{T|Z,W}(t|z,w) f_{Z|W}(z|\omega)dzdtdF_{W}(\omega)=(1-e^{-u})F_W(w),
\end{align*}
for all $w\in\mathcal{W}$. (This is just a reformulation of equation (3) in the main text in terms of densities.) Differentiating both sides of the equation with respect to $w$ and dividing by $dF_W(w)/dw$, we obtain 
$$B_0(\psi_0 |w) +\int  B_1(\psi_1(z),z| w)dz =(1-e^{-u}),$$
where 
\begin{align*}
B_0(t|w)&= \int_0^t \int I(z>s) f_{T|Z,W}(s|z,w)f_{Z|W}(z|w)dzds;\\
B_1(t,z|w)&=\int_0^t  I(z\le s) f_{T|Z,W}(s|z,w) f_{Z|W}(z|w)ds.
\end{align*}
Now, we can define an operator $A$ from $(\mathcal{P},\|\cdot\|_{\mathcal{P}})$ to $(L^2(W),\|\cdot\|_{L^2(W)})$ defined by
$$ A(\psi_u)(w) =B_0(\psi_0 |w) +\int  B_1(\psi_1(z),z| w)dz -(1-e^{-u}),$$
for all $w\in\mathcal{W}$ and $\psi_u=(\psi_0,\psi_1)\in \mathcal P$. By the above reasoning, we have $A(\psi_u)(w)=0$, for all $w\in\mathcal{W}$ for all $\psi_u=(\psi_0,\psi_1)\in \mathcal P$ satisfying equation \eqref{restatemaineq} in the main text. The system of equations $A(\psi_u)(w)=0$, for all $w\in\mathcal{W}$ is just a conditional version of \eqref{restatemaineq} in the main text. Remark that, by the fundamental theorem of calculus, $ B_0(\cdot|w)$ is differentiable with derivative
$$\frac{\partial B_0}{\partial  t}(t|w)= \int I(z>t) f_{T|Z,W}(t|z,w)f_{Z|W}(z|w)dz.$$ 
Similarly, for all $t\ne z$, $ B_1(\cdot,z|w)$ is differentiable at $t$ with derivative
$$\frac{\partial B_1}{\partial  t}(t,z|w)= I(z\le t) f_{T|Z,W}(t|z,w) f_{Z|W}(z|w).$$ 
Let us define 
\begin{align*}
A^\prime_{\varphi_u}(\psi_u)(w) 
&= \psi_0  \frac{\partial B_0}{\partial t}(\varphi_0(u)|w)+\int \psi_1(z) \frac{\partial B_1}{\partial t}(\varphi_1(z,u),z|w)dz.
\end{align*}
The mapping $A^\prime(\varphi_u)(\psi_u)$ is the Gâteaux derivative of $A$ at $\varphi_u$ in the direction $\psi_u$ when this Gâteaux derivative exists.
We make the following Assumption:
\begin{assumptionp}{S.1}\label{asss1}
The following holds:
\begin{itemize}
\item[(i)] The distribution of $(U,Z,W)$ is absolutely continuous with continuous density;
\item[(ii)] For every $\psi_u=(\psi_0,\psi_1)\in \mathcal P$, the operator $A$ is Gâteaux differentiable at $\varphi_u$ in the direction $\psi_u$ with Gâteaux derivative equal to $A^\prime_{\varphi_u}(\psi_u)$.
\item[(iii)] There exists constant $K_1>0$ such that $f_Z(z)/f_{Z|U}(z|u)\le K_1$ for all $z\in\mathcal{Z}$;
\item[(iv)] $P(Z>\varphi_0(u)|U=u)>0$.
\end{itemize}
\end{assumptionp}
\noindent Conditions (i), (ii), (iii) are regularity conditions. Condition (iv) says that there are some subjects with $U=u$ who are not treated (before the end of their spell).

\subsubsection{Proof of Theorem \ref{th.loc_ident}}
We gather here the body of the proof. It relies on Lemma \ref{l2loc} proved in Section \ref{auxlemloc}.  Take now a fixed $\psi_u=(\psi_0,\psi_1)\in\mathcal{P}$ such that equation \eqref{restatemaineq} in the main text holds. \\

The proof proceeds in two steps. First, we characterize the kernel of $A^\prime_{\varphi_u}$.
By Lemma \ref{l2loc}, we have, almost surely
\begin{align}\notag &A^\prime_{\varphi_u}(\psi_u)(W) \\
\notag &=e^{-u}\frac{\psi_0}{\varphi_0'(u)}   \P(Z>\varphi_0(u)|U=u,W)\\
\notag &\quad +\int e^{-u}\frac{\psi_1(z)}{\frac{\partial\varphi_1}{\partial u}(z,u)}I(z\le \varphi_0(u))f_{Z|U,W}(z |u,W)dz\\
\label{gateaux3} &= e^{-u}E\left[\left. I(Z>\varphi_0(u))\frac{\psi_0}{\varphi_0'(u)}+ I(Z\le\varphi_0(u))\frac{\psi_1(Z)}{\frac{\partial\varphi_1}{\partial u}(z,u)}\right|W, U=u\right], \end{align}
where $f_{Z|U,W}(\cdot|u,w)$ is the density of $Z$ given $U=u,W=w$.\\

Second, we consider two cases depending on $ \left\| A^\prime_{\varphi_u}(\psi_u)\right\|_{L^2(W)}$:\\

\noindent\textit{Case 1: $ \left\| A^\prime_{\varphi_u}(\psi_u)\right\|_{L^2(W)}>0$.} In this case, by the definition of the Gâteaux derivative, there exists $\epsilon>0$, such that, for all $\delta\in[-\epsilon,\epsilon]$, it holds that 
\begin{equation}\label{gateaux4} \left\|\frac{A(\varphi_u+\delta \psi_u) - A(\varphi_u)}{\delta}- A^\prime_{\varphi_u}(\psi_u) \right\|_{L^2(W)}\le \frac12 \left\| A^\prime_{\varphi_u}(\psi_u)\right\|_{L^2(W)}.\end{equation}
Let us now show that $\varphi$ is locally identified in $\mathcal{N}_{\epsilon}$. We reason by contradiction. If there existed $\delta \in [-\epsilon,\epsilon]$ such that $A(\varphi_u+\delta \psi_u) =A(\varphi_u)$, then equation \eqref{gateaux4} would yield
$$ \left\| A^\prime_{\varphi_u}(\psi_u) \right\|_{L^2(W)}\le \frac12 \left\| A^\prime_{\varphi_u}(\psi_u)\right\|_{L^2(W)},$$
which is false. Hence, there are no $\delta \in [-\epsilon,\epsilon]$ such that $A(\varphi_u+\delta \psi_u) =A(\varphi_u)$, which shows local identification of $\varphi(\cdot,u)$ in $$\mathcal{N}_\epsilon=\{(\varphi_0(u)+\delta \psi_0,\varphi_1(\cdot,u)+\delta \psi_1):\ \delta\in [-\epsilon,\epsilon]\}.$$

\noindent\tr{\textit{Case 2: $ \left\| A^\prime_{\varphi_u}(\psi_u)\right\|_{L^2(W)}=0$.} By equation \eqref{gateaux3} and Assumption \ref{Rglob}, $\left\| A^\prime_{\varphi_u}(\psi_u)\right\|_{L^2(W)}=0$  implies that 
\begin{equation}\label{case21}\P\left(\left.I(Z>\varphi_0(u))\frac{\psi_0}{\varphi_0'(u)}+ I(Z\le\varphi_0(u))\frac{\psi_1(Z)}{\frac{\partial\varphi_1}{\partial u}(z,u)} =0 \right|U=u\right)=1.\end{equation}
Since $\P(Z>\varphi_0(u)|U=u)>0$ (Assumption \ref{asss1} (iv)), equation \eqref{case21} yields that $\psi_0=0$, which also implies $\P( I(Z\le\varphi_0(u)) \psi_1(Z)=0|U=u)=1$. From this, we get 
\begin{align*}
\P( I(Z\le\varphi_0(u)) \psi_1(Z)\ne 0)&=\int I(z\le\varphi_0(u)) I(\psi_1(z)\ne 0)f_Z(z)dz\\
&= \int I(z\le\varphi_0(u)) I(\psi_1(z)\ne 0)f_{Z|U}(z|u)\frac{f_Z(z)}{f_{Z|U}(z|u)}dz\\
&\le K_1\int I(z\le\varphi_0(u)) I(\psi_1(z)\ne 0)f_{Z|U}(z|u)dz\\
&= K_1\P(I(Z\le\varphi_0(u)) \psi_1(Z)\ne 0|U=u)=0,
\end{align*}
where the inequality follows from Assumption \ref{asss1} (iii). This implies $I(Z\le\varphi_0(u)) \psi_1(Z)=0\ a.s.$ As a result, for all $\delta \in \R$, we have 
\begin{align*}
&(\varphi_0(u)+\delta\psi_0)I(Z> (\varphi_0(u)+\delta\psi_0))+ (\varphi_1(Z,u)+\delta\psi_1(Z))I(Z\le (\varphi_0(u)+\delta\psi_0))\\
&= \varphi_0(u)I(Z> \varphi_0(u))+ (\varphi_1(Z,u)+\delta\psi_1(Z))I(Z\le \varphi_0(u))\\
&=\varphi_0(u)I(Z> \varphi_0(u))+ \varphi_1(Z,u)I(Z\le \varphi_0(u))\\
&=\varphi(Z,u)\ a.s.
\end{align*}
Hence for all $\delta\in\R$,  the ``induced" mapping of $(\varphi_0(u)+\delta \psi_0,\varphi_1(\cdot,u)+\delta \psi_1)$ is equal to $\varphi(Z,u)$ almost surely. This shows local identification of $\varphi(\cdot,u)$ in $\mathcal{N}_\epsilon$ for any $\epsilon >0$ (all elements in $$\mathcal{N}_\epsilon=\{((\varphi_0(u)+\delta \psi_0,\varphi_1(\cdot,u)+\delta \psi_1):\ \delta\in [-\epsilon,\epsilon]\}$$ have an ``induced" mapping equal to $\varphi(Z,u)$ a.s.).}

\subsubsection{Auxiliary lemmas concerning Theorem \ref{th.loc_ident}}\label{auxlemloc}
The next lemma justifies the existence of $f_{T|Z,W}$ under our assumptions.
\begin{Lemma}\label{justifl} Under the assumptions of Theorem \ref{th.loc_ident}, we have 
$$f_{T|Z,W}(t|z,w)=\lambda(z,t) f_{U|Z,W}\left(\left.\int_0^t\lambda(z,s)ds\right|z,w\right),$$
where $f_{U|Z,W}(\cdot|z,w)$ is the density of $U$ given $Z=z$, $W=w$.
\end{Lemma}
\begin{Proof}
\tr{First note that 
\begin{align*}
&\P\left(T\le t|Z=z,W=w\right)\\
&=\P\left(\varphi(z,U)\le t|Z=z,W=w\right)\\
&=\P\left(U\le (\varphi(z,\cdot))^{-1}(t)|Z=z,W=w\right)\\
&=\P\left(\left. U\le \int_0^t\lambda(z,s)ds\right|Z=z,W=w\right).
\end{align*}
Hence, by the chain rule, the function $t\mapsto \P\left(T\le t|Z=z,W=w\right)$ is differentiable with derivative
$$f_{T|Z,W}(t|z,w)=\lambda(z,t) f_{U|Z,W}\left(\left.\int_0^t\lambda(z,s)ds\right|z,w\right).$$}
\end{Proof}
The following lemma helps to characterize the kernel of $A_{\varphi_u}$.
\begin{Lemma}\label{l2loc} Under the assumptions of Theorem 3.2, the following holds:
\begin{itemize}
\item[(i)] For all $w\in\mathcal{W}$, $ \varphi_0'(u)\frac{\partial B_0}{\partial t}(\varphi_0(u)|w)= e^{-u}\P(Z>\varphi_0(u)|U=u,W=w).$
\item[(ii)] For all $w\in\mathcal{W}$, $$\P\left(\left.\frac{\partial \varphi_1}{\partial u}(Z,u)\frac{\partial B_1}{\partial t}(\varphi_1(Z,u),z|w)= e^{-u}I(Z\le \varphi_0(u))f_{Z|U,W}(Z |u,w)\right|W=w\right)=1,$$
where $f_{Z|U,W}(\cdot|u,w)$ is the density of $Z$ given $U=u,W=w$.
\end{itemize}
\end{Lemma}

\begin{Proof} 
 \tr{First, we show (i). For all $t\in\R,w\in\mathcal{W}$, we have 
\begin{align*}B_0(t|w)&= \P(T\le t, D=0|W=w)\\
&=  \P(T\le t, Z>T|W=w)\\
&= \P(\varphi_0(U)\le t, Z>\varphi_0(U)|W=w)\\
&= \P(U\le \varphi_0^{-1}(t), Z>\varphi_0(U)|W=w)\\
&=\int_0^{\varphi_0^{-1}(t)}\P(Z>\varphi_0(u)|U=u,W=w)e^{-u}du,
\end{align*}
where in the second equality we used the fact that the events $\{Z>T\}$ and $\{Z>\varphi_0(U)\}$ are equivalent (see the proof of Theorem \ref{th.main_eq}), the third equality follows from equation \eqref{main_eq} in the main text and in the last line we used Bayes' theorem and the fact that $U|W=w\sim \text{Exp}(1)$. Next, recall that $\varphi_0^{-1}(t)=\int_0^{t}\lambda(\infty,s)ds$. Hence, since $\lambda(\infty,\cdot)$ is continuous by Assumption \ref{RI} (ii), by the inverse function theorem, $\varphi_{0}^{-1}$ is differentiable with derivative equal to $\lambda(\infty,\cdot)$. By the chain rule, this yields that 
$$\frac{\partial B_0}{\partial t}(t|w)=\lambda(\infty,t)\P(Z>\varphi_0(u)|U=u,W=w)e^{-u}.$$
We get $$ \frac{\partial B_0}{\partial t}(\varphi_0(u)|w)= \lambda(\infty,\varphi_0(u))\P(Z>\varphi_0(u)|U=u,W=w)e^{-u}.$$
To conclude, notice that $\int_0^{\varphi_0(u)}\lambda(\infty,t)dt =\varphi_{0}^{-1}(\varphi_0(u))=u$. Hence, differentiating both sides of $\int_0^{\varphi_0(u)}\lambda(\infty,t)dt =u$, by the chain rule, we have $\varphi_0'(u)\lambda(\infty,\varphi_0(u))=1$, which yields the result.}\\

Next, we prove (ii). Let us consider $z\ne \varphi_0(u)$. Remark that we have
\begin{align*}B_1(\varphi_1(z,u),z|w)&=I(z\le \varphi_1(z,u)) \P(T\le \varphi_1(z,u)|Z=z,W=w)f_{Z|W}(z|w)\\
&=I(z\le \varphi_1(z,u)) \P(\varphi_1(z,U)\le \varphi_1(z,u)|Z=z,W=w)f_{Z|W}(z|w)\\
&= I(z\le \varphi_1(z,u)) \P(U\le u|Z=z,W=w)f_{Z|W}(z|w).
\end{align*}
Next, recall that $\varphi_1(z,\cdot)$ is strictly increasing and $\varphi_1(z,\varphi_0^{-1}(z))=z$. Hence, we have
$I(z\le \varphi_1(z,u))=I(\varphi_1(z,\varphi_0^{-1}(z))\le \varphi_1(z,u))=I(\varphi_0^{-1}(z)\le u)$.  This leads to 
$$B_1(\varphi_1(z,u),z|w)=I(z\le \varphi_0(u))  \P(U\le u|Z=z,W=w)f_{Z|W}(z|w) $$
Next remark that, since $\varphi_1(z,\varphi_0^{-1}(z))=z$, $\varphi_1$ is strictly increasing and we took $z\ne \varphi_0(u)$, we have $\varphi_1(z,u)\ne z$. By Assumption 2.2 (ii) and the fundamental theorem of calculus, this implies that $\Lambda(z,\cdot)$ is continuously differentiable in a neighborhood around $\varphi_1(z,u)$. Hence, by the inverse function theorem, $\varphi_1(z,\cdot)=\Lambda(z,\cdot)^{-1}$ is differentiable at $u$ with derivative denoted by $\frac{\partial \varphi_1}{\partial u}(z,u)$. Note also that, since $\varphi_0$ is continuous, $I(z\le \varphi_0(\cdot))$ is constant in a neighborhood around $u$ (and therefore differentiable). The above discussion and the chain rule imply that $B_1(\varphi_1(z,\cdot),z|w)=I(z\le \varphi_0(\cdot))  \P(U\le \cdot|Z=z,W=w)f_{Z|W}(z|w)$ is differentiable at $u$ with derivative
 $$\frac{\partial \varphi_1}{\partial u}(z,u)\frac{\partial B_1}{\partial t}(\varphi_1(z,u),z|w)=  I(z\le \varphi_0(u)) f_{U|Z,W}(u|z,w)f_{Z|W}(z|w)$$ almost everywhere. Using Bayes' theorem and the fact that $U|W=w\sim \text{Exp}(1)$, we get $f_{U|Z,W}(u|z,w)f_{Z|W}(z|w)=e^{-u}f_{Z|U,W}(z |u,w) $, we 
get $$\frac{\partial \varphi_1}{\partial u}(z,u)\frac{\partial B_1}{\partial t}(\varphi_1(z,u),z|w)=  e^{-u}I(z\le \varphi_0(u))f_{Z|U,W}(z|u,w).$$
To conclude just notice that $\P(Z=\varphi_0(U)|W)=0\ \text{a.s.}$ since the distribution of $(U,Z,W)$ is continuous by Assumption \ref{asss1}.
\end{Proof}

\subsection{On Theorem \ref{globident}} In Section \ref{secexpr}, we justify the existence of $f_{\epsilon|D,W}(\cdot|0,w)$ and $f_{\epsilon|D,Z,W}(\cdot|1,z,w)$, in Section \ref{proofglobid}, we prove Theorem \ref{globident} and, in Section\ref{suffglobident}, we give sufficient conditions for Assumption \ref{as.global}.

\subsubsection{Expressions for $f_{\epsilon|D,W}(\cdot|0,w)$ and $f_{\epsilon|D,Z,W}(\cdot|1,z,w)$} \label{secexpr}
In this section, we justify the existence of $f_{\epsilon|D,W}(\cdot|0,w)$ and $f_{\epsilon|D,Z,W}(\cdot|1,z,w)$ under our assumptions by expressing them as functions of $f_{T|Z,W}$ and $f_{Z|W}$ (the existence of $f_{T|Z,W}$ is justified by Lemma \ref{justifl}).
We have the following expressions for $f_{\epsilon|D,W}(\cdot|0,w)$ and $f_{\epsilon|D,Z,W}(\cdot|1,z,w)$:
\begin{align*}
f_{\epsilon|D,W}(e|0,w)&=\int I(z>e+\varphi_0(u))f_{T|Z,W}(e+\varphi_0(u)|z,w)f_{Z|W}(z|w)dz;\\
f_{\epsilon|D,Z,W}(e|1,z,w)&=I(z\le e+\varphi_1(z,u))f_{T|Z,W}(e+\varphi_1(z,u)|z,w)f_{Z|W}(z|w)dz.
\end{align*}
\subsubsection{Proof of Theorem \ref{globident}}\label{proofglobid}
Let $(\psi_0,\psi_1)\in\mathcal{P}$ be such that $$F_0(\psi_0,w) +F_1(\psi_1,w)=(1-e^{-u})F_W(w),\text{for all}\ w\in\mathcal{W}.$$ This and Theorem 3.1 yield , 
\begin{align}
\notag &E[I(T\le\varphi_0(u),D=0,W\le w)- I(T\le\psi_0,D=0,W\le w)]\\
\notag &\quad+E[I(T\le \varphi_1(Z,u),D=1, W\le w)-I(T\le \psi_1(Z), D=1,W\le w)]=0\\
\notag &\Leftrightarrow \int_0^w E[I(T\le\varphi_0(u))- I(T\le\psi_0)|D=0,W=\omega ]\P(D=0|W=\omega)dF_W(\omega)\\
\notag &\quad + \int_0^w E[E[I(T\le \varphi_1(Z,u))-I(T\le \psi_1(Z))|D=1,W=\omega]]\P(D=1|W=\omega)dF_W(\omega) =0.
\end{align}
Differentiating with respect to to $w$ and dividing by $dF_W(w)/dw$ both sides of the previous equation, we obtain
\begin{align}\notag &E[I(T\le\varphi_0(u))- I(T\le\psi_0)|D=0,W=w ]\P(D=0|W=w)\\
\label{reecriture}&\quad + E[I(T\le \varphi_1(Z,u))-I(T\le \psi_1(Z))|D=1,W=w] \P(D=1|W=w)=0
\end{align}
Let now $\Delta_0=\psi_0-\varphi_0(u)$ and $\Delta_1=\psi_1-\varphi_1(\cdot,u)$.
We have
\begin{align}
\notag &E[I(T\le\varphi_0(u))- I(T\le\psi_0)|D=0,W=w]\\
\notag &=E[I(\epsilon\le 0)- I(\epsilon\le\Delta_0)|D=0,W=w]\\
\label{reec1}&=\int_0^1\Delta_0f_{\epsilon|D,W}(\delta\Delta_0|0,w)d\delta
\end{align}
and
\begin{align}
\notag &E[I(T\le \varphi_1(Z,u))-I(T\le \psi_1(Z,u))|Z,D=1,W=w ]\\
\notag &=E[I(\epsilon\le 0)-I(\epsilon\le  \Delta_1(Z))|Z,D=1,W=w ]\\
\label{reec2}&=E\left[\left.\int_0^1\Delta_1(Z)f_{\epsilon|D,Z,W}(\delta\Delta_1(Z)|1,Z,W)d\delta\right|D=1,W=w \right].
\end{align}
Together \eqref{reecriture}, \eqref{reec1} and \eqref{reec2} yield that we have
\begin{align*}
&E\biggr[(1-D)\Delta_0\int_0^1f_{\epsilon|D,W}(\delta\Delta_0|0,w)d\delta\\
&\quad +D\Delta_1(Z)\int_0^1f_{\epsilon|D,Z,W}(\delta\Delta_1(Z)|1,Z,W)d\delta\biggr|W=w\biggr]=0\\
& \Leftrightarrow E[(\Delta_0(1-D)+\Delta_1(Z)D) \omega_\Delta(Z,D,W)|W=w]=0.
\end{align*}
This implies
\begin{equation}\label{asnov}\Delta_0(1-D)+\Delta_1(Z)D=0\ \text{a.s.},\end{equation}
by Assumption \ref{as.global} Since $\P(D=0)>0$ (Assumption \ref{Rglob}), we obtain $\Delta_0=0$. This and \eqref{asnov} yield
$\P(\Delta_1(Z)D=0)=1,$ which yields $\P(\Delta_1(Z)\ne 0|D=1)=0$ by Bayes' rule and the fact that $\P(D=1)>0$ (Assumption \ref{Rglob}).
 From this we get 
\begin{align*}
\P( I(Z\le\varphi_0(u)) \Delta_1(Z)\ne 0)
&=\int I(z\le\varphi_0(u)) I(\Delta_1(z)\ne 0)f_Z(z)dz\\
&= \int  I(\Delta_1(z)\ne 0)f_{Z|D}(z|1)\frac{f_Z(z)I(z\le\varphi_0(u))}{f_{Z|D}(z|1)}dz\\
&\le K\int  I(\Delta_1(z)\ne 0)f_{Z|D}(z|1)dz\\
&\le K\P( \Delta_1(Z)\ne 0|D=1)=0,
\end{align*}
where the first inequality is due to Assumption \ref{Rglob}.
As a result, we have 
\begin{align*}
\psi(Z,u)&=\psi_0I(Z>\psi_0)+ \psi_1(Z)I(Z\le \psi_0)\\
&= (\varphi_0(u)+\Delta_0)I(Z> \varphi_0(u))+ (\varphi_1(Z,u)+\Delta_1(Z))I(Z\le  \varphi_0(u))\\
&=\varphi_0(u)I(Z> \varphi_0(u))+ \varphi_1(Z,u)I(Z\le \varphi_0(u))\\
&=\varphi(Z,u)\ a.s.
\end{align*}

\subsubsection{Sufficient conditions for Assumption \ref{as.global}}\label{suffglobident}
Call $\{\Omega,\mathcal{F},\P\}$ the probability space underlying our analysis.  Let $\P_{Z,D|W}(\cdot|w)$ be the probability measure associated with the distribution of $(Z,D)$ given $W=w$, that is for any Borel set $B$ in $\R^2$, $\P_{Z,D|W}(B|w)=\P((Z,D)\in B|W=w)$. For $\Delta=(\Delta_0,\Delta_1)\in\mathcal{P}$, let us define the following probability distributions:
$$\P_{Z,D|W}^{\Delta}(B|w)=K_{w,\Delta}\int_{(z,d)\in B}\omega_{\Delta}(z,d,w)d\P_{Z,D|W}(z,d|w),$$
 to every Borel set $B$ in $\R^2$, where 
$$K_{w,\Delta}=\left[\int\omega_{\Delta}(z,d,w)d\P_{Z,D|W}(z,d|w)\right]^{-1}$$
guarantees that $\P_{Z,D|W}^{\Delta}(\cdot|w)$ is indeed a probability distribution. We assume below that the ``constant" $K_{w,\Delta}$ exists for relevant $\Delta$. We can also define $\P_{Z,D}^{\Delta}$, the unconditional version of $\P_{Z,D|W}^{\Delta}(\cdot|w)$ defined as
$$\P_{Z,D}^{\Delta}(B)= \int \P_{Z,D|W}^{\Delta}(B|w)dF_W(w) $$
for all Borel set $B$ in $\R^2$.
We have the following Lemma which relates Assumption \ref{as.global} to bounded completeness conditions.
\begin{Lemma}
Assume that Assumption \ref{Rglob} in the main text holds. Suppose also that for all $(\psi_0,\psi_1)\in\mathcal{P}$ and $\Delta=(\varphi_0(u)-\psi_0,\varphi_1(\cdot,u)-\psi_1)\ne 0$, the following holds
\begin{itemize}
\item[(i)] $K_{w,\Delta}$ exists.
\item[(ii)] The family of distributions $\{\P_{Z,D|W}^{\Delta}(\cdot|w)\}_{w\in\mathcal{W}}$ is boundedly complete, that is for all bounded $f:\R_+\times\{0,1\}\mapsto \R$,
\begin{align*}&\int f(z,d) d\P_{Z,D|W}^{\Delta}(z,d|W)=0,\ a.s.\\
&\Rightarrow  \P_{Z,D}^{\Delta}(f(Z,D)=0)=1.\end{align*}
\item[(iii)] For every Borel set $B$ in $\R^2$, $\P_{Z,D}^{\Delta}(B)=0\Rightarrow \P_{Z,D}(B):=\P((Z,D)\in B)=0$. 
\end{itemize}Then, Assumption \ref{as.global} holds.
\end{Lemma}
\begin{Proof}
Take $(\psi_0,\psi_1)\in\mathcal{P}$ and define $\Delta=(\varphi_0(u)-\psi_0,\varphi_1(\cdot,u)-\psi_1).$ Assume that
$$E[(\Delta_0(1-D)+\Delta_1(Z)D) \omega_\Delta(Z,D,W)|W]=0\ a.s.$$
This implies 
$$\int (\Delta_0(1-d)+\Delta_1(z)d) d\P_{Z,D|W}^{\Delta}(z,d|W)=0,\ a.s.$$
This yields $\P_{Z,D}^{\Delta}(\Delta_0(1-D)+\Delta_1(Z)D=0)$ by condition (ii). Since, by condition (iii), for any Borel set $B$ in $\R^2$, we have $\P_{Z,D}^{\Delta}(B)=0\Rightarrow \P_{Z,D}(B)=0$, this implies that $\P(\Delta_0(1-D)+\Delta_1(Z)D\ne 0)=\P_{Z,D}^{\Delta}(\Delta_0(1-D)+\Delta_1(Z)D\ne 0)=0$. 
\end{Proof}
\subsection{Proof of equations \eqref{F0cens} and \eqref{F1cens} in the main text } \label{proofcens}
We only show \eqref{F1cens} since the proof of \eqref{F0cens} is similar. We have, for all $\phi:\mathcal{Z}\mapsto [0,c_0)$ and $w\in\mathcal{W}$, 
$$ E\left[\frac{\delta}{G(Y)}I( Y\le \phi(\tilde Z),\tilde D=1,W\le w)\right]= E\left[\frac{\delta}{G(Y)}I( T\le \phi(Z),D=1,W\le w)\right]$$
since, when $\delta=1$, it holds that $Y=T$ and $\tilde D=D$, and, when $D=1$, $\tilde Z=Z$. This yields
\begin{align*}
  &E\left[\frac{\delta}{G(T)}I( T\le \phi(Z),D=1,W\le w)\right]\\
 &=E\left[ E\left[\left.\frac{\delta}{G(T)}I( T\le \phi(Z),D=1,W\le w)\right|T,Z,D,W\right]\right]\\
  &=E\left[ E\left[\left.\frac{\delta}{G(T)}\right|T,Z,D,W\right]I( T\le \phi(Z),D=1,W\le w)\right]\\
  &=E\left[ E\left[\left.\frac{1}{G(T)}\right|T,Z,D,W,\delta=1\right]\P(\delta=1|T,Z,D,W)I( T\le \phi(Z),D=1,W\le w)\right]\\
  &=E\left[ E\left[\left.\frac{1}{G(T)}\right|T,Z,D,W,\delta=1\right]\P( C\ge T|T,Z,D,W)I( T\le \phi(Z),D=1,W\le w)\right]\\
  &=E\left[ \frac{G(T)}{G(T)}I( T\le \phi(Z),D=1,W\le w)\right]\\
  &=E\left[ I( T\le \phi(Z),D=1,W\le w)\right]=F_1(\psi,w),
\end{align*}
where we used $\P( C\ge T|T,Z,D,W)=G(T)$ which is a consequence of Assumption \ref{as.censoring}.
\section{Proof of Theorem 4.1} \label{S.1}

The body of the proof of Theorem \ref{AN} is given in Section \ref{subsec.pAN}. It relies on technical lemmas proved in Sections \ref{subsec.UC} and \ref{subsec.WC}. For $\theta \in \Theta$, we use the notation $\bar L(\theta)= (nm)^{-1}\sum_{i=1}^n \sum_{j=1}^mp(u_j) M_\theta(u_j,W_i)^2$.

\subsection{Main proof}\label{subsec.pAN}

\textbf{Step 1: Consistency.}  By Assumption \ref{paramregul} (iii) and the Lebesgue dominated convergence theorem, $L(\cdot)$ is continuous. Moreover, Assumption \ref{paramid} implies that $L$ is uniquely minimized at $\theta_*$. Additionally, by Lemma \ref{convM} (iii), we have \begin{equation}\label{convL} \sup_{\theta\in\Theta}\left|\widehat{L}(\theta)-L(\theta)\right|=o_P(1).\end{equation}
Hence, by Theorem 2.1 in \cite{newey1994large}, it holds that $\widehat{\theta}- \theta=o_P(1)$. \\

\noindent \textbf{Step 2: Asymptotic expansion.} Consider a sequence of $\theta\in\Theta$ such that $\theta\to\theta_*$. Remark that $M_{\theta_*} \equiv 0$ and that 
 \begin{align}
\notag \widehat{L}( \theta)-\bar L(\theta)&= \frac{1}{nm}\sum_{i=1}^n \sum_{j=1}^m p(u_j)\left[\left(\widehat{M}_{ \theta}-M_{ \theta}\right)(u_j,W_i)\right]^2\\
 \label{TMD1}&\quad + \frac{2}{nm}\sum_{i=1}^n \sum_{j=1}^m p(u_j) M_{\theta}(u_j,W_i)\left[\left(\widehat{M}_{\theta}-M_{ \theta}\right)(u_j,W_i)\right].
 \end{align}
Moreover, since $\bar L(\theta_*)=0$, it holds that 
$$\frac{1}{nm}\sum_{i=1}^n \sum_{j=1}^mp(u_j) \left[\left(\widehat{M}_{ \theta}-M_{  \theta}\right)(u_j,W_i)\right]^2=\widehat{L}(\theta_*)+I_1(\theta)+I_2 (\theta),$$ where 
  \begin{align*}
 I_1(\theta)&= \frac{1}{nm}\sum_{i=1}^n \sum_{j=1}^mp(u_j) \left[\left(\widehat{M}_{ \theta}-M_{ \theta}-\widehat{M}_{\theta_*}\right)(u_j,W_i) \right]^2, \\
I_2 (\theta)&= \frac{2}{nm}\sum_{i=1}^n \sum_{j=1}^mp(u_j) \widehat{M}_{\theta_*}(u_j,W_i)\left[\left(\widehat{M}_{  \theta}-M_{ \theta}-\widehat{M}_{\theta_*}\right)(u_j,W_i) \right].
 \end{align*}
By Lemma \ref{convM} (ii) we have $|I_1(\theta)|=o_P(n^{-1})$. Then, by the inequality of Cauchy-Schwarz, it holds that 
$$I_2(\theta)\le 2 \sqrt{\frac{1}{nm}\sum_{i=1}^n \sum_{j=1}^mp(u_j)^2 \widehat{M}_{\theta_*}(u_j,W_i)^2}\sqrt{\frac{1}{nm}\sum_{i=1}^n \sum_{j=1}^m\left(\widehat{M}_{  \theta}-M_{ \theta}-\widehat{M}_{\theta_*}\right)(u_j,W_i)^2 },$$
which is $o_P(n^{-1})$ by Lemma \ref{convM} (i) and (ii).
As a result, we obtain
\begin{equation}\label{TMD2}\frac{1}{nm}\sum_{i=1}^n \sum_{j=1}^m p(u_j)\left[\left(\widehat{M}_{ \theta}-M_{  \theta}\right)(u_j,W_i)\right]^2= \widehat{L}(\theta_*)+o_P(n^{-1}).
\end{equation}

Using a Taylor expansion of $\theta\mapsto M_{\theta}(u_j,W_i)$ around $\theta_*$ and the uniform boundedness of the first two derivatives of $\theta\mapsto M_{\theta}(u,w)$, we get 
\begin{align}
 \notag &\frac{1}{nm}\sum_{i=1}^n \sum_{j=1}^m p(u_j)M_{{\theta}}(u_j,W_i)\left[\left(\widehat{M}_{{\theta}}-M_{{\theta}}\right)(u_j,W_i)\right]\\
 \notag &=( \theta-\theta_*)^\top \frac{1}{nm}\sum_{i=1}^n \sum_{j=1}^m p(u_j)\frac{\partial M_{{\theta_*}}}{\partial \theta}(u_j,W_i)\left[\left(\widehat{M}_{{\theta}}-M_{{\theta}}\right)(u_j,W_i)\right]+o_P(||{\theta}-\theta_*||^2)\\
\label{TMD3} &=\frac12( \theta-\theta_*)^\top n^{-1/2}\widehat{V}+o_P(n^{-1/2}||\theta-\theta_*||+||{\theta}-\theta_*||^2),
 \end{align}
 where the last equality follows from Lemma \ref{convM} (ii) and $\widehat{V}$ is defined in Lemma \ref{weakM}. 
Equations \eqref{TMD1}, \eqref{TMD2} and \eqref{TMD3} yield
\begin{align}\label{AE0}
\widehat{L}({\theta})-\bar L({\theta}) -  \widehat{L}(\theta_*)&=({\theta}-\theta_*)^\top n^{-1/2} \widehat{V}+ o_P(n^{-1/2}||\theta-\theta_*||+||{\theta}-\theta_*||^2+n^{-1}).
\end{align}
Using a second-order Taylor expansion of $\bar L(\cdot)$ around $\theta_*$, we finally get 
\begin{align}\notag \widehat{L}(\theta)-\widehat{L}(\theta_*)&=(\theta-\theta_*)^\top n^{-1/2}\widehat{V}+ \frac12(\theta-\theta_*)^\top \nabla^2\bar L(\theta_*)(\theta-\theta_*)\\
&\quad +o_P(n^{-1/2}||\theta-\theta_*||+||{\theta}-\theta_*||^2+n^{-1}),\label{AE}\end{align}
where the remainder term is controlled using the uniform boundedness of the second differential of $\theta\mapsto M_{\theta}(u,w)$.\\

 \noindent\textbf{Step 3: Rate of convergence.} Using \eqref{AE} with $\theta=\widehat \theta$ and the fact that $ \widehat{L}( \widehat\theta)-\widehat{L}(\theta_*)\le 0$, we obtain
 \begin{align*}\notag  \frac12(\widehat \theta-\theta_*)^\top \nabla^2\bar L(\theta_*)(\widehat\theta-\theta_*)+  (\widehat\theta-\theta_*)^\top n^{-1/2}\widehat{V}+o_P(n^{-1/2}||\widehat \theta-\theta_*||+||\widehat {\theta}-\theta_*||^2+n^{-1})\le 0. 
\end{align*}
 By Lemma \ref{weakM}, $(\widehat \theta-\theta_*)^\top n^{-1/2}\widehat{V}=O_P(n^{-1/2} ||\widehat \theta-\theta_*||)$. Then, since by the law of large numbers $\nabla^2\bar L(\theta_*)\xrightarrow{\P}\nabla^2L(\theta_*)$ and $\nabla^2L(\theta_*)$ is positive definite, we get 
 $$||\widehat {\theta}-\theta_*||^2= O_P(n^{-1/2}||\widehat \theta-\theta_*||)+o_P(||\widehat {\theta}-\theta_*||^2+n^{-1}),$$
 which yields $||\widehat {\theta}-\theta_*||=O_P(n^{-1/2})$. \\
  
\noindent \textbf{Step 4: Conclusion.} Let $K=\{h=\theta-\theta_*,\ \theta \in \Theta\}$. Equation \eqref{AE} and the fact that $\Theta$ is compact imply that, for any $h\in K$, we have 
\begin{align*}n\left[\widehat{L}(\theta_*+n^{-1/2}h)-\widehat{L}(\theta_*)\right]&=h^\top  \widehat{V} + \frac12 h^\top \nabla^2\bar L(\theta_*)h+R(h),
\end{align*}
where $\sup\limits_{h\in K}|R(h)|=o_P(1)$.
This implies 
\begin{align*}\inf_{h\in K}n\left[\widehat{L}(\theta_*+n^{-1/2}h)-\widehat{L}(\theta_*)\right]&\le \inf_{h\in K}\widehat{U}(h)+o_P(1),
\end{align*}
where  $\widehat U(h)= h^\top \widehat{V}+ (1/2) h^\top \nabla^2\bar L(\theta_*)h$.
Let $\widehat{h}_n=n^{1/2}(\widehat{\theta}-\theta_*)$. The sequence $\widehat{h}_n$ is uniformly tight since $ \sqrt{n}(\widehat{\theta}-\theta_*)=O_P(1)$. We have $$\widehat{h}_n\in\argmin\limits_{h\in K} n\left[\widehat{L}(\theta_*+n^{-1/2}h)-\widehat{L}(\theta_*)\right].$$ 
As a result 
$$\widehat{U}(\widehat{h}_n) \le \inf_{h\in K}\widehat U(h) +o_P(1).$$ Note that by the law of large numbers, $\nabla^2\bar L(\theta_*)\xrightarrow{\P}\nabla^2L(\theta_*)$. Hence, by Lemma \ref{weakM} and Slutsky's theorem, $\widehat{U}(\cdot)$ converges in distribution to $U(\cdot)$, where $U(h)=h^\top V +(1/2)h^\top \nabla^2L(\theta_*)h$ with $V$ defined in Lemma \ref{weakM}. The criterion function $U(h)$ is uniquely maximized at $-\nabla^2L(\theta_*)^{-1}V$, which is tight. Hence, by Theorem 14.1 in \cite{kosorok2008introduction}, $\widehat{h}_n$ converges in distribution to $-\nabla^2L(\theta_*)^{-1}V$, which is a normal random variable with mean zero and variance given by $\nabla^2L(\theta_*)^{-1} \mbox{Var}(V) \nabla^2 L(\theta_*)^{-1}$.

\subsection{Uniform convergence results}\label{subsec.UC}

\begin{Lemma}\label{convF}Let $u\in [0,u_m)$. Under the assumptions of Theorem \ref{AN}, the following holds:
\begin{itemize} 
\item[(i)] $\sup\limits_{w\in\mathcal{W},\theta \in \Theta}\left|(\widehat{F}_0-F_0)(\varphi_{\theta0}(u),w)\right|=O_P(n^{-1/2}).$
\item[(ii)] $\lim\limits_{\epsilon \to 0 }\sup\limits_{w\in\mathcal{W},\ \theta \in\Theta:\ ||\theta-\theta_*||\le \epsilon}n^{1/2}\left|(\widehat{F}_0-F_0)(\varphi_{\theta0}(u),w)-(\widehat{F}_0-F_0)(\varphi_{\theta_*0}(u),w)\right|=o_P(1).$
\item[(iii)] $\sup\limits_{w\in\mathcal{W},\theta \in \Theta}\left|(\widehat{F}_1-F_1)(\varphi_{\theta1}(\cdot,u),w)\right|=O_P(n^{-1/2}).$
\item[(iv)] $\lim\limits_{\epsilon \to 0 }\sup\limits_{w\in\mathcal{W},\ \theta \in\Theta:\ ||\theta-\theta_*||\le \epsilon}n^{1/2}\left|(\widehat{F}_1-F_1)(\varphi_{\theta1}(\cdot,u),w)-(\widehat{F}_1-F_1)(\varphi_{\theta_*1}(\cdot,u),w)\right|=o_P(1).$
\end{itemize}
\end{Lemma}
\begin{Proof} The proof of the first two statements is similar to that of the last two results. Hence, we only prove (iii) and (iv). 

When $\delta_i=1$, it holds that $Y_i=T_i$ and $\tilde D_i=D_i$, and, when $D_i=1$, $\tilde Z_i=Z_i$. Hence, we have 
\begin{align*}
&\widehat{F}_1(\varphi_{\theta1}(\cdot, u),w)-F_1(\varphi_{\theta1}(\cdot,u),w)=\widehat{I}_1(\theta,u,w) - \widehat{I}_2(\theta,u,w),\end{align*}
where
\begin{align*}
\widehat{I}_1(\theta,u,w)&=\frac{1}{n}\sum_{i=1}^n  K_i(\theta,w) -E\left[K(\theta,w)\right]
;\\\widehat{I}_2(\theta,u,w)&=\frac{1}{n}\sum_{i=1}^n  K_i(\theta,w)\frac{\widehat{G}(T_i)-G(T_i)}{\widehat{G}(T_i)},
\end{align*}
with 
\begin{align*}
K(\theta,w) &=\frac{\delta}{G(T)}I(T\le \varphi_{\theta1}(Z,u),D=1, W\le w);\\
K_i(\theta,w) &=\frac{\delta_i}{G(T_i)}I(T_i\le \varphi_{\theta1}(Z_i,u),D_i=1, W_i\le w).
\end{align*}
 Let $t_0=\sup\limits_{\theta\in \Theta,z\in\mathcal{Z}}\varphi_{\theta1}(z,u)$. The class 
\begin{equation}\label{class}\begin{array}{lll}\R_+\times\R_+\times\mathcal{Z}\times \{0,1\}\times \mathcal{W} & \mapsto &\R\\
(t,c,z,d,w) & \mapsto &\frac{I(t\le c)}{G(t\min(t,t_0))}I(t\le \varphi_{\theta 1}(z,u),  d=1,w\le \omega),
\end{array}\end{equation}
$\theta\in\Theta, w\in \mathcal{W}$, is Donsker because it is the product of uniformly bounded Donsker classes (Corollary 9.32  in \cite{kosorok2008introduction}). 
Then, by the Donsker theorem, it holds that
\begin{align}\label{I1i}\sup\limits_{w\in\mathcal{W},\theta \in \Theta}\left|\widehat{I}_1(\theta,u,w)\right|
&=O_P\left(n^{-1/2}\right).\end{align}
 Using Theorem 1.1 in \cite{gill1983large} and the fact that $t_0<c_0$,
\begin{equation}\label{KMU}\sup_{s\in[0,t_0]}\left|\widehat{G}(s) -G(s)\right| =O_P\left(n^{-1/2}\right).\end{equation}
As a result, since $G(t_0)>0$, we get 
$$\sup\limits_{w\in\mathcal{W},\ \theta \in\Theta}\left|\widehat{I}_2(\theta,u,w)\right| =O_P\left(n^{-1/2}\right).$$
This and \eqref{I1i} yield (iii). 

Let us now prove (iv). First, we claim that
\begin{equation}\label{se}\lim_{\epsilon\to 0}\sup_{||\theta-\theta_*||\le \epsilon}\sup\limits_{w\in\mathcal{W}}n^{1/2}\left|\widehat{I}_1(\theta,u,w)-\widehat{I}_1(\theta_*,u,w)\right|=o_P(1).\end{equation}
Since the class \eqref{class} is Donsker, by the analysis of the proof of Lemma 3 in \cite{brown2002weighted}, \eqref{se} is a consequence of
\begin{align}\label{suffse}
\lim_{\epsilon\to 0}\sup_{w\in\mathcal{W}, \theta\in\Theta:\ ||\theta-\theta_*||\le \epsilon} E\left[|K(\theta,w)-K(\theta_*,w)|^2\right]=0,\end{align}
where $K(\theta,w) =\frac{\delta}{G(T)}I(Y\le \varphi_{\theta1}(Z,u),D=1, W\le w)$.  
Notice that \begin{align*}
&E[\{K(\theta,w)-K(\theta_*,w)\}^2]\\
&=E[K(\theta,w)\{K(\theta,w)-K(\theta_*,w)\}] -E[K(\theta_*,w)\{K(\theta,w)-K(\theta_*,w)\}].
\end{align*}
Remark that $\sup_{\theta\in\Theta, w\in\mathcal{W}}|K(\theta,w)|\le \bar C$ a.s., where $\bar C=1/G(t_0)$.  For all $\theta\in\Theta, w\in\mathcal{W}$, this yields,
\begin{align*}
&\left|E[K(\theta,w)\{K(\theta,w)-K(\theta_*,w)\}]\right|\\
&\le\bar CE[\left|K(\theta,w)-K(\theta_*,w)\right|]\\
&\le \bar C^2E\left[\left| I(T\le \varphi_{\theta1}(Z,u),D=1 )-I(Y\le \varphi_{\theta_*1}(Z,u), D=1)\right|\right]\\
&=\bar  C^2E\left[I(T\in [\min(\varphi_{\theta1}(Z,u),\varphi_{\theta_*1}(Z,u)), \varphi_{\theta1}(Z,u)\vee  \varphi_{\theta_*1}(Z,u)], D=1)\right] \\
&\le \bar C^2E\left[I(T\in [  \varphi_{\theta_*1}(Z,u)-C_u||\theta-\theta_*||,   \varphi_{\theta_*1}(Z,u)+C_u||\theta-\theta_*||], D=1)\right]\\
&\le \bar  C^2E\left[\P(T\in [  \varphi_{\theta_*1}(Z,u)-C_u||\theta-\theta_*||,   \varphi_{\theta_*1}(Z,u)+C_u||\theta-\theta_*||]|Z,D=1)\right],
\end{align*}
where $C_u$ is defined in Assumption \ref{paramregul} (vi).
A similar bound holds for $E[K(\theta,w)\{K(\theta,w)-K(\theta_*,w)\}] $.
Hence, \eqref{suffse} follows from the fact that the density of $T$ given $Z,D=1$ is uniformly bounded by Assumption \ref{paramregul} (vii). This shows \eqref{se}.

Next, by the inequality of Cauchy-Schwarz, we have
\begin{align*}
&\lim\limits_{\epsilon\to 0} \sup\limits_{w \in \mathcal{W}, \theta\in\Theta:\ ||\theta-\theta_*||\le \epsilon}\left|\widehat{I}_2(\theta,u,w)-\widehat{I}_2(\theta_*,u,w)\right|^2\\
&\le \lim\limits_{\epsilon\to 0} \sup\limits_{w \in \mathcal{W}, \theta\in\Theta:\ ||\theta-\theta_*||\le \epsilon}\frac1n\sum_{i=1}^n \left|K_i (\theta,w)- K_i (\theta_*,w)\right|\\
&\quad \times \lim\limits_{\epsilon\to 0} \sup\limits_{w \in \mathcal{W}, \theta\in\Theta:\ ||\theta-\theta_*||\le \epsilon}\frac1n\sum_{i=1}^n  \left\{\frac{I(T_i\le t_0)}{\widehat{G}(T_i)}\left[\widehat{G}(T_i)-G(T_i)\right]\right\}^2 \left|K_i (\theta,w)- K_i (\theta_*,w)\right|\\
&\le \lim\limits_{\epsilon\to 0} \sup\limits_{w \in \mathcal{W}, \theta\in\Theta:\ ||\theta-\theta_*||\le \epsilon}\frac1n\sum_{i=1}^n \left|K_i (\theta,w)- K_i (\theta_*,w)\right|\times O_P\left(n^{-1}\right),
\end{align*}
where the last equality is due to \eqref{KMU} and $\sup_{\theta\in\Theta, w\in\mathcal{W}}|K(\theta,w)|\le \bar C$ a.s. Next, it holds that
$$|K(\theta,w)-K(\theta_*,w)|=K(\theta,w)[1-K(\theta_*,w)]+ K(\theta_*,w)[1-K(\theta,w)].
$$
Hence, the class 
\begin{align*}
\begin{array}{lll}\R_+\times\R_+\times\mathcal{Z}\times \{0,1\}\times \mathcal{W} & \mapsto &\R\\
(t,c,z,d,w)  & \mapsto &\frac{I(t\le c)}{G(\min(t,t_0)}I( d=1,w\le \omega)|I(t\le \varphi_{\theta 1}(z,u))-I(t\le \varphi_{\theta_* 1}(z,u)|,
\end{array}
\end{align*}
$ \theta\in\Theta, \omega \in\mathcal{W}$,
is Glivenko-Cantelli because it is the sum and the product of uniformly bounded Glivenko-Cantelii classes (Corollary 9.27 in \cite{kosorok2008introduction}). Moreover, 
\begin{align*}
&\lim\limits_{\epsilon\to 0} \sup\limits_{w \in \mathcal{W}, \theta\in\Theta:\ ||\theta-\theta_*||\le \epsilon} E\left[|K(\theta,w)-K(\theta_*,w)|\right] =0\end{align*}
by the arguments used to prove \eqref{suffse}. Hence, by the Glivenko-Cantelli theorem, we have 
\begin{align*}
\lim\limits_{\epsilon\to 0} \sup\limits_{w \in \mathcal{W}, \theta\in\Theta:\ ||\theta-\theta_*||\le \epsilon}\frac1n\sum_{i=1}^n |K_i(\theta,w)-K_i(\theta_*,w)|
=o_P(1).
\end{align*}
It implies that $\lim\limits_{\epsilon\to 0} \sup\limits_{w \in \mathcal{W}, \theta\in\Theta:\ ||\theta-\theta_*||\le \epsilon}|\widehat{I}_2(\theta,u,w)-\widehat{I}_2(\theta_*,u,w)|=o_P(n^{-1/2})$.
This and \eqref{se} result in (iv). 
\end{Proof}

\begin{Lemma}\label{convM}
Let $u\in [0,u_m)$. Under the assumptions of Theorem \ref{AN}, the following holds:
\begin{itemize} 
\item[(i)]$\sup\limits_{w\in\mathcal{W},\theta \in \Theta}\left|(\widehat{M}_\theta-M_\theta)(u,w)\right|=O_P(n^{-1/2}).$
\item[(ii)] $\sup\limits_{w\in\mathcal{W}}\lim\limits_{\theta \to\theta_*}n^{1/2}\left|(\widehat{M}_\theta-M_\theta)(u,w)-\widehat{M}_{\theta_*}(u,w)\right|=o_P(1).$
\item[(iii)] $\sup\limits_{\theta \in \Theta}\left|\widehat{L}(\theta)-L(\theta)\right|=o_P(1).$
\end{itemize}
\end{Lemma}
\begin{Proof} It holds that 
\begin{align*}
(\widehat{M}_\theta-M_\theta)(u,w)&=(\widehat{F}_0-F_0)(\varphi_{\theta0}(u),w) + (\widehat{F}_1-F_1)(\varphi_{\theta1}(\cdot,u),w)\\
&\quad -(1-e^{-u})(\widehat{F}_W-F_W)(w).
\end{align*}
Moreover, since the class $\{\omega\mapsto I(\omega\le w),\ w\in\mathcal{W}\}$ is Donsker, we have 
 $$\sup\limits_{w\in\mathcal{W},\theta \in \Theta}\left|(\widehat{F}_W-F_W)(w)\right|=O_P(n^{-1/2}).$$
 This and Lemma \ref{convF} (i) and (iii) yield (i). Statement (ii) follows from $M_{\theta_*}\equiv 0$, 
 \begin{align*}
(\widehat{M}_\theta-M_\theta)(u,w)-\widehat{M}_{\theta_*}(u,w)&=
(\widehat{F}_0-F_0)(\varphi_{\theta0}(u),w)-(\widehat{F}_0-F_0)(\varphi_{\theta_*0}(u),w)\\
&\quad + (\widehat{F}_1-F_1)(\varphi_{\theta1}(u),w)-(\widehat{F}_1-F_1)(\varphi_{\theta_*1}(\cdot,u),w)
\end{align*}
and Lemma \ref{convF} (ii) and (iv). 

Finally, by (i), we have $\sup\limits_{\theta \in \Theta}|\widehat{L}(\theta)-\bar L(\theta)|=o_P(1)$. Therefore, it suffices to show that $\sup\limits_{\theta \in \Theta}|\bar L(\theta)-L(\theta)|=o_P(1)$ to obtain (iii). This holds because the class $\{w\mapsto M_{\theta}(u,w),\theta\in\Theta\}$ is Glivenko-Cantelli by Theorem 2.7.11 in \cite{VdVW} since $\theta\mapsto M_\theta(u,w)$ is Lipschitz by Assumption \ref{paramregul} (iii). 
\end{Proof}

\subsection{A weak convergence result} \label{subsec.WC}

Let $X_i=(Y_i, \delta_i,\tilde Z_i, \tilde D_i, W_i)$ and $\mathcal{X}$ be the support of $X_1$. 

\begin{Lemma}\label{weakF} Let $u \in [0,u_m)$.   Under the assumptions of Theorem \Ref{AN}, there exist mappings $m_{0uw}, m_{1uw}$ from $\mathcal{X}^2$ to $\R$ such that $m_{0uw}(x_1,x_2)$ and $m_{1uw}(x_1,x_2)$ are bounded uniformly in $w\in\mathcal{W}$, $x_1,x_2\in\mathcal{X}$ and 
\begin{itemize}
\item[(i)] $\sup\limits_{w\in\mathcal{W}}\left|(\widehat{F}_0-F_0)(\varphi_{\theta_*0}(u),w)-n^{-2}\sum_{i=1}^n\sum_{j=1}^n  m_{0uw}(X_i,X_j)\right|=o_P(n^{-1/2})$;
\item[(ii)] $\sup\limits_{w\in\mathcal{W}}\left|(\widehat{F}_1-F_1)(\varphi_{\theta_*1}(\cdot,u),w)-n^{-2}\sum_{i=1}^n\sum_{j=1}^n m_{1uw}(X_i,X_j)\right|=o_P(n^{-1/2})$.
\end{itemize}
\end{Lemma}
\begin{Proof} The proofs of (i) and (ii) are similar. Hence, we only show (ii). We have 
\begin{align*}
(\widehat{F}_1-F_1)(\varphi_{\theta_*1}(\cdot,u),w)&=\widehat{I}_1(u,w) -\widehat{I}_2(u,w)+\widehat{I}_3(u,w), 
\end{align*}
where 
\begin{align*}
\widehat{I}_1(u,w)&=\frac{1}{n}\sum_{i=1}^n  \frac{\delta_i}{G(Y_i)} I(Y_i\le \varphi_{\theta_*1}(\tilde Z_i,u),\tilde D_i=1, W_i\le w)\\
&\quad- E\left[\frac{\delta}{G(Y)}I(T\le \varphi_{\theta_*1}(\tilde Z,u),\tilde D=1,W\le w)\right];\\
\widehat{I}_2(u,w)&=\frac{1}{n}\sum_{i=1}^n  \frac{\delta_i}{G^2(Y_i)} I(Y_i\le\varphi_{\theta_*1}(\tilde Z_i,u),\tilde D_i=1, W_i\le w)\left[\widehat{G}(Y_i)-G(Y_i)\right];\\
\widehat{I}_3(u,w)&=\frac{1}{n}\sum_{i=1}^n  \frac{\delta_i}{\widehat{G}(Y_i)G^2(Y_i)} I(Y_i\le \varphi_{\theta_*1}(\tilde Z_i,u),\tilde D_i=1, W_i\le w)\left[\widehat{G}(Y_i)-G(Y_i)\right]^2.
\end{align*}
Let $t_0=\sup_{ z\in \mathcal{Z}}\varphi_{\theta_*1} (z,u_m)$. By Theorem 1.1 in \cite{gill1983large} and the fact that $t_0<c_0$, 
$$ \sup_{s\in[0,t_0]}\left|\widehat{G}(s) -G(s)\right| =O_P\left(n^{-1/2}\right).$$
As a result, since $G(t_0)>0$, we get 
$$\sup_{w\in\mathcal{W}} |I_3(u,w)|=O_P\left(n^{-1}\right).$$
Next, by Theorem 1 in \cite{lo1986product}, there exist mean zero bounded stochastic processes $\{\xi_j(s)\}_{s\in\R_+}$ depending on $(Y_j,\delta_j)$ such that $$\sup_{s\in[0,t_0]}\left|\widehat{G}(s)-G(s)-\frac1n\sum_{j=1}^n \xi_j(s)\right|=o_P\left(n^{-1/2}\right).$$
Hence, because $G(t_0)>0$, we have
\begin{align*}
 I_2(u,w)
 &=\frac{1}{n^2}\sum_{i=1}^n \sum_{j=1}^n \frac{\delta_i}{G^2(Y_i)} I(Y_i\le \varphi_{\theta_*1}(\tilde Z_i,u),\tilde D_i=1, W_i\le w)\xi_j(Y_i) + o_P\left(\frac1{\sqrt{n}}\right).
\end{align*}
Letting  \begin{align*}m_{1uw}(X_i,X_j) &=\frac{\delta_i}{G(Y_i)} I(Y_i\le \varphi_{\theta_*1}(\tilde Z_i,u),\tilde D_i=1, W_i\le w)\\
&\quad  - E\left[\frac{\delta}{G(Y)}I(T\le \varphi_{\theta_*1}(\tilde Z,u),\tilde D=1,W\le w)\right]\\
&\quad - \frac{\delta_i}{G^2(Y_i)} I(Y_i\le \varphi_{\theta_*1}(\tilde Z_i,u),\tilde D_i=1, W_i\le w)\xi_j(Y_i),\end{align*} we obtain the result.
\end{Proof}

\begin{Lemma}\label{weakM} Under thz assumptions of theorem \ref{AN}, there exists a zero-mean Gaussian random variable $V$ such that $$\widehat{V}=\frac{2}{n^{1/2}m}\sum_{k=1}^m \sum_{i=1}^n p(u_k)\frac{\partial M_{\theta_*}}{\partial \theta}(u_k,W_i)\widehat{M}_{\theta_*}(u_k,W_i)\xrightarrow{d}V.$$
\end{Lemma}

\begin{Proof} By Lemma \ref{weakF}, we have
$$\widehat{V}= \sqrt{n}\frac{1}{n^3}\sum_{i,j,k=1,\dots, n}m(X_i,X_j,X_k) +o_P(n^{-1/2}),$$
where 
\begin{align*} 
m(X_i,X_j,X_k)  &=\frac2m \sum_{\ell=1}^m \sum_{i,j,k=1,\dots, n}\Big\{p(u_\ell)\frac{\partial M_{\theta_*}}{\partial \theta}(u_\ell,W_k)\\
&\quad \times \left[m_{0u_\ell}(X_i,X_j) +m_{1u_\ell}(X_i,X_j) - (1-e^{-u_\ell})(I(W_i\le W_k) -F_W(W_k))\right]\Big\}.
\end{align*}
Let $\widetilde m$ be the symmetrization of $m$ defined by  \begin{align*}6\widetilde m(X_i,X_j,X_k)&= m(X_i,X_j,X_k)+ m(X_i,X_k,X_j)+ m(X_j,X_i,X_k)\\
&\quad +m(X_j,X_k,X_i)+m(X_k,X_i,X_j)+m(X_k,X_j,X_i).\end{align*}
It holds that 
$$\widehat{V}= \sqrt{n}\frac{1}{n^3}\sum_{i,j,k=1,\dots, n}\widetilde m(X_i,X_j,X_k) +o_P(n^{-1/2}).$$
Moreover, since $\widetilde m$ is bounded, we have 
$$\widehat{V}= \sqrt{n}\frac{1}{n^3}\sum_{i<j<k\in\{1,\dots n\}}\frac16 \widetilde m(X_i,X_j,X_k) +o_P(n^{-1/2}).$$
The term $n^{-3}\sum_{i<j<k\in\{1,\dots n\}}\frac16 \widetilde m(X_i,X_j,X_k) $ is a U-statistic with bounded kernel. Hence, the result follows from the central limit theorem for U-statistics (Theorem 1.1 in \cite{bose2018u}).
\end{Proof}

\section{Proof of Theorem \ref{ANB}} \label{S.2}

The body of the proof of Theorem \ref{ANB} is given in Section \ref{subsec.pANB}. It relies on technical lemmas proved in Sections \ref{subsec.UCB} and \ref{subsec.WCB}. Let $\widehat{L}_b, \widehat{M}_{b}$ be defined similarly as $\widehat{L},\widehat{M}$, except that they are based on the bootstrap data rather than on the original data. 


\subsection{Main proof}\label{subsec.pANB}

\noindent\textbf{Step 1: Consistency.}  Using Lemma \ref{convMB2} (iii) and following arguments similar to step 1 of the proof of Theorem 4.1, we get that $\widehat\theta_b-\theta_* \xrightarrow{\P} 0$ \\

\noindent \textbf{Step 2: Asymptotic expansion.} 
For $\theta \in \Theta$, let $\bar{L}_b(\theta)= \frac{1}{nm}\sum_{i=1}^n \sum_{j=1}^mp(u_j) M_\theta(u_j,W_{bi})^2$. Consider a sequence of $\theta\in\Theta$ such that $\theta\to\theta_*$. All the following statements are in $P$-probability with respect to the original sample. Remark that
 \begin{align}
\notag \widehat{L}_b( \theta)-\bar{L}_b(\theta)&= \frac{1}{nm}\sum_{i=1}^n \sum_{j=1}^m p(u_j)\left[\left(\widehat{M}_{ b \theta}- M_{ \theta}\right)(u_j,W_{bi})\right]^2\\
 \label{TMDB1}&\quad + \frac{2}{nm}\sum_{i=1}^n \sum_{j=1}^m   M_{\theta}(u_j,W_{bi})p(u_j)\left[\left(\widehat{M}_{b\theta}- M_{ \theta}\right)(u_j,W_{bi})\right].
 \end{align}
Moreover,  it holds that \begin{align*}&\frac{1}{nm}\sum_{i=1}^n \sum_{j=1}^mp(u_j) \left[\left(\widehat{M}_{ b\theta}- M_{  \theta}\right)(u_j,W_{bi})\right]^2\\
&=\frac{1}{nm}\sum_{i=1}^n \sum_{j=1}^mp(u_j) \widehat{M}_{ b\theta_*}(u_j,W_{bi})^2+I_{b1}(\theta)+I_{b2}(\theta),\end{align*} where 
  \begin{align*}
 I_{b1}(\theta)&= \frac{1}{nm}\sum_{i=1}^n \sum_{j=1}^mp(u_j) \left[\left(\widehat{M}_{ b\theta}-M_{ \theta}-\widehat{M}_{b\theta_*}\right)(u_j,W_{bi}) \right]^2\\
I_{b2} (\theta)&= \frac{2}{nm}\sum_{i=1}^n \sum_{j=1}^mp(u_j) \widehat{M}_{b\theta_*}(u_j,W_{bi}) \left[\left(\widehat{M}_{  b\theta}-M_{ \theta}-\widehat{M}_{b\theta_*}\right)(u_j,W_{bi}) \right].
 \end{align*}
By Lemma \ref{convMB2} (ii), we have $|I_{b1}(\theta)|=o_{P^*}(n^{-1})$. Then, by the inequality of Cauchy-Schwarz, it holds that 
\begin{align*}
I_{b2}(\theta)&\le 2 \sqrt{\frac{1}{nm}\sum_{i=1}^n \sum_{j=1}^mp(u_j)^2\widehat{M}_{b\theta_*}(u_j,W_{bi})^2} \\
& \quad \times \sqrt{\frac{1}{nm}\sum_{i=1}^n \sum_{j=1}^m\left[\left(\widehat{M}_{  b\theta}-M_{ \theta}-\widehat{M}_{b\theta_*}\right)(u_j,W_{bi}) \right]^2},\end{align*}
which is $o_{P^*}(n^{-1})$\ by Lemma \ref{convMB2}.
 As a result, we obtain
 \begin{align}\notag &\frac{1}{nm}\sum_{i=1}^n \sum_{j=1}^mp(u_j) \left[\left(\widehat{M}_{ b\theta}-M_{  \theta}\right)(u_j,W_{bi})\right]^2\\&=\frac{1}{nm}\sum_{i=1}^n \sum_{j=1}^mp(u_j) \widehat{M}_{ b\theta_*}(u_j,W_{bi})^2+o_{P^*}(n^{-1}). \label{TMDB2}\end{align}
Then, using a Taylor expansion of $\theta\mapsto M_{\theta}(u_j,W_{bi})$ around $\theta_*$ and Lemma \ref{convMB2}, it holds that
 \begin{align}
 \notag &\frac{1}{nm}\sum_{i=1}^n \sum_{j=1}^m p(u_j)M_{{\theta}}(u_j,W_{bi})\left[\left(\widehat{M}_{b{\theta}}-M_{{\theta}}\right)(u_j,W_{bi})\right]\\
\notag &= (\theta -\theta_*)^\top \frac{1}{nm}\sum_{i=1}^n \sum_{j=1}^m p(u_j)\frac{\partial M_{{\theta_*}}}{\partial \theta}(u_j,W_{bi})\left[\left(\widehat{M}_{b{\theta_*}}-M_{{\theta_*}}\right)(u_j,W_{bi})\right]+o_{P^*}(||{\theta}-\theta_*||^2)\\
&= (\theta -\theta_*)^\top  n^{-1/2}\widehat{V}_b +o_{P^*}(n^{-1/2}||\theta-\theta_*||+||{\theta}-\theta_*||^2),
\label{TMDB3} 
 \end{align}
 where $\widehat{V}_b$ is defined in Lemma \ref{weakMB}. 
Next, \eqref{TMDB1}, \eqref{TMDB2}, \eqref{TMDB3} yield
 \begin{align*}
 \notag \widehat{L}_b( \theta)-\bar{L}_b(\theta)&=
 \frac{1}{nm}\sum_{i=1}^n \sum_{j=1}^mp(u_j) \widehat{M}_{ b\theta_*}(u_j,W_{bi})^2\\
 &\quad +(\theta -\theta_*)^\top n^{-1/2}\widehat{V}_b+o_{P^*}(n^{-1/2}||{\theta}-\theta_*||+||{\theta}-\theta_*||^2+n^{-1}).
 \end{align*}
Therefore, we get, by second-order Taylor expansions,
\begin{align*}
&\widehat{L}_b( \theta) -\widehat{L}_b( \widehat\theta)+o_{P}(n^{-1/2}||{\theta}-\theta_*||+||{\theta}-\theta_*||^2+n^{-1})\\
&= (\theta -\widehat \theta)^\top n^{-1/2}\widehat{V}_b+\frac12 (\theta-\theta_*)^\top\nabla^2\bar{L}_b(\theta_*)(\theta- \theta_*)- \frac12 (\widehat\theta-\theta_*)^\top\nabla^2\bar{L}_b(\theta_*)(\widehat\theta-\theta_*)\\
&= (\theta -\widehat \theta)^\top n^{-1/2}\widehat{V}_b+\frac12 (\theta-\widehat \theta)^\top\nabla^2\bar{L}_b(\theta_*)(\theta- \widehat\theta)+(\theta-\widehat\theta)^\top\nabla^2\bar{L}_b(\theta_*)(\widehat\theta-\theta_*)\\
&=(\theta -\widehat \theta)^\top n^{-1/2}(\widehat{V}_b-\widehat V)+\frac12 (\theta-\widehat \theta)^\top\nabla^2\bar{L}_b(\theta_*)(\theta- \widehat\theta)+o_{P^*}(n^{-1/2}||\theta-\widehat\theta||),
\end{align*}
where in the last equality, we used $$\widehat\theta-\theta_*= -n^{-1/2} \nabla^2 \bar{L}(\theta_*)^{-1}\widehat V+o_P(n^{-1/2}) = {-n^{-1/2}} \nabla^2 \bar{L}_b(\theta_*)^{-1}\widehat V+o_{P^*}(n^{-1/2}).$$ 
Since $\widehat \theta -\theta_*=O_{P^*}(n^{-1/2})$, we obtain
\begin{align}\notag
\widehat{L}_b( \theta) -\widehat{L}_b( \widehat\theta)&= (\theta -\widehat \theta)^\top n^{-1/2}(\widehat{V}_b-\widehat V)+\frac12 (\theta-\widehat \theta)^\top\nabla^2\bar{L}_b(\theta_*)(\theta- \widehat\theta)\\
\label{AEB} &\quad+o_{P^*}(n^{-1/2}||{\theta}-\widehat\theta||+||{\theta}-\widehat\theta||^2+n^{-1}).
\end{align}
Since \eqref{AE} and \eqref{AEB} are similar, the rest of the proof follows the same arguments as steps 3 and 4 of Theorem \ref{AN} and is therefore omitted. 

\subsection{Uniform convergence results for the bootstrap}\label{subsec.UCB}

\begin{Lemma}\label{convFB}Let $u\in [0,u_m)$. Under the assumptions of Theorem \ref{ANB}, in $P$-probability with respect to the original sample, the following holds:
\begin{itemize} 
\item[(i)] $\sup\limits_{w\in\mathcal{W},\theta \in \Theta}\left|(\widehat{F}_{b0}-\widehat F_{0})(\varphi_{\theta0}(u),w)\right|=O_{P^*}(n^{-1/2})$;
\item[(ii)] $\lim\limits_{\epsilon \to 0 }\sup\limits_{w\in\mathcal{W},\ \theta \in\Theta:\ ||\theta-\widehat\theta||\le \epsilon}n^{1/2}\left|(\widehat{F}_ {b0}-\widehat F_0)(\varphi_{\theta0}(u),w)-(\widehat{F}_{b0}-\widehat F_0)(\varphi_{\widehat\theta0}(u),w)\right|=o_{P^*}(1)$;
\item[(iii)] $\sup\limits_{w\in\mathcal{W},\theta \in \Theta}\left|(\widehat{F}_{b1}- \widehat F_1)(\varphi_{\theta1}(\cdot,u),w)\right|=O_{P^*}(n^{-1/2})$;
\item[(iv)] $\lim\limits_{\epsilon \to 0 }\sup\limits_{w\in\mathcal{W},\ \theta \in\Theta:\ ||\theta-\widehat\theta||\le \epsilon }n^{1/2}\left|(\widehat{F}_{b1}-\widehat F_1)(\varphi_{\theta1}(\cdot,u),w)-(\widehat{F}_{b1}-\widehat F_1)(\varphi_{\widehat\theta1}(\cdot,u),w)\right|=o_{P^*}(1).$
\end{itemize}
\end{Lemma}

\begin{Proof} The proof is similar to that of Lemma \ref{convF}. Hence, we just mention the key differences here. The usual Donsker theorem is replaced by the Donsker theorem for the bootstrap (Theorem 3.6.3 in \cite{VdVW}). The uniform convergence for the bootstrap Kaplan-Meier estimator is obtained through Theorem 1 in \cite{lo1986product}, the Donsker theorem for the bootstrap and the continuous mapping theorem.
\end{Proof}

\begin{Lemma}\label{convMB}Let $u\in [0,u_m)$. Under the assumptions of Theorem \Ref{ANB}, in $P$-probability with respect to the original sample, the following holds:
\begin{itemize} 
\item[(i)]$\sup\limits_{w\in\mathcal{W},\theta \in \Theta}\left|(\widehat{M}_{b\theta}-\widehat{M}_\theta)(u,w)\right|=O_{P^*}(n^{-1/2})$;
\item[(ii)] $\lim\limits_{\epsilon \to 0 }\sup\limits_{w\in\mathcal{W},\ \theta \in\Theta:\ ||\theta-\widehat\theta|| \le \epsilon}n^{1/2}\left|(\widehat{M}_{b\theta}-\widehat{M}_\theta)(u,w)-(\widehat{M}_{b\widehat\theta}-\widehat{M}_{\widehat\theta})(u,w)\right|=o_{P^*}(1)$;
\item[(iii)] $\sup_{\theta\in\Theta} |\widehat{L}_b(\theta)-\widehat{L}(\theta)|=o_{P^*}(1).$
\end{itemize}
\end{Lemma}
\begin{Proof} The proof is similar to that of Lemma \ref{convM}. The only difference is that the usual Donsker theorem is replaced by the Donsker theorem for the bootstrap (Theorem 3.6.3 in \cite{VdVW}) and the usual Glivenko-Cantelli theorem is replaced by the Glivenko-Cantelli theorem for the bootstrap.
\end{Proof}

\begin{Lemma}\label{convMB2}Let $u\in [0,u_m)$. in $P$-probability with respect to the original sample, the following holds:
\begin{itemize} 
\item[(i)]$\sup\limits_{w\in\mathcal{W},\theta \in \Theta}\left|(\widehat{M}_{b\theta}-M_\theta)(u,w)\right|=O_{P^*}(n^{-1/2})$;
\item[(ii)] $\lim\limits_{\epsilon \to 0 }\sup\limits_{w\in\mathcal{W},\ \theta \in\Theta:\ ||\theta-\theta_*|| \le \epsilon} n^{1/2}\left|(\widehat{M}_{b\theta}-M_\theta)(u,w)-\widehat{M}_{b\theta_*}(u,w)\right|=o_{P^*}(1)$;
\item[(iii)] $\sup_{\theta\in\Theta} |\widehat{L}_b(\theta)-L(\theta)|=o_{P^*}(1)$
\end{itemize}
\end{Lemma}
\begin{Proof} The results are a direct consequence of the triangle inequality, Lemmas \ref{convM} and \ref{convMB} and the fact that $\widehat{\theta}-\theta_*=o_P(1)$.
\end{Proof}

\subsection{A weak convergence result for the bootstrap}\label{subsec.WCB}

Let $X_{bi} = (Y_{bi}, \delta_{bi},\tilde Z_{bi},\tilde D_{bi},W_{bi})$.

\begin{Lemma}\label{weakFB} Let $u \in [0,u_0)$.  Under the assumptions of Theorem \ref{ANB}, there exist mappings $m_{b0uw}, m_{b1uw}$ from $\mathcal{X}^2$ to $\R$ such that $m_{b0uw}(x_1,x_2)$ and $m_{b1uw}(x_1,x_2)$ are bounded uniformly in $w\in\mathcal{W}$, $x_1,x_2\in\mathcal{X}$ and, in $P$-probability with respect to the original sample,
\begin{align*}(i) \quad &\sup_{w\in\mathcal{W}}\left|(\widehat{F}_{b0}-\widehat{F}_0)(\varphi_{\widehat\theta0}(u),w)-n^{-2}\sum_{i=1}^n\sum_{j=1}^n  \Big(m_{b0uw}(X_{bi},X_{bj})-m_{b0uw}(X_{i},X_{j}) \Big)\right|\\
&=o_{P^*}(n^{-1/2}).;\end{align*} 
\begin{align*}(ii) \quad &\sup_{w\in\mathcal{W}}\left|(\widehat{F}_{b1}-\widehat F_1)(\varphi_{\widehat\theta1}(\cdot,u),w)-n^{-2}\sum_{i=1}^n\sum_{j=1}^n \Big(m_{b1uw}(X_{bi},X_{bj})-m_{b1uw}(X_{i},X_{j})\Big) \right|\\
&=o_{P^*}(n^{-1/2}).\end{align*}

\end{Lemma}

\begin{Proof} The proof is similar to that of Lemma \ref{weakF} except that we use the second statement in  Theorem 1 in \cite{lo1986product} rather than the first statement in the same theorem.
\end{Proof}

\begin{Lemma}\label{weakMB} Let \begin{align*}\widehat{V}_b&=\frac{2}{n^{1/2}m}\sum_{k=1}^m \sum_{i=1}^n p(u_k)\frac{\partial M_{\theta_*}}{\partial \theta}(u_k,W_{bi})\widehat{M}_{b\theta_*}(u_k,W_{bi}).\end{align*} Under Assumptions 3.5, 4.1, and 4.2, $$\widehat{V}_b-\widehat V\xrightarrow{d}V,$$
conditionally on the original sample in $P$-probability, where $\widehat V$ and $V$ are defined in Lemma \ref{weakM}.
\end{Lemma}
\begin{Proof} The proof is similar to that of Lemma \ref{weakM}, except that we replace the central limit theorem for U-statistics by the bootstrap central limit theorem for U-Statistics (\cite{bickel1981asymptotic}).
\end{Proof}

\section{Additional simulation results}\label{S.3}

\subsection{Simulations for the Weibull design without censoring.}
\tr{Table \ref{tab:simW1} reports simulation results for the Weibull design without censoring. The results are qualitatively similar to that of the Weibull design with censoring as in Table 1 of the main text. There is however some improvement when there is no censoring, which is to be expected.}
\afterpage{
\clearpage
  \thispagestyle{empty}
  \begin{landscape}
\begin{table}[!p]
\centering
\adjustbox{max width=1.1\textwidth}{
\begin{tabular}{l | c c c c | c c c c | c c c c | c c c c}
\hline \hline
~ & $\hat{\theta}_{00}$ & $\hat{\theta}_{10}$ & $\hat{\theta}_{01}$ & $\hat{\theta}_{11}$ & $\hat{\theta}_{00}$ & $\hat{\theta}_{10}$ & $\hat{\theta}_{01}$ & $\hat{\theta}_{11}$ & $\hat{\theta}_{00}$ & $\hat{\theta}_{10}$ & $\hat{\theta}_{01}$ & $\hat{\theta}_{11}$ & $\hat{\theta}_{00}$ & $\hat{\theta}_{10}$ & $\hat{\theta}_{01}$ & $\hat{\theta}_{11}$ \\ \hline \hline
\multicolumn{17}{c}{\textsc{No Censoring}}\\ \hline
\multicolumn{17}{l}{\textbf{$n = 500$}}\\ \hline
~ & \multicolumn{4}{c|}{$\alpha = 0.25, \beta = 1, \bar{D} = 0.46$} & \multicolumn{4}{c|}{$\alpha = 0.75, \beta = 1, \bar{D} = 0.49$}& \multicolumn{4}{c|}{$\alpha = 0.25, \beta = 0.5, \bar{D} = 0.40$} & \multicolumn{4}{c}{$\alpha = 0.75, \beta = 0.5, \bar{D} = 0.43$}\\
 \hline
Bias & 0.013 & 0.033 & 0.015 & -0.004 & 0.010 & 0.069 & 0.017 & -0.014 & 0.036 & 0.353 & 0.019 & -0.082 & 0.020 & 0.462 & 0.020 & -0.141 \\ 
 SE & 0.105 & 0.424 & 0.115 & 0.432 & 0.104 & 0.864 & 0.123 & 0.462 & 0.165 & 3.308 & 0.131 & 0.769 & 0.151 & 3.352 & 0.135 & 0.747 \\ 
 $90\%$ & 0.899 & 0.857 & 0.907 & 0.823 & 0.867 & 0.825 & 0.871 & 0.783 & 0.730 & 0.657 & 0.872 & 0.674 & 0.705 & 0.619 & 0.829 & 0.623 \\ 
 $95\%$ & 0.940 & 0.957 & 0.955 & 0.885 & 0.930 & 0.914 & 0.941 & 0.881 & 0.845 & 0.822 & 0.929 & 0.782 & 0.803 & 0.745 & 0.911 & 0.763 \\ 
 $99\%$ & 0.985 & 0.996 & 0.991 & 0.980 & 0.991 & 0.992 & 0.987 & 0.966 & 0.930 & 0.971 & 0.992 & 0.920 & 0.943 & 0.925 & 0.992 & 0.940 \\ 
  \hline

\multicolumn{17}{l}{\textbf{$n = 1000$}}\\ \hline
~ & \multicolumn{4}{c|}{$\alpha = 0.25, \beta = 1, \bar{D} = 0.46$} & \multicolumn{4}{c|}{$\alpha = 0.75, \beta = 1, \bar{D} = 0.49$}& \multicolumn{4}{c|}{$\alpha = 0.25, \beta = 0.5, \bar{D} = 0.40$} & \multicolumn{4}{c}{$\alpha = 0.75, \beta = 0.5, \bar{D} = 0.42$}\\
 \hline
Bias & 0.011 & -0.003 & 0.004 & 0.010 & 0.010 & 0.007 & 0.004 & 0.000 & 0.023 & 0.149 & 0.004 & 0.010 & 0.016 & 0.065 & 0.007 & -0.054 \\ 
 SE & 0.076 & 0.309 & 0.077 & 0.284 & 0.076 & 0.317 & 0.083 & 0.297 & 0.115 & 2.822 & 0.088 & 0.509 & 0.113 & 0.748 & 0.095 & 0.512 \\ 
 $90\%$ & 0.884 & 0.918 & 0.911 & 0.925 & 0.896 & 0.906 & 0.908 & 0.910 & 0.773 & 0.672 & 0.890 & 0.750 & 0.740 & 0.654 & 0.882 & 0.741 \\ 
 $95\%$ & 0.954 & 0.954 & 0.954 & 0.965 & 0.950 & 0.959 & 0.949 & 0.963 & 0.863 & 0.841 & 0.935 & 0.866 & 0.839 & 0.809 & 0.924 & 0.843 \\ 
 $99\%$ & 0.999 & 0.998 & 0.994 & 0.996 & 0.996 & 0.997 & 0.985 & 0.996 & 0.976 & 0.965 & 0.985 & 0.948 & 0.946 & 0.964 & 0.975 & 0.961 \\ 
  \hline

\multicolumn{17}{l}{\textbf{$n = 3000$}}\\ \hline
~ & \multicolumn{4}{c|}{$\alpha = 0.25, \beta = 1, \bar{D} = 0.46$} & \multicolumn{4}{c|}{$\alpha = 0.75, \beta = 1, \bar{D} = 0.49$}& \multicolumn{4}{c|}{$\alpha = 0.25, \beta = 0.5, \bar{D} = 0.40$} & \multicolumn{4}{c}{$\alpha = 0.75, \beta = 0.5, \bar{D} = 0.42$}\\
 \hline
Bias & 0.003 & 0.002 & 0.002 & 0.005 & 0.002 & 0.005 & 0.002 & 0.005 & 0.007 & -0.002 & 0.003 & 0.014 & 0.005 & 0.018 & 0.003 & -0.007 \\ 
 SE & 0.042 & 0.171 & 0.044 & 0.163 & 0.042 & 0.174 & 0.049 & 0.172 & 0.061 & 0.390 & 0.050 & 0.265 & 0.063 & 0.405 & 0.055 & 0.271 \\ 
 $90\%$ & 0.893 & 0.913 & 0.923 & 0.906 & 0.889 & 0.915 & 0.904 & 0.919 & 0.832 & 0.827 & 0.876 & 0.869 & 0.813 & 0.811 & 0.876 & 0.854 \\ 
 $95\%$ & 0.961 & 0.969 & 0.962 & 0.951 & 0.956 & 0.960 & 0.973 & 0.962 & 0.912 & 0.912 & 0.951 & 0.924 & 0.905 & 0.902 & 0.947 & 0.936 \\ 
 $99\%$ & 0.990 & 0.996 & 0.996 & 0.998 & 0.992 & 0.998 & 0.992 & 0.996 & 0.990 & 0.984 & 0.996 & 0.986 & 0.988 & 0.992 & 0.988 & 0.995 \\ 
  \hline
\end{tabular}
}
\caption{Simulation results under the Weibull design without censoring.}
\label{tab:simW1}
\end{table}
  \end{landscape}
  \clearpage
}

\subsection{Simulations results for the log-normal design} 

Results from the estimation for the log-normal model using our method are given in Tables \ref{tab:sim1} (without censoring) and \ref{tab:sim2} (with censoring) below. Let us first comment on Table \ref{tab:sim1}. Around $60\%$ of observations are treated. When the instrument is strong ($\beta=1$), we see that the results are satisfactory in terms of bias, even when $n=500$. The coverage probability improves when $n$ grows and are even almost nominal when $n=3000$. The results for a weaker instrument ($\beta=0.5$) exhibit good performance in terms of bias when $n=1000$ or $n=3000$. The coverage probabilities are further away from the nominal level in this case. The results in Table \ref{tab:sim2} show that censoring deteriorates the performance, as can be expected.

\afterpage{
\clearpage
    \thispagestyle{empty}
    \begin{landscape}
\begin{table}[!p]
\centering
\adjustbox{max width=1.1\textwidth}{
\begin{tabular}{l | c  c  c  c | c c c c | c c c c | c c c c}
\hline \hline
~ & $\hat{\theta}_{00}$ & $\hat{\theta}_{10}$  & $\hat{\theta}_{01}$ & $\hat{\theta}_{11}$ & $\hat{\theta}_{00}$ & $\hat{\theta}_{10}$  & $\hat{\theta}_{01}$ & $\hat{\theta}_{11}$ & $\hat{\theta}_{00}$ & $\hat{\theta}_{10}$  & $\hat{\theta}_{01}$ & $\hat{\theta}_{11}$ & $\hat{\theta}_{00}$ & $\hat{\theta}_{10}$  & $\hat{\theta}_{01}$ & $\hat{\theta}_{11}$ \\ \hline \hline
\multicolumn{17}{c}{\textsc{No Censoring}}\\ \hline
\multicolumn{17}{l}{\textbf{n = 500}}\\ \hline
~ & \multicolumn{4}{c}{$\alpha = 0.25, \beta = 1, \bar{D} = 0.577$} & \multicolumn{4}{c}{$\alpha = 0.75, \beta = 1, \bar{D} = 0.636$}& \multicolumn{4}{c}{$\alpha = 0.25, \beta = 0.5, \bar{D} = 0.538$} & \multicolumn{4}{c}{$\alpha = 0.75, \beta = 0.5, \bar{D} = 0.599$}\\
  \hline
Bias & -0.002 & -0.073 & -0.001 & 0.023 & -0.001 & -0.084 & -0.003 & 0.034 & -0.018 & -0.194 & -0.005 & 0.011 & -0.008 & -0.222 & 0.001 & 0.031 \\ 
  SE & 0.094 & 0.478 & 0.083 & 0.181 & 0.099 & 0.474 & 0.091 & 0.207 & 0.131 & 0.974 & 0.101 & 0.306 & 0.145 & 0.885 & 0.101 & 0.345 \\ 
  $90\%$ & 0.874 & 0.714 & 0.897 & 0.776 & 0.894 & 0.711 & 0.896 & 0.833 & 0.709 & 0.508 & 0.854 & 0.604 & 0.763 & 0.548 & 0.887 & 0.721 \\ 
  $95\%$ & 0.921 & 0.812 & 0.963 & 0.862 & 0.942 & 0.808 & 0.952 & 0.928 & 0.792 & 0.651 & 0.904 & 0.761 & 0.814 & 0.711 & 0.946 & 0.818 \\ 
  $99\%$ & 0.983 & 0.971 & 0.995 & 0.976 & 0.987 & 0.942 & 0.993 & 0.987 & 0.898 & 0.976 & 0.974 & 0.948 & 0.891 & 0.982 & 0.993 & 0.924 \\ 
   \hline

\multicolumn{17}{l}{\textbf{n = 1000}}\\ \hline
~ & \multicolumn{4}{c}{$\alpha = 0.25, \beta = 1, \bar{D} = 0.576$} & \multicolumn{4}{c}{$\alpha = 0.75, \beta = 1, \bar{D} = 0.635$}& \multicolumn{4}{c}{$\alpha = 0.25, \beta = 0.5, \bar{D} = 0.537$} & \multicolumn{4}{c}{$\alpha = 0.75, \beta = 0.5, \bar{D} = 0.598$}\\
  \hline
Bias & -0.006 & -0.024 & -0.000 & 0.009 & -0.005 & -0.018 & 0.000 & 0.011 & -0.013 & -0.081 & 0.001 & 0.015 & -0.007 & -0.078 & 0.005 & 0.015 \\ 
  SE & 0.067 & 0.311 & 0.055 & 0.118 & 0.071 & 0.265 & 0.059 & 0.111 & 0.094 & 0.720 & 0.066 & 0.209 & 0.105 & 0.644 & 0.064 & 0.189 \\ 
  $90\%$ & 0.891 & 0.804 & 0.932 & 0.859 & 0.890 & 0.790 & 0.937 & 0.878 & 0.709 & 0.537 & 0.845 & 0.643 & 0.760 & 0.548 & 0.931 & 0.792 \\ 
  $95\%$ & 0.928 & 0.876 & 0.961 & 0.929 & 0.953 & 0.860 & 0.974 & 0.943 & 0.796 & 0.679 & 0.923 & 0.785 & 0.835 & 0.690 & 0.972 & 0.867 \\ 
  $99\%$ & 0.989 & 0.994 & 0.996 & 0.991 & 0.997 & 0.968 & 0.994 & 0.996 & 0.920 & 0.971 & 0.978 & 0.945 & 0.952 & 0.923 & 0.988 & 0.955 \\ 
   \hline

\multicolumn{17}{l}{\textbf{n = 3000}}\\ \hline
~ & \multicolumn{4}{c}{$\alpha = 0.25, \beta = 1, \bar{D} = 0.576$} & \multicolumn{4}{c}{$\alpha = 0.75, \beta = 1, \bar{D} = 0.635$}& \multicolumn{4}{c}{$\alpha = 0.25, \beta = 0.5, \bar{D} = 0.538$} & \multicolumn{4}{c}{$\alpha = 0.75, \beta = 0.5, \bar{D} = 0.598$}\\
  \hline
Bias & -0.001 & -0.007 & -0.000 & 0.002 & -0.001 & -0.005 & 0.000 & 0.002 & -0.003 & -0.022 & 0.001 & 0.004 & -0.002 & -0.017 & 0.002 & 0.003 \\ 
  SE & 0.038 & 0.161 & 0.032 & 0.064 & 0.040 & 0.143 & 0.035 & 0.062 & 0.055 & 0.358 & 0.039 & 0.111 & 0.060 & 0.323 & 0.036 & 0.096 \\ 
  $90\%$ & 0.896 & 0.910 & 0.911 & 0.894 & 0.879 & 0.882 & 0.923 & 0.909 & 0.816 & 0.786 & 0.879 & 0.811 & 0.810 & 0.759 & 0.911 & 0.825 \\ 
  $95\%$ & 0.950 & 0.959 & 0.954 & 0.946 & 0.935 & 0.938 & 0.966 & 0.957 & 0.879 & 0.907 & 0.931 & 0.909 & 0.896 & 0.830 & 0.966 & 0.875 \\ 
  $99\%$ & 0.987 & 0.992 & 0.995 & 0.991 & 0.990 & 0.999 & 0.995 & 0.994 & 0.971 & 0.997 & 0.988 & 0.996 & 0.971 & 0.989 & 0.994 & 0.974 \\ 
   \hline
\end{tabular}
}
\caption{Simulation results under the log-normal design without censoring.}
\label{tab:sim1}
\end{table}
    \end{landscape}
    \clearpage
}
    
  \afterpage{
\clearpage
    \thispagestyle{empty}
    \begin{landscape}
\begin{table}[!p]
\centering
\adjustbox{max width=1.3\textwidth}{
\begin{tabular}{l | c  c  c  c | c c c c | c c c c | c c c c}
\hline \hline
~ & $\hat{\theta}_{00}$ & $\hat{\theta}_{10}$  & $\hat{\theta}_{01}$ & $\hat{\theta}_{11}$ & $\hat{\theta}_{00}$ & $\hat{\theta}_{10}$  & $\hat{\theta}_{01}$ & $\hat{\theta}_{11}$ & $\hat{\theta}_{00}$ & $\hat{\theta}_{10}$  & $\hat{\theta}_{01}$ & $\hat{\theta}_{11}$ & $\hat{\theta}_{00}$ & $\hat{\theta}_{10}$  & $\hat{\theta}_{01}$ & $\hat{\theta}_{11}$ \\
\hline \hline
 \multicolumn{17}{c}{\textsc{Censoring}}\\ \hline
\multicolumn{17}{l}{\textbf{n = 500}}\\ \hline
~ & \multicolumn{4}{c}{$\alpha = 0.25, \beta = 1, \bar{D} = 0.514, \bar{\delta} = 0.796$} & \multicolumn{4}{c}{$\alpha = 0.75, \beta = 1, \bar{D} = 0.565, \bar{\delta} = 0.794$} & \multicolumn{4}{c}{$\alpha = 0.25, \beta = 0.5, \bar{D} = 0.467, \bar{\delta} = 0.791$} & \multicolumn{4}{c}{$\alpha = 0.75, \beta = 0.5, \bar{D} = 0.516, \bar{\delta} = 0.789$}\\
  \hline
Bias & -0.009 & -0.093 & 0.010 & 0.032 & -0.018 & -0.097 & 0.015 & 0.035 & -0.022 & -0.194 & -0.002 & 0.036 & -0.026 & -0.192 & 0.015 & 0.046 \\ 
  SE & 0.356 & 0.609 & 0.355 & 0.394 & 0.360 & 0.575 & 0.350 & 0.412 & 0.358 & 1.156 & 0.345 & 0.504 & 0.369 & 0.992 & 0.358 & 0.524 \\ 
  $90\%$ & 0.715 & 0.521 & 0.716 & 0.642 & 0.706 & 0.611 & 0.709 & 0.677 & 0.687 & 0.518 & 0.717 & 0.572 & 0.637 & 0.457 & 0.699 & 0.606 \\ 
  $95\%$ & 0.764 & 0.751 & 0.775 & 0.758 & 0.792 & 0.770 & 0.797 & 0.788 & 0.782 & 0.804 & 0.834 & 0.772 & 0.782 & 0.679 & 0.818 & 0.787 \\ 
  $99\%$ & 0.971 & 0.957 & 0.975 & 0.945 & 0.986 & 0.984 & 0.963 & 0.979 & 0.963 & 0.979 & 0.975 & 0.956 & 0.979 & 0.963 & 0.989 & 0.987 \\ 
   \hline

\multicolumn{17}{l}{\textbf{n = 1000}}\\ \hline
~ & \multicolumn{4}{c}{$\alpha = 0.25, \beta = 1, \bar{D} = 0.513, \bar{\delta} = 0.797$} & \multicolumn{4}{c}{$\alpha = 0.75, \beta = 1, \bar{D} = 0.565, \bar{\delta} = 0.795$} & \multicolumn{4}{c}{$\alpha = 0.25, \beta = 0.5, \bar{D} = 0.465, \bar{\delta} = 0.792$} & \multicolumn{4}{c}{$\alpha = 0.75, \beta = 0.5, \bar{D} = 0.515, \bar{\delta} = 0.79$}\\
  \hline
Bias & -0.003 & -0.034 & -0.000 & 0.005 & -0.018 & -0.021 & -0.004 & 0.029 & -0.015 & -0.063 & -0.001 & 0.024 & -0.020 & -0.098 & -0.007 & 0.033 \\ 
  SE & 0.367 & 0.484 & 0.347 & 0.382 & 0.354 & 0.444 & 0.355 & 0.358 & 0.375 & 0.821 & 0.367 & 0.411 & 0.364 & 0.814 & 0.349 & 0.418 \\ 
  $90\%$ & 0.707 & 0.607 & 0.715 & 0.689 & 0.725 & 0.588 & 0.708 & 0.715 & 0.701 & 0.464 & 0.700 & 0.583 & 0.692 & 0.431 & 0.703 & 0.638 \\ 
  $95\%$ & 0.779 & 0.841 & 0.824 & 0.800 & 0.776 & 0.754 & 0.794 & 0.822 & 0.766 & 0.722 & 0.763 & 0.727 & 0.768 & 0.643 & 0.787 & 0.777 \\ 
  $99\%$ & 0.941 & 0.951 & 0.961 & 0.936 & 0.977 & 0.962 & 0.967 & 0.991 & 0.977 & 0.985 & 0.924 & 0.910 & 0.965 & 0.984 & 0.952 & 0.993 \\ 
   \hline

\multicolumn{17}{l}{\textbf{n = 3000}}\\ \hline
~ & \multicolumn{4}{c}{$\alpha = 0.25, \beta = 1, \bar{D} = 0.513, \bar{\delta} = 0.796$} & \multicolumn{4}{c}{$\alpha = 0.75, \beta = 1, \bar{D} = 0.564, \bar{\delta} = 0.795$} & \multicolumn{4}{c}{$\alpha = 0.25, \beta = 0.5, \bar{D} = 0.466, \bar{\delta} = 0.791$} & \multicolumn{4}{c}{$\alpha = 0.75, \beta = 0.5, \bar{D} = 0.516, \bar{\delta} = 0.79$}\\
  \hline
Bias & -0.008 & -0.007 & -0.007 & 0.001 & -0.014 & -0.002 & 0.004 & 0.000 & -0.027 & -0.014 & 0.007 & 0.019 & -0.018 & -0.001 & 0.001 & 0.023 \\ 
  SE & 0.357 & 0.384 & 0.351 & 0.362 & 0.359 & 0.374 & 0.349 & 0.351 & 0.354 & 0.488 & 0.356 & 0.369 & 0.365 & 0.449 & 0.357 & 0.362 \\ 
  $90\%$ & 0.715 & 0.688 & 0.707 & 0.713 & 0.718 & 0.698 & 0.714 & 0.715 & 0.724 & 0.647 & 0.715 & 0.684 & 0.726 & 0.608 & 0.730 & 0.704 \\ 
  $95\%$ & 0.740 & 0.781 & 0.744 & 0.769 & 0.764 & 0.779 & 0.766 & 0.765 & 0.756 & 0.828 & 0.769 & 0.789 & 0.773 & 0.771 & 0.786 & 0.786 \\ 
  $99\%$ & 0.938 & 0.974 & 0.938 & 0.926 & 0.925 & 0.926 & 0.939 & 0.928 & 0.940 & 0.992 & 0.958 & 0.943 & 0.935 & 0.945 & 0.947 & 0.955 \\ 
   \hline
\end{tabular}
}
\caption{Simulation results under the log-normal design with censoring.}
\label{tab:sim2}
\end{table}
    \end{landscape}
    \clearpage
}
\section{Bootstrap confidence intervals for the empirical application} \label{S.4}

The 95\% confidence intervals are computed using 500 bootstrap draws. In each bootstrap sample, the optimization algorithm is started at the value of the estimator in the original sample. The confidence intervals are pointwise, in the sense that for all times equal to $1,2,\dots,770$ days, we compute the estimated values of the hazard rate (or difference of the hazard rates) in each bootstrap sample and then consider the $0.024$ and $0.976$ quantiles of the empirical distribution of these  estimated values of the hazard rate. It can be seen that the treatment significantly increases the hazard rate for all times below $400$ days with the Weibull model and all times below 770 days with the log-normal  model.

\begin{figure}[H]
\begin{minipage}{0.48\textwidth}
  \centering
  \includegraphics[width=80mm]{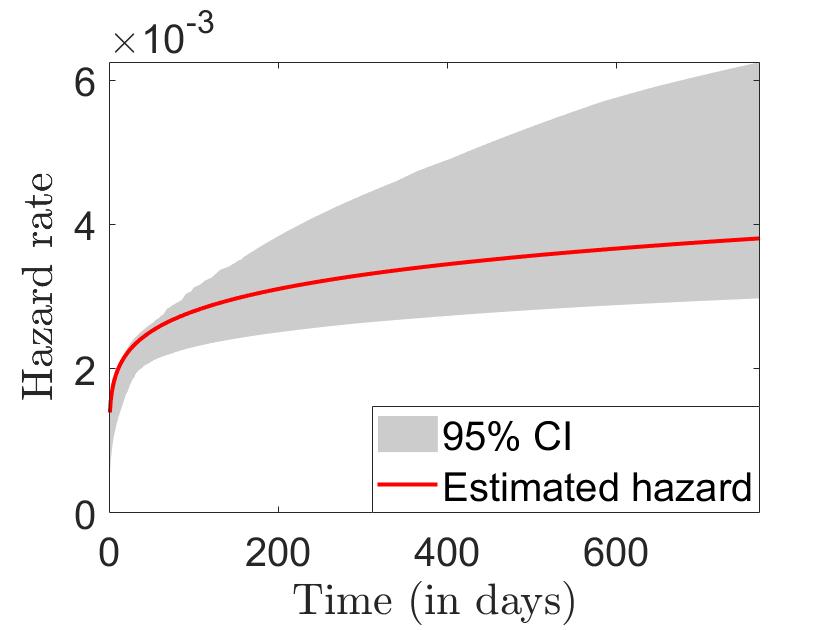} 
  \caption{Bootstrap 95\% confidence intervals for the hazard rate when subjects are never treated with the Weibull model.}
   \label{fig:cov_design1} 
    \end{minipage}
    \quad\quad
    \begin{minipage}{0.48\textwidth}
        \centering
   \includegraphics[width=80mm]{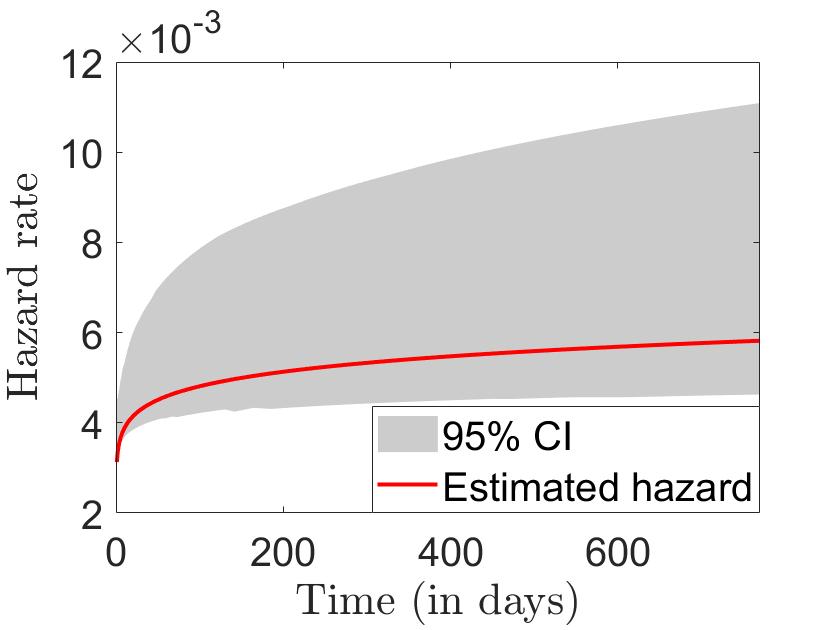}
  \caption{Bootstrap 95\% confidence intervals for the hazard rate when subjects are treated at time 0 with the Weibull model.}
      \label{fig:cov_design2}
      \end{minipage}
\end{figure}

\begin{figure}[H]
\begin{minipage}{0.48\textwidth}
  \centering
  \includegraphics[width=80mm]{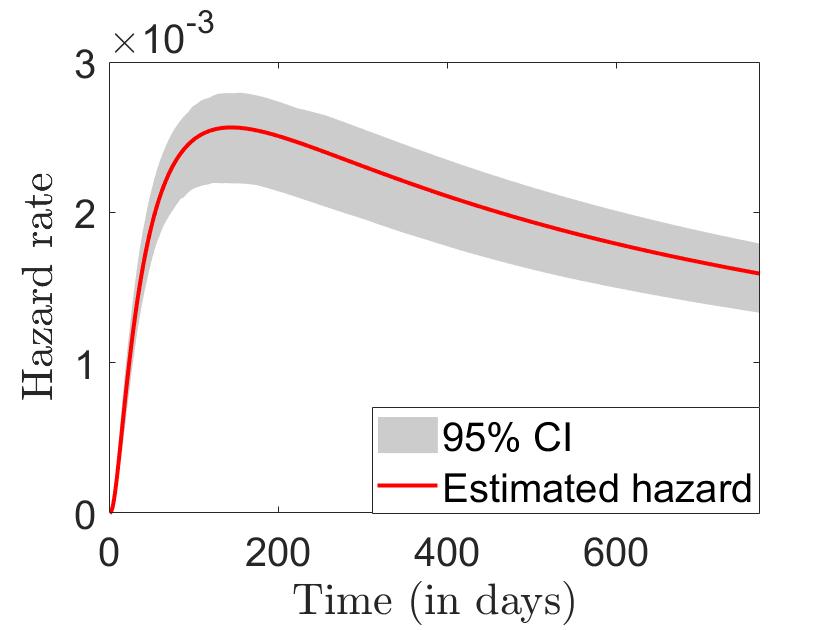}
  \caption{Bootstrap 95\% confidence intervals for the hazard rate when subjects are never treated with the log-normal model.}
   \label{fig:cov_design1} 
    \end{minipage}
    \quad\quad
    \begin{minipage}{0.48\textwidth}
        \centering
   \includegraphics[width=80mm]{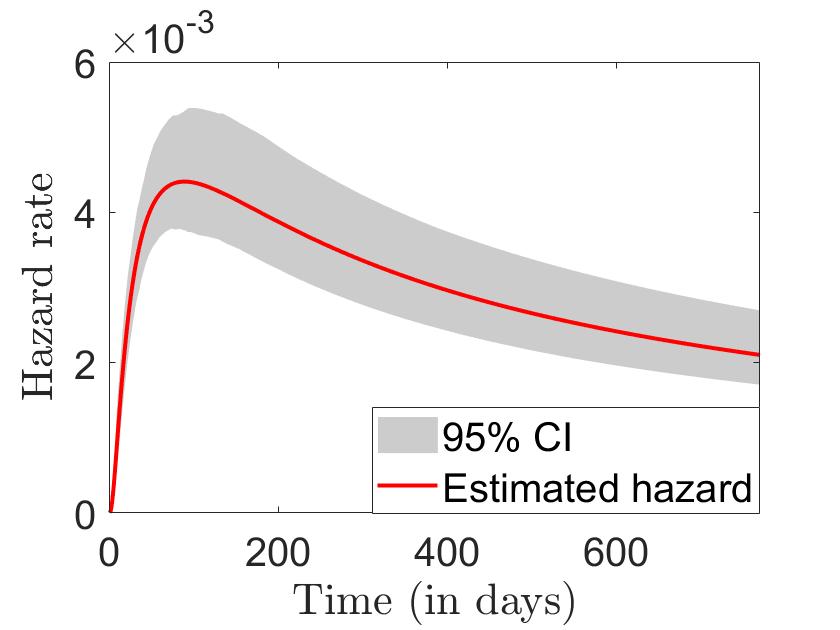}
  \caption{Bootstrap 95\% confidence intervals for the hazard rate when subjects are treated at time 0 with the log-normal model.}
      \label{fig:cov_design2}
      \end{minipage}
\end{figure}

\begin{figure}[H]
\begin{minipage}{0.48\textwidth}
  \centering
  \includegraphics[width=80mm]{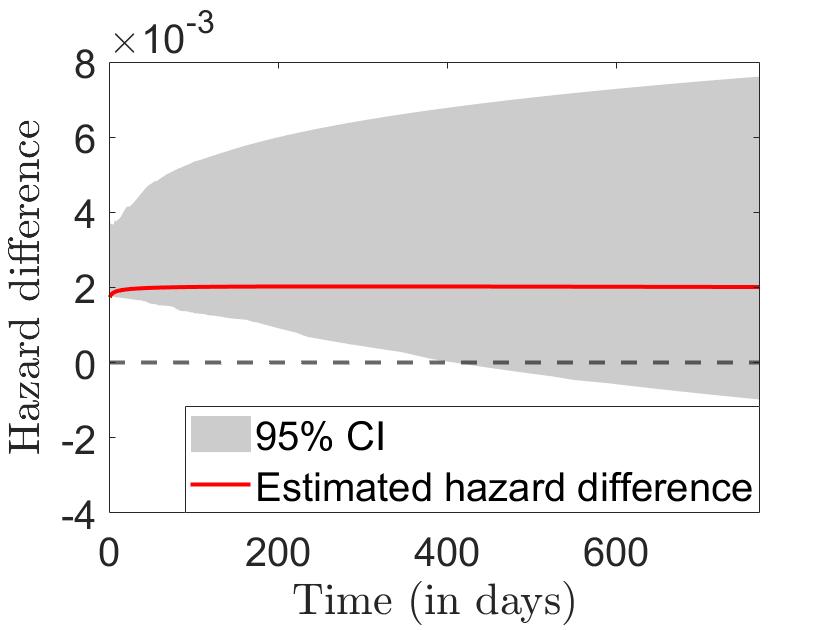}
  \caption{Bootstrap 95\% confidence intervals for the difference of the hazard rates of the treated at time 0 and the never treated with the Weibull model.}
   \label{fig:cov_design1} 
    \end{minipage}
    \quad\quad
    \begin{minipage}{0.48\textwidth}
        \centering
   \includegraphics[width=80mm]{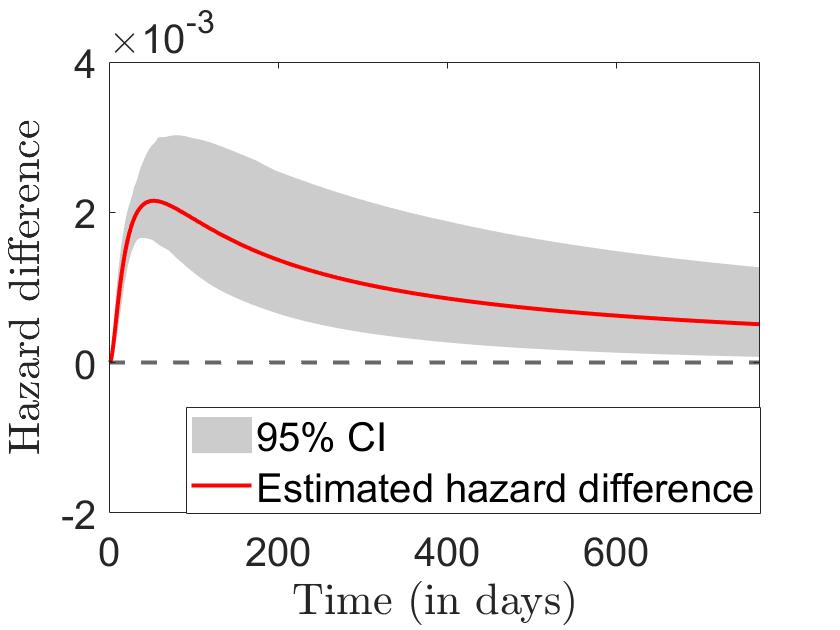}
  \caption{Bootstrap 95\% confidence intervals for the difference of the hazard rates of the treated at time 0 and the never treated with the log-normal model.}
      \label{fig:cov_design2}
      \end{minipage}
\end{figure}

\end{document}